\documentclass{article}
  \usepackage{amsfonts,amssymb,mathrsfs,pb-diagram,makeidx}
  \usepackage[all]{xy}
\def \1{{\bf 1}}
\def \a{{{\mathfrak a}}}

\def \al{\alpha}

\def \Ad{{\rm Ad}}

\def \b{{{\mathfrak b}}}

\def \bs{\backslash}

\def \bwedge{\left.\bigwedge\!\right.}

\def \C{{\mathbb C}}

\def \CC{{\mathscr C}}

\def \CD{{\mathscr D}}
\def \CE{{\mathscr E}}
\def \CEE{{\mathscr E}}

\def \CG{{\mathscr G}}

\def \CO{{\cal O}}
\def \cos{{\rm cos \hspace{2pt}}}
\def \coth{{\rm coth}}
\def \CP{{\mathscr P}}

\def \d{\delta}
\def \D{\Delta}
\def \det{{\rm det \hspace{2pt}}}
\def \diag{{\rm diag}}

\def \eac{{\acute{\rm e}}}
\def \eacit{{\acute{\it e}}}
\def \ell{{\rm ell}}
\def \ep{\varepsilon}
\def \exp{{\rm exp}}

\def \g{{{\mathfrak g}}}
\def \ga{\gamma}
\def \Ga{\Gamma}
\def \gen{{\rm gen}}

\def \Gr{{\rm Gr \hspace{2pt}}}
\def \h{{{\mathfrak h}}}

\def \Im{{\rm Im \hspace{2pt}}}

\def \Ind{{\rm Ind \hspace{2pt}}}

\def \k{{{\mathfrak k}}}

\def \la{\lambda}

\def \La{\Lambda}
\def \li{{\rm li}}
\def \Lie{{\rm Lie}}

\def \m{{{\mathfrak m}}}
\def \Mat{{\rm Mat}}
\def \mod{{\rm mod}}
\def \n{{{\mathfrak n}}}

\def \N{\mathbb N}
\def \nreg{\rm nreg}

\def \O{{\rm O}}
\def \om{\omega}
\def \Om{\Omega}
\def \ox{\otimes}
\def \p{{{\mathfrak p}}}

\def \prf{\noindent{\bf Proof: }}

\def \PGL{{\rm PGL}}

\def \qed{\ifmmode\eqno $\square$ \else\noproof\vskip 12pt plus 3pt minus 9pt \fi}
 \def\noproof{{\unskip\nobreak\hfill\penalty50\hskip2em\hbox{}%
     \nobreak\hfill $\square$ \parfillskip=0pt%
     \finalhyphendemerits=0\par}}
\def \Q{\mathbb Q}

\def \R{{\mathbb R}}
\def \Re{{\rm Re \hspace{1pt}}}

\def \ra{\rightarrow}
\def \Ra{\Rightarrow}

\def \reg{{\rm reg}}

\def \si{\sigma}

\def \sin{{\rm sin \hspace{2pt}}}
\def \SL{{\rm SL}}
\def \sl{{\rm sl}}
\def \SO{{\rm SO}}

\def \SS{{\rm S}}

\def \tanh{{\rm tanh}}
\def \Th{\Theta}
\def \th{\theta}
\def \tr{{\hspace{1pt}\rm tr\hspace{2pt}}}

\def \vol{{\rm vol}}

\def \x{\times}
\def \z{{\mathfrak z}}

\def \Z{\mathbb Z}
\def \={\ =\ }

\newcommand{\rez}[1]{\frac{1}{#1}}

\renewcommand{\matrix}[4]{\left( \begin{array}{cc}#1 & #2 \\ #3 & #4 \end{array}
            \right)}
\newcommand{\matrixfour}[4]
	   {\left( \begin{array}{cccc} 
	       #1 & \  & \  & \  \\ 
	       \  & #2 & \  & \  \\  
	       \  & \  & #3 & \  \\  
	       \  & \  & \  & #4 \end{array}\right)}
\newcommand{\matrixtwo}[2]{\matrix {#1}{}{}{#2}}
\newtheorem{theorem}{Theorem}[section]

\newtheorem{lemma}[theorem]{Lemma}

\newtheorem{proposition}[theorem]{Proposition}

\begin{document}

 \title{A prime geodesic theorem for $\SL_4$}
 \author{Anton Deitmar \& Mark Pavey}
 \date{}
 \maketitle
 
\tableofcontents

\newpage
\section*{Introduction}
The prime geodesic
theorem (PGT) gives a growth asymptotic for the number of closed geodesics counted by their lengths
\cite{Gangolli,Hejhal,Knieper,Koyama,LuoSarnak,Margulis,Pollicott,Zelditch}.
Classically, it was known for manifolds of strictly negative sectional curvature.
In \cite{primgeo}, one finds a generalization to higher rank locally symmetric manifolds, counting only regular geodesics.
In \cite{primgeoII} the PGT is extended to non-regular geodesics, under the proviso of a geometric regularity condition on the manifold.

The current paper is dedicated to quotients of the symmetric space of $\SL_4(\R)$.
In this situation the regularity condition fails.
As a consequence, the analytical difficulties grow exponentially.
The group $\SL_4(\R)$ is the ``biggest'' semisimple group for which the unitary dual is known.
This enables us to give prime geodesic theorems with remainder term.

In \cite{class,classnc,primgeo,primgeoII}, the prime geodesic theorem is applied to derive class number asymptotics.
The last open case for these is the case of complex quartic fields for which the PGT required is the one of $\SL_4(\R)$.
As there are additional algebraic obstacles, this application of the PGT is not straightforward and will be subject of a subsequent paper.

In what follows we give a summary of the arguments and techniques used in the proof of our main results.  In particular we shall point out the difficulties that arise in applying the methods of \cite{primgeo} to our situation.  Section \ref{ch:EulerChar} provides some necessary background and preliminary results.  In Section \ref{ch:Ruelle} we introduce the zeta functions we shall be studying and prove some of their analytic properties.  The results of Section \ref{ch:EulerChar} are made use of here in the definitions of the zeta functions and proofs of their analytic properties.  In Section \ref{ch:PGT} we apply standard techniques from analytic number theory to the results of Section \ref{ch:Ruelle} in order to prove the Prime Geodesic Theorem required in our context.  

Our strategy is that we define a suitable zeta function and use the trace formula with special test-functions to derive analytic continuation as in \cite{geom}, so that we can use the techniques of the proof of the prime number theorem.
It is in the analytic continuation that we encounter the first obstacle.  We say that an element $g\in G$ is \emph{weakly neat} if the adjoint $\Ad(g)$ has no non-trivial roots of unity as eigenvalues.  A subgroup of $G$ is weakly neat if every element is.  The results of \cite{geom} use the assumptions that the group $G$ has trivial centre and the group $\Ga$ is weakly neat.  Note that, if $G$ has trivial centre then $\Ga$ weakly neat implies $\Ga$ torsion free, since any non weakly neat torsion elements must be central.  To generalise the results of \cite{geom} to the case where $\Ga$ is not weakly neat requires two things.  Firstly, we need to modify the definition of the zeta function from that given in \cite{geom}.  Secondly, the definition of the zeta function includes the first higher Euler characteristic of the spaces $X_{\Ga_{\ga}}$ for $\ga\in\CE(\Ga)$.  In \cite{geom} these are only defined for $\Ga_{\ga}$ torsion free so we need to broaden the definition.

In Section \ref{ch:EulerChar} we give a definition of first higher Euler characteristics which is general enough for our application.  We further prove that its value is always positive in the cases we consider.  This fact is needed in the proof of the Prime Geodesic Theorem.  In \cite{geom} it is shown that the position of the poles and zeros of the generalised Selberg zeta function depend on the Lie algebra cohomology of the irreducible, unitary representations of $\SL_4(\R)$.  Using a result of Hecht and Schmid (\cite{HechtSchmid83}) it suffices to look at the infinitesimal characters of the irreducible, unitary representations of $\SL_4(\R)$.  For this purpose, in Section \ref{ch:EulerChar} we also describe the unitary dual of $\SL_4(\R)$ using a result of Speh (\cite{Speh81}) and give a result about the infinitesimal characters of certain elements of this set.

In Section \ref{ch:Ruelle} we define the generalised Selberg zeta function and use the trace formula to deduce its analytic properties.  In particular, we give a formula for its vanishing order at a given point and prove a functional equation, from which we can deduce that it is of finite order.  We then define the generalised Ruelle zeta function and prove that it is a finite quotient of generalised Selberg zeta functions.  In particular, in the cases we are interested in it has a zero at $s=1$ and all other poles and zeros in the half plane $\Re s\leq\frac{3}{4}$, and is furthermore of finite order.  We introduce the Ruelle zeta function as its logarithmic derivative is a Dirichlet series from whose properties we can prove the Prime Geodesic Theorem.

The definitions of the generalised Selberg and Ruelle zeta functions depend on the choice of a finite dimensional virtual representation of a particular closed subgroup $M$ of $\SL_4(\R)$.  The trace of this representation at elements of $\Ga$ appears in the Dirichlet series arising as the logarithmic derivative of the Ruelle zeta function.  We can choose such a representation $\tilde{\si}$ so that its trace is zero for all non-regular elements of $\Ga$.  

In Section \ref{ch:PGT} we apply the methods of \cite{Randol77} together with standard techniques of analytic number theory to prove a Prime Geodesic Theorem.  In carrying out the application to class numbers it is necessary to show that certain subsets of the summands contribute negligibly to the asymptotic.  This is done using a ``sandwiching" argument.  To carry out this ``sandwiching" argument we define a sequence of Dirichlet series, which are not connected to any zeta function, but whose analytic properties can also be deduced from the Selberg trace formula.  We use a version of the Wiener-Ikehara theorem to prove an asymptotic result for an increasing sequence of functions derived from these Dirichlet series, which we use to ``sandwich" the product over the elements  we are interested in against the sum over all elements which comes from the Ruelle zeta function.  Unfortunately the asymptotics we are able to derive from the Wiener-Ikehara theorem do not provide an error term like those we can deduce from the Ruelle zeta function using the methods of \cite{Randol77}.  This is why we lose the error term in Theorem~\ref{thm:PGT2}.

We finish this introduction with a few remarks about the limitations of the method used here in terms of further applications and mention a couple of other recent results in the same direction.  In order to be able to make use of the trace formula for compact spaces we have had to limit our sum over class numbers by means of the choice of a finite set of primes, as described above.  In \cite{classnc} Deitmar and Hoffmann have been able to use a different trace formula to prove that as $x\ra\infty$
$$
\sum_{R(\CO)\leq x}h(\CO)\sim\frac{e^{3x}}{3x},
$$
where the sum is over all isomorphism classes of orders in complex cubic fields.

In order to get the error term in the Prime Geodesic Theorem we have made use of the classification of the unitary dual of $\SL_4(\R)$.  At present the unitary dual is not known for any higher dimensional groups so a Prime Geodesic Theorem with error term is not possible using our methods.  Finally we mention that the correspondence between geodesics and orders actually works by identifying primitive geodesics with fundamental units in orders.  By Dirichlet's unit theorem, an order in a number field $F$ has a unique fundamental unit (up to inversion and multiplication by a root of unity) only if $F$ is real quadratic, complex cubic or totally complex quartic.  Hence an asymptotic of our form can be proven only in these three cases.

\newpage
\section{Euler Characteristics and Infinitesimal Characters}
    \label{ch:EulerChar}

In this section we introduce some concepts and prove some results which will be needed for our consideration of the zeta functions in the next section.

\subsection{The unitary duals of $\SL_2(\R)$ and $\SL_4(\R)$}

The \emph{unitary dual} \index{unitary dual} $\hat{G}$ \index{$\hat{G}$} of a locally compact group $G$ is the set of all equivalence classes of irreducible unitary representations of $G$.  If $\pi$ is an irreducible, unitary representation of $G$ then we shall write, by slight abuse of notation, $\pi\in\hat{G}$.
Let $G_1=\SL_2(\R)$.  We define the following subgroups: let $K_1$ be the maximal compact subgroup $\SO(2)$; let $M_1=\{\pm I_2\}$.
let $A_1$ be the group of diagonal matrices $\diag(a,a^{-1}$ where $a>0$, 
and let $N_1$ be the group of upper triangular matrices in $G_1$ with ones on the diagonal.
Let $\g_1$, $\k_1$, $\a_1$, $\n_1$ be the respective complexified Lie algebras.  Let $\rho_1\in\a_1^*$ be defined by $\rho(\diag(a,-a))=a$,
and let $P_1=M_1 A_1 N_1$ be the parabolic subgroup of $G_1$ with split torus $A_1$ and unipotent radical $N_1$.
We will use the notation of \cite{Knapp}, Chap. II.
For $n\geq 2$ we have the discrete series representations \index{representation!discrete series} $\CD_n^+$ \index{$\CD_n^\pm$} and $\CD_n^-$ of $G_1$.  
We also denote the two limit of discrete series representations \index{representation!limit of discrete series} by $\CD_1^+$ and $\CD_1^-$.

Let $\tau^+, \tau^-$ be the two characters of $M-1$.
For $\nu\in\a_1^*$ we define
$
\CP^{\pm,\nu}=\Ind_{P_1}^{G_1}(\tau^{\pm}\ox\nu).
$
For $x\in\R$ define the \emph{principal series representations}\index{representation!principal series}
$$
\CP^{\pm,ix}=\CP^{\pm,ix\rho_1}.
$$
Then the representations $\CP^{\pm,ix}$ \index{$\CP^{\pm,ix}$} are unitary and are all irreducible except for $\CP^{-,0}\cong\CD_1^+\oplus\CD_1^-$.  The only equivalences among the representations $\CP^{\pm,\nu}$ are that $\CP^{+,ix}$ and $\CP^{-,ix}$ are equivalent to $\CP^{+,-ix}$ and $\CP^{-,-ix}$ respectively for all $x\in\R$.
We also define the \emph{complementary series representations}\index{representation!complementary series}
$
\CC^x=\CP^{+,x\rho_1},
$
for $0<x<1$.  
We have the following classification theorem:

\begin{theorem}
\label{thm:SL2Dual}
\index{unitary dual!of $\SL_2(\R)$}
The unitary dual of $\SL_2(\R)$ consists of

a) the trivial representation;

b) the principal series representations $\CP^{+,ix}$ for $x\in\R$ and $\CP^{-,ix}$ for $x\in\R\smallsetminus\{0\}$;

c) the complementary series $\CC^x$ for $0<x<1$;

d) the discrete series $\CD_n^+$ and $\CD_n^-$ for $n\geq 2$ and limits of discrete series $\CD_1^+$ and $\CD_1^-$.

The only equivalences among these representations are $\CP^{\pm,ix}\cong\CP^{\pm,-ix}$ for all $x\in\R$.

\end{theorem}
\prf
\cite{Knapp}, Theorem 16.3.
\qed

Now let $G=\SL_4(\R)$ and $K=\SO(4)$, so $K$ is a maximal compact subgroup of $G$.  Let $P'=M'A'N'$ be a parabolic subgroup of $G$ with split component $A'$ and unipotent radical $N'$.  Let $\g$ and $\a'$ be the complexified Lie algebras of $G$ and $A'$ and let $\a'^*$ be the complex dual of $\a'$.  Let $\rho'$ be the half sum of the positive roots of the system $(\g',\a')$.  Then we can define induced representations in an entirely analogous way to that used for for $G_1=\SL_2(\R)$ above.  For an irreducible, unitary representation $\tau$ of $M'$ and for $\nu\in\a'^*$ we write the corresponding induced representation of $G$ as $\Ind_{P'}^G(\tau\ox\nu)$.
If $\nu=ix\rho'$ for some $x\in\R$ then the induced representation $\Ind_{P'}^G(\tau\ox\nu)$ is unitary with respect to the inner product
$
\langle f,g\rangle=\int_{K} \langle f(k),g(k)\rangle_{\tau}\ dk.
$
In this case we call $\Ind_{P'}^G(\tau\ox\nu)$ a principal series representation.  For other choices is $\nu\in\a'^*$ it may be possible to make $\Ind_{P'}^G(\tau\ox\nu)$ unitary with respect to a different inner product, in which case we call $\Ind_{P'}^G(\tau\ox\nu)$ a complementary series representation.  There are also certain irreducible, unitary subrepresentations of induced representations (see \cite{Speh81}, p121) which are called \emph{limit of complementary series representations}\index{representation!limit of complementary series}.  These limit of complementary series representations can often be realised as induced representations from a parabolic subgroup $P''\supset P'$.

We define the following subgroups of $G$.  Let
\begin{eqnarray*}
M & = & \SS\matrixtwo{\SL_2^{\pm}(\R)}{\SL_2^{\pm}(\R)} \\
  & \cong & \left\{(x,y)\in \Mat_2(\R)\times \Mat_2(\R)|\begin{array}{c}\det(x),\det(y)=\pm 1 \\ \det(x)\det(y)=1 \end{array}\right\}.
\index{$M$}
\end{eqnarray*}
Let $A$ be the group of all diagonal matrices $\diag(a,a,a^{-1},a^{-1})$, where $a>0$. Finally, let
$
N=\matrix{I_2}{\Mat_2(\R)}{0}{I_2}.\index{$N$}
$
Let $P$ be the parabolic subgroup of $G$ with Langlands decomposition $P=MAN$.  For $m_1, m_2\in\N$, we denote by $\bar{\pi}_{m_1,m_2}$ the representation of $M$ induced from the representation $\CD_{m_1}^+\ox\CD_{m_2}^+$ of $\SL_2(\R)\x\SL_2(\R)$.  For $m\in\N$, let $\bar{\pi}_m=\bar{\pi}_{m,m}$ and let $I_m =\Ind_P^G \left(\bar{\pi}_m \otimes \frac{1}{2} \rho_P \right)$.  The representations $I_m$ each have a unique irreducible quotient known as the Langlands quotient (see \cite{Knapp}, Theorem 7.24).  We denote the Langlands quotient of $I_m$ by $\pi_m$.\index{$\pi_m$}

We have the following classification theorem:

\begin{theorem}
\label{thm:SL4Dual}
\index{unitary dual!of $\SL_4(\R)$}
The unitary dual of $\SL_4(\R)$ consists of

a) the trivial representation;

b) principal series representations;

c) complementary series representations $\Ind^G_P\left(\bar{\pi}_m\ox t\rho_P\right)$, for $m\in\N$ and $0<t<\rez{2}$;

d) complementary series representations induced from parabolics other than $P=MAN$;

e) limit of complementary series representations, which are irreducible unitary subrepresentations of $I_m$, for $m\in\N$;

f) the family of representations $\pi_m$, indexed by $m\in\N$.
\end{theorem}
\prf
This follows from \cite{Speh81}, Theorem 5.1, where the unitary dual of $\SL_4^{\pm}(\R)=\{X\in\Mat_4(\R):\det X=\pm 1\}$ is given.
\qed

\subsection{Euler characteristics}
\label{sec:EulerChars}
Let $G$ be a real reductive group and suppose that there is a finite subgroup $E$ of the centre of $G$ and a reductive and Zariski-connected linear group $\CG$ such that $G/E$ is isomorphic to a subgroup of $\CG(\R)$ of finite index.  Note these conditions are satisfied whenever $G$ is a Levi component of a connected semisimple group with finite centre.    Let $K$ be a maximal compact subgroup of $G$ and let $\Ga$ be a discrete, cocompact subgroup of $G$.  Let $X_{\Ga}$ be the locally symmetric space $\Ga\bs G/K$.  If $\Ga$ is torsion-free we define the \emph{first higher Euler characteristic} \index{first higher Euler characteristic}\index{Euler characteristic!first higher} of $\Ga$ to be
$$
\chi_1(X_{\Ga}) = \chi_1(\Ga) = \sum_{j=0}^{\dim X_{\Ga}} (-1)^{j+1}jh^j(X_{\Ga}),\index{$\chi_1(\Ga)$}
$$
where $h^j(X_{\Ga})$ is the $j$th Betti number of $X_{\Ga}$.  We want to generalise this definition to cases where the group $\Ga$ is not necessarily torsion free.  We prove below a proposition which allows us broaden the definition to all cases required by our applications.

Let $\th$ be the Cartan involution fixing $K$ pointwise.  There exists a $\th$-stable Cartan subgroup $H=AB$ of $G$, where $A$ is a connected split torus and $B\subset K$ is a Cartan of $K$.  We assume that $A$ is central in $G$.  Let $C$ denote the centre of $G$.  Then $C\subset H$ and we write $B_C$, $\Ga_C$ for $B\cap C$ and $\Ga\cap C$ respectively.  Let $G^1$ be the derived group of $G$ and let $\Ga^1=G^1\cap\Ga C$.  We note in particular that, since $G=G^1C$, we have $\Ga\subset\Ga^1C\subset\Ga^1AB$.  Let $\Ga_A=A\cap\Ga_C B_C$ \index{$\Ga_A$} be the projection of $\Ga_C$ to $A$.  Then $\Ga_A$ is discrete and cocompact in $A$ (see \cite{Wolf62}, Lemma 3.3).

Let $\g_{\R}$ be the real Lie algebra Lie($G$) with polar decomposition $\g_{\R}=\k_{\R}\oplus\p_{\R}$.  Let $b$ be a fixed nondegenerate invariant bilinear form on the Lie algebra $\g_{\R}$ which is negative definite on $\k_{\R}$ and positive definite on $\p_{\R}$.  Then the form $-b(X,\th(Y))$ is positive definite and thus defines a left invariant metric on $G$.  Unless otherwise stated all Haar measures will be normalised according to the Harish-Chandra normalisation given in \cite{HarishChandra75}, \S 7.  Note that this normalisation depends on the choice of the bilinear form $b$ on $\g$.
On the space $G/H$ we have a pseudo-Riemannian structure given by the form $b$.  The Gauss-Bonnet construction generalises to pseudo-Riemannian structures to give an Euler-Poincar$\eac$ measure $\eta$ on $G/H$.  Define a (signed) Haar measure on $G$ by
$$
\mu_{EP}=\eta\ox(\textnormal{Haar measure on }H).
$$
The \emph{Weyl group} \index{Weyl group} $W=W(G,H)$ is defined to be the quotient of the normaliser of $H$ in $G$ by the centraliser.  It is a finite group generated by elements of order two.  We define the \emph{generic Euler characteristic} \index{generic Euler characteristic}\index{Euler characteristic!generic} by
$$
\chi_{\gen}(\Ga)=\chi_{\gen}(X_{\Ga})=\frac{\mu_{EP}(\Ga\bs G)}{|W|}.\index{$\chi_{\gen}(\Ga)$}
$$

\begin{proposition}
\label{pro:ECharWellDef}
Assume $\Ga$ is torsion free, $A$ is central in $G$ of dimension one and $\Ga'\subset\Ga$ is of finite index in $\Ga$.  Then $A/\Ga_A$ acts freely on $X_{\Ga}$ and $\chi_{\gen}(X_{\Ga})=\chi_1(\Ga)\vol(A/\Ga_A)$.  It follows that
$$
\chi_1(\Ga)=\chi_1(\Ga')\frac{\left[\Ga_A:\Ga'_A\right]}{\left[\Ga:\Ga'\right]}.
$$
\end{proposition}
\prf
The group $A_{\Ga}=A/\Ga_A$ acts on $\Ga\bs G/B$ by multiplication from the right.  We claim that this action is free, i.e., that it defines a fibre bundle
$$
A_{\Ga}\ra\Ga\bs G/B\ra\Ga\bs G/H.
$$
To see this let $\Ga xaB=\Ga xB$ for some $a\in A$ and $x\in G$.  Then $a=x^{-1}\ga xb$ for some $\ga\in \Ga$ and $b\in B$.  Since $\Ga\subset\Ga^1C\subset\Ga^1AB$ we can write $\ga$ as $\ga^1 a_{\ga}b_{\ga}$, with $\ga^1\in\Ga^1$ and $a_{\ga}\in A$ and $b_{\ga}\in B_C=B\cap C$.  It follows that $a_{\ga}\in\Ga_A$.  Since $A$ is central in $G$, we can write $\ga^1=aa_{\ga}^{-1}xb^{-1}x^{-1}b_{\ga}^{-1}$.  Since $\ga^1\in G^1$ and $aa_{\ga}^{-1}\in A\subset C$, we must have $aa_{\ga}^{-1}=1$, so $a=a_{\ga}\in \Ga_{A}$, which implies the claim.
In the same way we see that we get a fibre bundle
\begin{equation}
\label{eqn:TBundle}
A_{\Ga}\ra\Ga\bs G/K\ra A\Ga\bs G/K.
\end{equation}

Let $\chi$ be the usual Euler characteristic.  From \cite{HopfSamelson41} we take the equation
$\chi(K/B)=|W|$.  Using multiplicativity of Euler numbers in the fibre bundle
$
K/B\ra A\Ga\bs G/B\ra A\Ga\bs G/K
$
we get
\begin{eqnarray*}
\chi_{\gen}(X_{\Ga}) & = & \vol(A/\Ga_A)\frac{\eta(\Ga\bs G/H)}{|W|} \\
  &   & \\
  & = & \vol(A/\Ga_A)\frac{\chi(\Ga\bs G/H)}{\chi(K/B)} \\
  &   & \\
  & = & \vol(A/\Ga_A)\frac{\chi(A\Ga\bs G/B)}{\chi(K/B)} \\
  &   & \\
  & = & \vol(A/\Ga_A)\chi(A\Ga\bs G/K) \\
  & = & \vol(A/\Ga_A)\chi(A\bs X_{\Ga}).
\end{eqnarray*}
It remains to show that $\chi(A\bs X_{\Ga})=\chi_1(\Ga)$.

Let $\a_{\R}$ and $\g^1_{\R}$ be the Lie algebras of $A$ and $G^1$ respectively.  Then we can write
\begin{equation}
\label{eqn:gDecomp}
\g_{\R}=\a_{\R}\oplus\g^1_{\R}\oplus\z_{\R},
\end{equation}
where $\z_{\R}$ is central in $\g_{\R}$.  Let $X$ be the bi-invariant vector field on $\Ga\bs G/K$ generating the $A_{\Ga}$ action.  Then we can consider $X$ as an element of $\a_{\R}$ under the decomposition (\ref{eqn:gDecomp}).  Considering the dual of (\ref{eqn:gDecomp}), we can identify $\a_{\R}^*$ with a subspace of $\g_{\R}^*$.  Since $\a$ is central in $\g$, a non-zero element of $\a_{\R}^*$ gives us an $A_{\Ga}$-invariant, closed 1-form $\om$ on $\Ga\bs G/K$ such that for every $\Ga gK\in\Ga\bs G/K$ we have $\om(\Ga gK)(X)\neq 0$.  Since $A_{\Ga}\cong\R/\Z$ is connected and compact, the cohomology of the de Rham complex $\Ga\bs G/K$ coincides with the cohomology of the subcomplex of $A_{\Ga}$-invariants $\Om(\Ga\bs G/K)^{A_{\Ga}}$.  Using local triviality of the bundle one sees that
$$
\Om(\Ga\bs G/K)^{A_{\Ga}}=\pi^*\Om(A\Ga\bs G/K)\oplus\pi^*\Om(A\Ga\bs G/K)\wedge\om,
$$
where $\pi^*$ denotes the projection map.  Set $C_0=\pi^*\Om(A\Ga\bs G/K)$ and $C_1=C_0\oplus C_0\wedge\om=\Om(\Ga\bs G/K)^{A_{\Ga}}$.  Then
$$
H^p(C_1)=H^p(C_0)\oplus H^{p-1}(C_0)
$$
and so
\begin{eqnarray}
\label{eqn:ECharEq}
\chi_1(C_1) & = & \sum_{p\geq 0} (-1)^{p+1} p\dim H^p(C_1) \nonumber \\
  & = & \sum_{p\geq 0} (-1)^{p+1} p\left(\dim H^p(C_0)+\dim H^{p-1}(C_0)\right) \nonumber \\
  & = & \sum_{p\geq 0} (-1)^{p+1} (p-(p+1))\dim H^p(C_0) \nonumber \\
  & = & \sum_{p\geq 0} (-1)^p \dim H^p(C_0) \nonumber \=  \chi(C_0).
\end{eqnarray}
This gives the required result.
\qed

Let $\Ga$ be a discrete, cocompact subgroup of $G$.  
Then $\Ga$ has a torsion-free subgroup $\Ga'$ of finite index. Suppose that the torus $A$ is central in $G$ and of dimension one.  Define $\Ga_A$ and $\Ga'_A$ as above.

We define the first higher Euler characteristic \index{first higher Euler characteristic}\index{Euler characteristic!first higher} of $\Ga$ as
\begin{equation}
\label{eqn:fhEC}
\chi_1(X_{\Ga})=\chi_1(\Ga)=\chi_1(\Ga')\frac{\left[\Ga_A:\Ga'_A\right]}{\left[\Ga:\Ga'\right]}.\index{$\chi_1(\Ga)$}
\end{equation}
Propostition \ref{pro:ECharWellDef} shows that this is well-defined.  We note that in the case that $\Ga$ itself is torsion free, this definition of first higher Euler characteristic agrees with the one given above.

\subsection{Euler-Poincar\'e functions}
For the next definition we assume $G$ to be a semisimple real reductive group of inner type and we fix a maximal compact subgroup $K$.  We further assume that $G$ contains a compact Cartan subgroup.  Let $\g_{\R}$ be the real Lie algebra Lie($G$) with polar decomposition $\g_{\R}=\k_{\R}\oplus\p_{\R}$, and write $\g=\k\oplus\p$ for its complexification.  Recall that we are using Harish-Chandra's Haar measure normalisation as given in \cite{HarishChandra75}, \S 7 and this normalisation depends on the choice of an invariant bilinear form $b$ on $\g$, which for our purposes in this section we shall leave arbitrary.

A \emph{virtual representation} \index{representation!virtual} $\si$ of a group is a formal difference of two representations $\si=\si^+-\si^-$, which is called finite dimensional if both $\si^+$ and $\si^-$ are.  Two virtual representations $\si=\si^+-\si^-$ and $\tau=\tau^+-\tau^-$ of a group are said to be isomorphic if there is an isomorphism $\si^+\oplus\tau^-\cong\tau^+\oplus\si^-$.  The trace and determinant of a virtual representation $\si=\si^+-\si^-$ are defined by $\tr\si=\tr\si^+-\tr\si^-$ and $\det\si=\det\si^+/\det\si^-$.  The dimension of $\si$ is defined as $\dim\si=\dim\si^+-\dim\si^-$.

If $V$ is a representation space with $\Z$-grading then we shall consider it naturally as a virtual representation space by $V^+=V_{{\rm even}}$ and $V^-=V_{{\rm odd}}$.  In particular, if $V$ is a subspace of $\g$ we shall always consider the exterior product $\bwedge^*V$ as a virtual representation $\bwedge^*V=\bwedge^{{\rm even}}V-\bwedge^{{\rm odd}}V$ with respect to the adjoint representation.  We consider symmetric powers and cohomology spaces similarly.

For a smooth function $f$ on $G$ of compact support and an irreducible unitary representation $(\pi,V_{\pi})\in\hat{G}$ define the operator
$$
\pi(f)=\int_G \pi(g)f(g) dg\index{$\pi(f)$}
$$
on $V_{\pi}$.  Since $f$ is smooth and has compact support, $\pi(f)$ is of trace class.

Let $(\si,V_{\si})$ be a finite dimensional virtual representation of $G$.  An \emph{Euler-Poincar$\eacit$ function} \index{Euler-Poincar$\eac$ function} $f_{\si}$ \index{$f_{\si}$} for $\si$ is a compactly supported, smooth function on $G$ such that $f_{\si}\left(kxk^{-1}\right)=f_{\si}\left(x\right)$ for all $k\in K$ and for every irreducible unitary representation $(\pi,V_{\pi})$ of $G$ the following identity holds:
\begin{equation}
\label{eqn:EPFnId}
\tr \pi\left(f_{\si}\right)=\sum_{p=0}^{\dim(\p)} (-1)^p \dim\left( V_{\pi}\ox\bwedge^p\p\ox V_{\si}\right)^K,
\end{equation}
where the superscript $K$ denotes the subspace of $K$ invariants.  By \cite{Labesse} such functions exist.  We note that the value of such functions depends on the choice of Haar measure.  We shall assume all Euler-Poincar$\eac$ functions are given with respect to the Harish-Chandra normalisation.  In the following lemmas we prove some of their properties.

\begin{lemma}
\label{lem:EPfnrestrict}
Let $G$ denote a semisimple real reductive group of inner type, with connected component $G^0$, maximal compact subgroup $K$ and compact Cartan subgroup $T\subset K$.  Let $G^+ =TG^0$.  Further let $\si$ be a finite dimensional representation of $G$, $\si^+ =\si|_{G^+}$ and $f_{\si}$ an Euler-Poincar$\eacit$ function for $\si$ on $G$.

Then $f_{\si}|_{G^+}$ is an Euler-Poincar$\eacit$ function for $\si^+$ on $G^+$.
\end{lemma}

\prf
This is Lemma 1.5 of \cite{geom}.
\qed

\begin{lemma}
\label{lem:EPfnproduct}
Let $H$, $H_1$, $H_2$ be real reductive groups of inner type such that $H=H_1\x H_2$.  Let $\si$ be an irreducible representation of $H$.  There exist irreducible representations $\si_1, \si_2$ of $H_1, H_2$ respectively such that $\si\cong \si_1\ox\si_2$ and let $f_{\si_i}$ be an Euler-Poincar$\eacit$ function for $\si_i$ on $H_i$.

Then $f_{\si}(h_1,h_2)=f_{\si_1}(h_1)f_{\si_2}(h_2)$ is an Euler-Poincar$\eacit$ function for $\si$ on $H$.
\end{lemma}
\prf
This is a straightforward computation.
\qed

\begin{lemma}
\label{lem:EPCentral}
Let $G$ denote a semisimple real reductive group of inner type, with maximal compact subgroup $K$ and compact Cartan subgroup $T\subset K$.  Let $g\in G$ be central.  Let $\si$ be a finite dimensional representation of $G$ and let $f_{\si},h_{\si}$ be Euler-Poincar$\eacit$ functions for $\si$ on $G$.  Then $f_{\si}(g)=h_{\si}(g)$. 
\end{lemma}

\prf
Since $g$ is central, the orbital integral can be written as
$$
\CO_g(f_{\si})=\int_{G_g\bs G}f_{\si}(x^{-1}gx)\,dx=f_{\si}(g),
$$
where $G_g$ denotes the centraliser of $g$ in $G$ (which in this case, since $g$ is central, is the whole of $G$).  Similarly $\CO_g(h_{\si})=h_{\si}(g)$.  By \cite{geom}, Proposition 1.4, the value of the orbital integral $\CO_g(f_{\si})$ of an Euler-Poincar$\eac$ function for $\si$ on $G$ depends only on $g$, not on the particular Euler-Poincar$\eac$ function chosen.  Hence, $f_{\si}(g)=\CO_g(f_{\si})=\CO_g(h_{\si})=h_{\si}(g)$, as claimed.
\qed

\begin{lemma}
\label{lem:EPfnvalue}
Let $g_1$ be an Euler-Poincar$\eacit$ function for the trivial representation on $\SL_2(\R)$.  Then $g_1(1)=g_1(-1)\in\R$.
\end{lemma}
\prf
For $v\in\R$ let $\CP^{+,iv}$ and $\CP^{-,iv}$ be principal series representations on $\SL_2(\R)$.  For $n\in\N, n\geq 2$ let $\CD^+_n$ and $\CD^-_n$ be discrete series representations on $\SL_2(\R)$ and write $\CD^{\pm}_n$ for $\CD^+_n\oplus\CD^-_n$.  By \cite{Knapp}, Theorem 11.6 there is a constant $M>0$ such that for any compactly supported, smooth function $f$ on $\SL_2(\R)$ we have
\begin{eqnarray*}
f(1) & = & M\sum_{n=2}^{\infty}(n-1)\tr\CD^{\pm}_n(f) \\
  & & + \frac{M}{4}\int_{-\infty}^{\infty}\tr\CP^{+,iv}(f)v\tanh\left(\frac{\pi v}{2}\right)+\tr\CP^{-,iv}(f)v\coth\left(\frac{\pi v}{2}\right)dv.
\end{eqnarray*}
Lemma 1.3 of \cite{geom} tells us that $\tr\CP^{\pm,iv}(g_1)=0$ for all $v\in\R$, so we have
\begin{equation}
\label{eqn:EPSL2}
g_1(1) = M\sum_{n=2}^{\infty}(n-1)\tr\CD^{\pm}_n(g_1).
\end{equation}
By the definition of an Euler-Poincar$\eac$ function, for all $n\geq 2$
$$
\tr\CD^{\pm}_n(g_1)=\sum_{p=0}^2(-1)^p\dim\left(\CD^{\pm}_n\ox\bwedge^p\p\right)^{\SO(2)},
$$
where $\p$ is the complex Lie algebra
$$
\p=\left\{\matrix{a}{b}{b}{-a}:a,b\in\C\right\}.
$$
Hence $g_1(1)\in\R$.

We want to know for which values of $n$ the trace $\tr\CD^{\pm}_n(g_1)$ is non-zero.  For this we need to know the $\SO(2)$-types of $\CD^{\pm}_n$ and $\p$.

For $l\in\Z$, define the one dimensional representation $\ep_l$ of $\SO(2)$ by
$$
\ep_l R(\th)=e^{il\th},
\quad {\rm where}\quad
R(\th)=\matrix{\cos\th}{-\sin\th}{\sin\th}{\cos\th}\in \SO(2).
$$
Note that $\ep_0$ is the trivial representation, which we shall denote by $triv$.

The unitary dual of $\SO(2)$ is the set $\{\ep_l:l\in\Z\}$.

We have the following isomorphisms of $\SO(2)$-modules:
$$
\bwedge^0 \p \=  triv, \quad
\bwedge^1 \p \=  \ep_2\oplus\ep_{-2}, \quad
\bwedge^2 \p \=  triv.
$$

\begin{lemma}
\label{lem:DiscSer}
For $n\in\N$ we have an isomorphism of $\SO(2)$-modules:
$$
\CD^{\pm}_n \cong \bigoplus_{{|j|\geq n}\atop{j\equiv n\,\mod 2}}\ep_j.
$$
\end{lemma}
\prf
Let $\tau_n$ be the unique $n$-dimensional representation of $\SL_2(\R)$ (see \cite{Knapp}, Chapter II $\S$1).  It follows from the definition of $\tau_n$ that the following isomorphism of $\SO(2)$-modules holds:
$$
\tau_{n} \cong \bigoplus_{{|j|\leq n-1}\atop{j\equiv n-1\ \mod(2)}}\ep_j.
$$

Let $P_1=M_1 A_1 N_1$ be the minimal parabolic of $\SL_2(\R)$.  Then the unitary dual of $M_1\cong \{1,-1\}$ consists of two one dimensional representations, which we denote by $1=triv$ and $-1$.  We denote by $\rho_1$ the character of $A_1$
$$
\rho_1\matrixtwo{a}{a^{-1}}=a,
$$
and write $\pi_{\pm 1,n-1}$ for $\Ind_{\bar{P}}\pm 1\otimes (n-1)\rho_{\bar{P}}$.  Then we have the following exact sequences of $\SL_2(\R)$-modules (see \cite{Knapp}, Chapter II \S5):
$$
0\rightarrow\CD^{\pm}_n\rightarrow\pi_{1,n-1}\rightarrow\tau_{n-1}\rightarrow 0,\ \ \ \ n\in\N,\ n\ {\rm even}
$$
and
$$
0\rightarrow\CD^{\pm}_n\rightarrow\pi_{-1,n-1}\rightarrow\tau_{n-1}\rightarrow 0,\ \ \ \ n\in\N,\ n\ {\rm odd}, n\neq 1
$$
and the isomorphism of $\SL_2(\R)$-modules
$$
\CD_1^{\pm}\cong\pi_{-1,0}.
$$
From the so-called \emph{compact picture} of induced representations given in \cite{Knapp}, Chapter VII \S1, it is fairly straightforward to compute the following isomorphisms of $\SO(2)$-modules, where the sums are over all integers with the same parity:
$$
\pi_{1,n-1} \ \cong \ \bigoplus_{j\ {\rm even}}\ep_j,\qquad
\pi_{-1,n-1} \ \cong \ \bigoplus_{j\ {\rm odd}}\ep_j.
$$
The lemma then follows easily by combining the various elements of the proof.
\qed

From the previous two lemmas we can see that $\tr\CD^{\pm}_n(g_1)$ is non-zero only if $n=2$ and in that case $\tr\CD^{\pm}_2(g_1)=-2$, so $g_1(1)=-2M$.

It remains to show that $g_1(-1)=g_1(1)$.  Let $R_z$ be the right multiplication operator of $\SL_2(\R)$ on the space $C^{\infty}_{c}\left(\SL_2(\R)\right)$ of smooth, compactly supported functions on $\SL_2(\R)$.  That is, $R_zg(x)=g(xz)$ for all $x$, $z\in\SL_2(\R)$ and $g\in C^{\infty}_{c}\left(\SL_2(\R)\right)$.  Let $\pi$ be an irreducible, unitary representation of $\SL_2(\R)$.  The matrix $-1=-I_2$ is central in $\SL_2(\R)$ so $\pi(-1)$ commutes with $\pi(x)$ for all $x\in\SL_2(\R)$ and hence, by Schur's Lemma (\cite{Knapp}, Proposition 1.5) is scalar.  This means that for any irreducible, unitary representation $\pi$ on $\SL_2(\R)$ we have $\tr\pi\left(R_{-1}g_1\right)=\pi(-1)\tr\pi(g_1)$.  Thus we get
$$
g_1(-1) \=  R_{-1}g_1(1) \=  M\sum_{n=2}^{\infty}(n-1)\CD^{\pm}_n(-1)\tr\CD^{\pm}_n(g_1) \=  \CD^{\pm}_2(-1)g_1(1),
$$
where $\CD^{\pm}_2(-1)$ is a scalar, which we now compute.  For $g\in\CD^+_2$ and $z\in\C$ we have the action
$$
\CD^+_2\matrix{a}{b}{c}{d}g(z)=(-bz+d)^{-2}g\left(\frac{az-c}{-bz+d}\right).
$$
Hence
$$
\CD^+_2\matrixtwo{-1}{-1}g(z)=(-1)^{-2}g\left(\frac{-z}{-1}\right)=g(z).
$$
Similarly we get
$$
\CD^-_2\matrixtwo{-1}{-1}g(z)=g(z)
$$
so $\CD^{\pm}_2(-1)=1$ and the lemma is proved.
\qed

\subsection{Euler characteristics in the case of $\SL_4(\R)$}
We now consider the particular case $G=\SL_4(\R)$.  Let $K=\SO(4)$, a maximal compact subgroup of $G$ and let $\Ga$ be a discrete, cocompact subgroup of $G$.  Let $A$ be the rank one torus
$$
A=\left\{ \matrixfour{a}{a}{a^{-1}}{a^{-1}}:a>0\right\},\index{$A$}
$$
and let $B$ be the compact group
$$
B=\matrixtwo {\SO(2)}{\SO(2)}.\index{$B$}
$$
$B$ is a compact Cartan subgroup of
\begin{eqnarray*}
M & = & \SS\matrixtwo{\SL_2^{\pm}(\R)}{\SL_2^{\pm}(\R)} \\
  & \cong & \left\{(x,y)\in \Mat_2(\R)\times \Mat_2(\R)|\begin{array}{c}\det(x),\det(y)=\pm 1 \\ \det(x)\det(y)=1 \end{array}\right\}.\index{$M$}
\end{eqnarray*}
Let
$$
N=\matrix{I_2}{\Mat_2(\R)}{0}{I_2}\index{$N$}
$$
and let $P=MAN$ be a parabolic subgroup of $G$ with Levi component $MA$ and unipotent radical $N$. Let
$$
A^-=\left\{ \matrixfour{a}{a}{a^{-1}}{a^{-1}}:0<a<1\right\},\index{$A^-$}
$$
be the negative Weyl chamber in $A$ with respect to the root system given by the choice of parabolic.  Let $\CE_P(\Ga)$ \index{$\CE_P(\Ga)$} be the set of $\Ga$-conjugacy classes of elements $\ga\in\Ga$ which are conjugate in $G$ to an element of $A^-B$.  We shall use this notation for the rest of this section.

We say $g\in G$ is \emph{regular} \index{regular} if its centraliser is a torus and \emph{non-regular} \index{non-regular} (or \emph{singular}) \index{singular} otherwise.  Clearly, for $\ga\in\Ga$ regularity is a property of the $\Ga$-conjugacy class $[\ga]$.  Let $[\ga]\in\CE_P(\Ga)$ and write $G_{\ga}$ for the centraliser of $\ga$ in $G$.  The element $\ga$ is conjugate in $G$ to an element $a_{\ga}b_{\ga}\in A^-B$ and we define the length $l_{\ga}$ of $\ga$ to be $l_{\ga}=b(\log a_{\ga},\log a_{\ga})^{1/2}$.  It follows that if $\ga$ is regular then $G_{\ga}\cong AB$.

Let $K_{\ga}$ be a maximal compact subgroup of $G_{\ga}$, containing $B$. Then the group $B$ is a Cartan subgroup of $K_{\ga}$, the product $AB$ is a $\th$-stable Cartan subgroup of $G_{\ga}$ and $A$ is central in $G_{\ga}$.  If we let $\Ga_{\ga}=\Ga\cap G_{\ga}$ denote the centraliser of $\ga$ in $\Ga$ then $\Ga_{\ga}$ is discrete and cocompact in $G_{\ga}$.  Let $\Ga'$ be a torsion free subgroup of finite index in $\Ga$ and let $\Ga_{\ga}'=\Ga'\cap G_{\ga}$.  Then $\Ga_{\ga}'$ is a torsion free subgroup of finite index in $\Ga_{\ga}$.  We define $\Ga_{\ga,A}$, $\Ga'_{\ga,A}$ as in Section \ref{sec:EulerChars}.  The first higher Euler characteristic $\chi_1(\Ga_{\ga})$ of $\Ga_{\ga}$ is then defined as in (\ref{eqn:fhEC}) as:
\begin{equation}
\label{eqn:EulerChar}
\chi_1(\Ga_{\ga})=\frac{\left[\Ga_{\ga,A}:\Ga_{\ga,A}'\right]}{\left[\Ga_{\ga}:\Ga_{\ga}'\right]}\chi_1(\Ga_{\ga}').
\end{equation}

\begin{lemma}
\label{lem:ECharEqn}
For $\ga\in\CE_P(\Ga)$ we have that
$$
\chi_1(\Ga_{\ga})=\frac{C_{\ga}\vol(\Ga_{\ga}\bs G_{\ga})}{l_{\ga_0}},
$$
where $C_{\ga}$ is an explicit constant depending only on $\ga$, which is equal to one when $\ga$ is regular, and $\ga_0$ is the primitive geodesic underlying $\ga$.
\end{lemma}

\prf
Since $\Ga_{\ga}'$ is torsion free we may take from the second proposition in section 2.4 of \cite{Deitmar95} the equation
$$
\chi_1(\Ga_{\ga}')=\frac{C_{\ga}\vol(\Ga_{\ga}'\bs G_{\ga})}{\la_{\ga}'},
$$
where $\la_{\ga}'$ denotes the volume of the maximal compact flat in $\Ga'\bs G/K$ containing $\ga$, and $C_{\ga}$ is an explicit constant depending only on $\ga$.  We note that in \cite{Deitmar95} there is an extra factor in the equation which does not show up here.  The reason is that in \cite{Deitmar95} a differential form, not a measure, was constructed.  The value of the constant $C_{\ga}$ is given in \cite{Deitmar95} in terms of the root system of $G_{\ga}$ with respect to a Cartan subgroup.  In the case that $\ga$ is regular, as noted above we have that $G_{\ga}=AB$ and hence $G_{\ga}$ is a Cartan subgroup of itself.  It is then easy to see that value of $C_{\ga}$ in this case is one.

We denote by $\la_{\ga}$ the volume of the maximal compact flat in $\Ga\bs G/K$ containing $\ga$.  The values of $\la_{\ga}$ and $\la_{\ga}'$ are given by the volumes of the images of $\Ga_{\ga,A}\bs A$ and $\Ga_{\ga,A}'\bs A$ respectively under their embeddings into $\Ga\bs G/K$ and $\Ga'\bs G/K$ respectively.  Since $A$ is one-dimensional, $\la_{\ga}=l_{\ga_0}$.  Furthermore, from the definition of the groups $\Ga_{\ga,A}$ and $\Ga_{\ga,A}'$ it follows that
$$
\left[\Ga_{\ga,A}:\Ga_{\ga,A}'\right]=\la_{\ga}'/\la_{\ga}=\la_{\ga}'/l_{\ga_0}.
$$
Since
$$
\vol(\Ga_{\ga}'\bs G_{\ga})=\left[\Ga_{\ga}:\Ga_{\ga}'\right]\vol(\Ga_{\ga}\bs G_{\ga}),
$$
the lemma follows from (\ref{eqn:EulerChar}).
\qed

\begin{proposition}
\label{thm:ECharPos}
For all $\ga\in\CE_P(\Ga)$ the Euler characteristic $\chi_1(\Ga_{\ga})$ is positive.  If $\ga$ is regular then
$$
\chi_1(\Ga_{\ga})=\frac{\left[\Ga_{\ga,A}:\Ga_{\ga,A}'\right]}{\left[\Ga_{\ga}:\Ga_{\ga}'\right]}.
$$
\end{proposition}

\prf
If $\ga$ is regular then $G_{\ga}\cong AB$, $K_{\ga}\cong B$ and $\Ga_{\ga}'$ is a complete lattice in $G_{\ga}$ so $X_{\Ga_{\ga}'}=\Ga_{\ga}'\bs G_{\ga}/K_{\ga}\cong S^1$, the unit circle in $\C$.  The Betti numbers $h^j(S^1)$ are equal to zero for $j\neq 0,1$ and one for $j=0,1$, hence $\chi_1(\Ga_{\ga}')=1$.  The claimed value for $\chi_1(\Ga_{\ga})$ then follows immediately from (\ref{eqn:EulerChar}).

Suppose $\ga$ is not regular.  Since $\ga\in\CE_P(\Ga)$, we have that $\ga$ is conjugate in $G$ to $a_{\ga}b_{\ga}\in A^-B\subset AM$.  Let $(\si,V_{\si})$ be a finite dimensional representation of $M$ and let $K_M$ be the maximal compact subgroup $\SS(\O(2)\x\O(2))$ of $M$,  which contains the compact Cartan $B$ of $M$.  We saw above that there exist Euler-Poincar$\eac$ functions of $\si$ on $M$.  Let $f_{\si}$ be one such.  We denote by $M_{\ga}$ the centraliser of $b_{\ga}$ in $M$ and by $\CO^M_{b_{\ga}}(f_{\si})$ the orbital integral
$$
\int_{M/M_{\ga}} f_{\si}\left(xb_{\ga}x^{-1}\right) dx.
$$

From Proposition 1.4 of \cite{geom} we get the equation
$$
\CO^M_{b_{\ga}}(f_{\si})=C_{\ga}\tr\si(b_{\ga}),
$$
where $C_{\ga}$ is the constant from Lemma \ref{lem:ECharEqn}.  Together with Lemma \ref{lem:ECharEqn} this gives us
$$
\chi_1(\Ga_{\ga})
  =\frac{\CO^M_{b_{\ga}}(f_{\si})\vol(\Ga_{\ga}\bs G_{\ga})}{l_{\ga_0}\tr\si(b_{\ga})}.
$$
Choosing $\si=1$, the trivial representation of $M$, this simplifies to
\begin{equation}
\label{eqn:EChar}
\chi_1(\Ga_{\ga})=\frac{\CO^M_{b_{\ga}}(f_1)\vol(\Ga_{\ga}\bs G_{\ga})}{l_{\ga_0}}.
\end{equation}

To complete the proof of the theorem we shall show that the orbital integral $\CO^M_{b_{\ga}}(f_1)$ is positive.  In the case we are considering
$$
b_{\ga}=\matrixtwo{\pm 1}{\pm 1},
$$
where $1$ denotes the identity matrix in $\SL_2(\R)$, hence it is central and $M_{\ga}=M$ so we have simply that
\begin{equation}
\label{eqn:OrbInt}
\CO^M_{b_{\ga}}(f_1)=f_1(b_{\ga}).
\end{equation}

The group $\bar{M}\cong\SL_2(\R)\x\SL_2(\R)$ is the connected component of $M\cong\SS(\SL_2^{\pm}(\R)\x\SL_2^{\pm}(\R))$.  $M$ has a maximal compact subgroup $K_M\cong\O(2)\x\O(2)$ and compact Cartan $T_M\cong\SO(2)\x\SO(2)$.  We have that $T_M \bar{M}=\bar{M}$.  Hence by Lemma \ref{lem:EPfnrestrict}, since $f_1$ is an Euler-Poincar$\eac$ function for the trivial representation on $M$ we have that $\bar{f}_1=f_1|_{\bar{M}}$ is an Euler-Poincar$\eac$ function for the trivial representation on $\bar{M}$.

Let $g_1$, $h_1$ be Euler-Poincar$\eac$ functions of the trivial representation on $\SL_2(\R)$.  By Lemma \ref{lem:EPfnproduct} the function $\tilde{f}_1=g_1 h_1$ is an Euler-Poincar$\eac$ function of the trivial representation on $\bar{M}$.

We recall that
$$
b_{\ga}=\matrixtwo{\pm 1}{\pm 1},
$$
which is central in $M$ and deduce from Lemma \ref{lem:EPCentral} that
$$
f_1(b_{\ga})=\bar{f}_1(b_{\ga})=\tilde{f}_1(b_{\ga})=g_1(\pm 1)h_1(\pm 1).
$$
From Lemmas \ref{lem:EPCentral} and \ref{lem:EPfnvalue} it follows that
$$
g_1(1)=g_1(-1)=h_1(1)=h_1(-1)\in\R
$$
and so $f_1(b_{\ga})$ is positive.
\qed

\subsection{The unitary dual of $K_M$}

Let $K_M=K\cap M\cong\SS(\O(2)\x\O(2))$.\index{$K_M$}  We shall need in the proof of later results to know the unitary dual $\widehat{K_M}$ of $K_M$, which is given in the following proposition.

First we define the following one dimensional representations of $\SO(2)$ and $\SO(2)\times \SO(2)$.
$$
\ep_l R(\th)=e^{il\th},\ \ \textnormal{for all}\ l\in\Z,
$$
$$
\ep_{l,k} \matrixtwo{R(\th)}{R(\eta)} = e^{i(l\theta + k\eta)},\ \ \textnormal{for all}\ l,k\in\Z,
$$
where
$$
R(\th)=\matrix{\cos\th}{-\sin\th}{\sin\th}{\cos\th}\in \SO(2).
$$
Note that $\ep_0$ and $\ep_{0,0}$ are the trivial representation of their respective groups.

\begin{proposition}
\label{pro:KMDual}
\index{unitary dual!of $K_M$}
For $l,k\in\Z$ not both zero there are unique representations $\d_{l,k}$ of $K_M$ on $\C^2$ with
$$
\d_{l,k}|_{\SO(2)\times \SO(2)}=\ep_{l,k}\oplus\ep_{-l,-k},
$$
and
$$
\d_{l,k}\matrixfour{-1}{1}{-1}{1}(z_1,z_2)=(z_2,z_1).
$$
We can also define the representation $\d$ of $K_M$ on $\C$ by
$$
\d(X,Y)(z)=(\det X)z=(\det Y)z.
$$
We have that $\widehat{K_M}=\{triv,\d,\d_{l,k}:l,k\in\Z\ \textrm{not both zero}\}$.
\end{proposition}
\prf
We have seen that $\widehat{\SO(2)}=\{\ep_l:l\in\Z\}$.
In general, for two locally compact groups $H$ and $K$, the map $(\si,\tau)\mapsto\si\otimes\tau$ defines an isomorphism $\hat{H}\times\hat{K}\cong\widehat{H\times K}$.  Thus the map from $\widehat{\SO(2)}\times\widehat{\SO(2)}$ to $\widehat{\SO(2)\times \SO(2)}$ given by $(\tau_1,\tau_2)\mapsto\tau_1\otimes\tau_2$ is an isomorphism.  Hence we have that $\widehat{\SO(2)\times \SO(2)}$ is the set
$$
\{\ep_l\ox\ep_k|\,l,k\in\Z\}=\{\ep_{l,k}|\,l,k\in\Z\}.
$$

Let
$$
T = \matrixfour{-1}{1}{-1}{1}
$$
and
$$
R(\th,\eta)=\matrixtwo{R(\th)}{R(\eta)}.
$$

For $l,k\in\Z$ we note that
$$
\d_{l,k}(T)\d_{l,k}(T)(z_1,z_2)=(z_1,z_2)=\d_{l,k}(T^2)(z_1,z_2)
$$
and
\begin{eqnarray*}
\d_{l,k}(T)\d_{l,k}(R(\th,\eta))\d_{l,k}(T)(z_1,z_2) & = & (e^{-i(l\th+k\eta)}z_1,e^{i(l\th+k\eta)}z_2) \\
    & = & \d_l(R(-\th,-\eta))(z_1,z_2) \\
    & = & \d_l(T R(\th,\eta) T)(z_1,z_2),
\end{eqnarray*}
so $\d_{l,k}$ is indeed a representation.  We shall show that the representations given in the proposition are in fact the only irreducible unitary representations of $K_M$.

Let $\tau\in\widehat{K_M}$.  Then by \cite{Knapp}, Theorem 1.12(d) the representation $\tau$ restricted to ${\SO(2)\x\SO(2)}$ is a direct sum of irreducible representations, that is
$$
\tau|_{\SO(2)\x\SO(2)}=\ep_{l_1,k_1}\oplus\ep_{l_2,k_2}\oplus\cdots\oplus\ep_{l_n,k_n},
$$
for some $l_1,\dots,l_n,k_1,\dots,k_n\in\Z$.  Let $v\in\tau$ be an $(\SO(2)\x\SO(2))$-eigenvector with eigenvalue $e^{i(l_1\th+k_1\eta)}$.  Then the equation
$$
TR(-\th,-\eta)=R(\th,\eta)T
$$
tells us that $\tau(T)v$ is also an $\SO(2)$-eigenvector, with eigenvalue $e^{-i(l_1\th+k_1\eta)}$.  If $\tau(T)v$ is a scalar multiple of $v$ then $l_1=k_1=0$ and the equation $T^2=I$ tells us that $\tau=triv$ or $\d$.  Otherwise $\tau|_{\SO(2)\x\SO(2)}=\ep_{l_1,k_1}\oplus\ep_{-l_1,-k_1}$ and $\tau=\d_{l_1,k_1}$.

For $l=k=0$ we have $\d_{0,0}\cong triv\oplus\d$.  For all other values of $l$ and $k$ the representation $\d_{l,k}$ is irreducible.  Also, for $l,k\in\Z$ we have that $\d_{l,k}$ is unitarily equivalent to $\d_{-l,-k}$.  There are no other equivalences between the representations $\d_{l,k}$.

The proposition follows.
\qed

\subsection{Infinitesimal characters}

Let $G$ be a connected reductive group with maximal compact subgroup $K$, real Lie algebra $\g_{\R}$ and complexified Lie algebra $\g$.  Let $\pi$ be a representation of $G$.  
Let $\pi_K$ be the $(\g,K)$-module of $K$-finite vectors in $\pi$.

Suppose that $\pi$ is an admissible representation of $G$ such that $\pi(Z)$ acts as a scalar on $\pi_K$ for all $Z$ in the centre $\z_G$ of the universal enveloping algebra of $\g$.  In particular this condition is satisfied when $\pi$ is irreducible admissible.  
Let then $\h$ be a Cartan subalgegra of $\g$ and let $\La_\pi\in\h^*$ be a representative of the infinitesimal character of $\pi$.

Let $G=\SL_4(\R)$ and the subgroups $K$ and $P=MAN$ be as above.  Let $\h$ be the diagonal subalgebra of $\g$.  In this case the Weyl group $W(\g,\h)$ acts on $\h$ by interchanging elements of the leading diagonal.  We shall see in the next section that the analytic properties of the zeta functions considered there are related to the infinitesimal characters of elements of $\hat{G}$.  The following proposition gives us information about these infinitesimal characters which will be required in the following section.

Let $\rho_P$ be the half-sum of the positive roots of the system $(\g,\a)$, where $\a$ is the complexified Lie algebra of $A$, so that
$$
\rho_P \matrixfour{a}{a}{-a}{-a} = 4a.
$$
Let $\si$ be a finite dimensional virtual representation of $M$, whose $K_M$ types are all contained in the set $\{triv,\d,\d_{l,k}:l,k\in\{0,2\}\}$.  Let $\hat{M}_{\si}$ be the subset of all $\xi\in\hat{M}$ such that $\tr\xi(f_{\si})=0$ for all Euler-Poincar$\eac$ functions $f_{\si}$ for $\si$ on $M$.  (Note that the value of $\tr\xi(f_{\si})=0$ depends only on $\xi$ and $\si$ and not on the choice of Euler-Poincar$\eac$ function $f_{\si}$.)  Let $\hat{G}_{\si}$ be the set of all elements of $\hat{G}$ except for: the trivial representation; representations induced from parabolic subgroups other than $P=MAN$ and representations induced from $\xi\in\hat{M}_{\si}$.  We define an order on the real dual space of $\a$ by $\la>\mu$ if and only if $(\la-\mu)=t\rho_P$ for some $t>0$.

\begin{proposition}
\label{pro:InfChars}
Let $\si$ be a finite dimensional virtual representation of $M$ and let $\pi\in\hat{G}_{\si}$.  Then the infinitesimal character $\La_{\pi}$ of $\pi$ satisfies
$$
\Re(w\La_{\pi})|_{\a}\geq-\rez{2}\rho_P\ \textrm{ or }\ -\rho_P\geq\Re(w\La_{\pi})|_{\a}
$$
for every $w\in W(\g,\h)$.
\end{proposition}
\prf
We shall prove the proposition by considering different cases in turn.  By Theorem \ref{thm:SL4Dual} we have the following cases to consider.  The case when $\pi$ is a principal series representation induced from $P=MAN$; the case when $\pi$ is a complementary series representation induced from $P=MAN$; the case when $\pi$ is a limit of complementary series representation and the case when $\pi$ is one of the representations $\pi_m$ for $m\in\N$.

First we consider the principal series representations.  Let $\pi=\pi_{\xi,\nu}=\Ind_P^G(\xi\ox\nu)$ be induced from $P$, where $\xi$ is an irreducible, unitary representation of $M \cong \SS(\SL_2^{\pm}(\R)\times \SL_2^{\pm}(\R))$ and $\nu\in\a^*$ is imaginary (ie. $\nu\in i\a_{\R}^*$, where $\a_{\R}^*$ is the real dual space of $\a$).

Let $\xi'$ be an irreducible subspace of $\xi |_{\SL_2(\R)\times \SL_2(\R)}$.  Then $\xi'$ has infinitesimal character $\La_{\xi'}$ and $\La_{\xi}=\La_{\xi'}$.  We have that $\xi'\cong\xi_1 \otimes \xi_2$ where $\xi_1$ and $\xi_2$ are irreducible unitary representations of $\SL_2(\R)$.  

To limit the possibilities for $\xi_1$ and $\xi_2$ that we need to consider we use the double induction formula (see \cite{Knapp}, (7.4)).
\begin{lemma}(Double induction formula)

Let $M_{\diamond}A_{\diamond}N_{\diamond}$ be a parabolic subgroup of $M$, so that $M_{\diamond}(A_{\diamond}A)(N_{\diamond}N)$ is a parabolic subgroup of $G$.  If $\si_{\diamond}$ is a unitary representation of $M_{\diamond}$ and $\nu_{\diamond}\in\a_{\diamond}^*=(\Lie_{\C}A_{\diamond})^*$, $\nu\in\a^*$, then there is a canonical equivalence of representations
$$
\Ind_{MAN}^G\left(\Ind_{M_{\diamond}A_{\diamond}N_{\diamond}}^M\left(\si_{\diamond}\ox\nu_{\diamond}\right)\ox\nu\right)\cong\Ind_{M_{\diamond}(A_{\diamond}A)(N_{\diamond}N)}^G\left(\si_{\diamond}\ox(\nu_{\diamond}+\nu)\right).
$$
\end{lemma}
\qed

We may assume that neither $\xi_1$ nor $\xi_2$ are induced since then, by the double induction formula, we would have the case that $\pi$ was induced from a parabolic other than $P=MAN$, which was excluded.  By Theorem \ref{thm:SL2Dual}, the remaining possibilities for $\xi_1$ and $\xi_2$ are the trivial representation and the discrete series and limit of discrete series representations.

Let $\La_{\pi}$ be the infinitesimal character of $\pi$ and $\La_{\xi_{i}}$ be the infinitesimal character of $\xi_{i}$, then $\La_{\xi}=\La_{\xi_1}+\La_{\xi_2}$.  Recall that $\h$ is the diagonal subalgebra of $\g$.  We lift $\La_{\xi}$ and $\nu$ to $\h$ by defining $\La_{\xi}$ to be zero on $\a$ and $\nu$ to be zero on the diagonal elements of $\m$, so that $\La_{\pi} = \La_{\xi}+\nu$ (\cite{Knapp}, Proposition 8.22).  The Weyl group $W=W(\g,\h)$ acts on $\La_{\pi}$.  Let $w\in W$, we will show that either $\tr\xi(f_{\si})=0$ or
$$
\Re(w\La_{\pi} |_{\a})\geq-\frac{1}{2}\rho_P
$$
or
$$
-\rho_P\geq\Re(w\La_{\pi} |_{\a}).
$$

First we take $\xi_1$ and $\xi_2$ to be the trivial representation.  This gives us that
$
\La_{\xi} (\diag({s},{-s},{t},{-t})) = s+t.
$
If
$
\nu (\diag({a},{a},{-a},{-a})) = \al a ,\ \ \ \ \al\in\textit{i}\R,
$
then
\begin{eqnarray}
\label{eqn:InfChar}
\lefteqn{\La_{\pi}(\diag({a},{b},{c},{-a-b-c}))} \nonumber \\
  & = & \La_{\xi}(\diag({\frac{a-b}{2}},{\frac{b-a}{2}},{c+\frac{a+b}{2}},{-c-\frac{a+b}{2}})) + \al\frac{a+b}{2} \\
  & = & \frac{a-b}{2} + \left(c+\frac{a+b}{2}\right) + \al\frac{a+b}{2} \nonumber \\
  & = & a+c+\frac{\al}{2}(a+b). \nonumber
\end{eqnarray}

If we let $w=1$ then from above we see that $\Re(\La_{\pi}|_{\a})=0$.  If we take $w$ to be the transposition interchanging $b$ and $c$ we get
$$
w\La_{\pi} (\diag({a},{b},{c},{-a-b-c})) = a+b+\frac{\al}{2}(a+c),
$$
the real part of which when restricted to $\a$ gives $\frac{1}{2} \rho_P$.  All other Weyl group elements are dealt with similarly and we see that $-\frac{1}{2}\rho_P\leq\Re(\La_{\pi}|_{\a})\leq\frac{1}{2} \rho_P$ for all $w\in W$.

It remains to consider the case when either or both of $\xi_1$ and $\xi_2$ are a discrete series representation or limit of discrete series representation with parameter $m_i$.  We set $\CD_0^+=\CD_0^- = triv$ and let $\xi_i = \CD_{m_i}^+$ or $\xi_i = \CD_{m_i}^-$, where $m_i \geq 0$, and $m_1$ and $m_2$ are not both zero.  From (\ref{eqn:EPFnId}) we get
\begin{equation}
\label{eqn:TraceXif1}
\tr\xi(f_{\si})=\sum_{p=0}^{\dim \p_M}\ (-1)^p\ \dim \left(V_{\xi}\ox \bwedge^p \p_M \ox V_{\si}\right)^{K_M}.
\end{equation}
We shall use Proposition \ref{pro:KMDual} to examine the $K_M$ types to limit the possibilities for $m_1$, $m_2$ for which $\tr\xi(f_{\si})\neq 0$.

\begin{lemma}
\label{lem:XiKMTypes}
For $m_1$, $m_2$ both non-zero we have the following isomorphism of $K_M$-modules:
$$
V_{\xi}\cong \bigoplus_{{j\geq m_1, j\equiv m_1\mod 2} \atop {k\geq m_2, k\equiv m_2\mod 2}}\d_{j,k}\oplus\d_{-j,k}.
$$
If $m_1=0$ we have
$$
V_{\xi}\cong \bigoplus_{k\geq m_2, k\equiv m_2\mod 2}\d_{0,k},
$$
and for $m_2=0$ we have
$$
V_{\xi}\cong \bigoplus_{j\geq m_1, j\equiv m_1\mod 2}\d_{j,0}.
$$
\end{lemma}
\prf
This follows from Lemma \ref{lem:DiscSer}.
\qed

\begin{lemma}
\label{lem:pMKMTypes1}
We have the following isomorphisms of $K_M$-modules:
\begin{eqnarray*}
  \bwedge^0 \p_M & \cong & triv \\
  \bwedge^1 \p_M & \cong & \d_{2,0}\oplus\d_{0,2} \\
  \bwedge^2 \p_M & \cong & \d\oplus\d\oplus\d_{2,2}\oplus\d_{2,-2} \\
  \bwedge^3 \p_M & \cong & \d_{2,0}\oplus\d_{0,2} \\
  \bwedge^4 \p_M & \cong & triv.
\end{eqnarray*}
\end{lemma}
\prf
$K_M$ acts on $\p_M$ by the adjoint representation and we can compute
$$
\p_M \cong \d_{2,0}\oplus\d_{0,2}.
$$
The other isomorphisms follow straightforwardly from this.
\qed

Let $v=v_1\ox v_2\ox v_3\in V_{\xi}\ox\bwedge^*\p_M\ox V_{\si}$, where the $v_i$'s all lie in a single $K_M$ type of their respective representation spaces.  Lemma~\ref{lem:pMKMTypes1} tells us that $v_2$ is in one of the following $K_M$ types: $triv$, $\d$, $\d_{2,0}$, $\d_{0,2}$, $\d_{2,2}$, $\d_{2,-2}$ by our assumption on $\si$, the possibilities for the $K_M$ type containing $v_3$ are also the same.  It follows that $\tr\xi(f_{\si})$ is non-zero only if $m_1,m_2\in\{0,2,4\}$.

It follows from the exact sequences in the proof of Lemma \ref{lem:DiscSer} that $\CD^+_m\oplus\CD^-_m\subset\CP^{\pm,(m-1)\rho_1}$, where the index on $\CP$ is $+$ if $m$ is even and $-$ if $m$ is odd.  The infinitesimal character of $\CP^{\pm,(m-1)\rho_1}$ is simply $(m-1)\rho_1$ (see \cite{Knapp}, Proposition 8.22), hence it follows that the infinitesimal characters of $\CD^+_m$ and $\CD^-_m$ are both equal to $(m-1)\rho_1$.  This gives us in the cases where $\tr\xi(f_{\si})$ is non-zero:
\begin{equation}
\label{eqn:infchar}
\La_{\xi} (\diag({s},{-s},{t},{-t})) = (m_1 -1)s+(m_2 -1)t,\ \ m_1,m_2\in\{0,2,4\}.
\end{equation}
By computing the action of the different Weyl group elements we can see that in all cases either $w\La_{\xi}|_{\a}\geq-\rez{2}\rho_P$ or $w\La_{\xi}|_{\a}\leq-\rho_P$.  Since $\La_{\pi}=\La_{\xi}+\nu$ and $\nu$ is imaginary we see that $\Re(\La_{\pi}|_{\a})=\La_{\xi}|_{\a}$, so the claim follows.

The complementary series representations are dealt with similarly.  We recall from Theorem \ref{thm:SL4Dual} that the complementary series induced from $P$ are $\pi=\Ind^G_P\left(\bar{\pi}_m\ox t\rho_P\right)$, for $m\in\N$ and $0<t<\rez{2}$.  We may argue as above to find the possibilities for $m$ such that $\tr\bar{\pi}_m(f_{\si})\neq 0$.  In this way we see that there are two possibilities for $m$ for which $\tr\bar{\pi}_m(f_{\si})\neq 0$, namely $m=2$ and $m=4$.  In the first case, $w\La_{\pi}|_{\a}\geq -\rez{2}\rho_P$, for all $w\in W(\g,\h)$.  In the second case we have
$$
\La_{\pi}(\diag({a},{b},{c},{-a-b-c}))=(3+2t)a+2tb+3c.
$$
If $w\in W(\g,\h)$ is the element which swaps the first and fourth diagonal entries then $w\La_{\pi}|_{\a}=-\frac{3}{2}\rho_P$.  In all other cases we have $w\La_{\pi}|_{\a}\geq -\rez{2}\rho_P$.

The limit of complementary series representations are closely related to the family of representations $\pi_m$, $m\in\N$, so we shall deal with them together.

For $m\in\N$, we denote by $\bar{\pi}_m$ the representation of $M$ induced from the representation $\CD_m^+\otimes\CD_m^+$ of $\SL_2(\R)\x\SL_2(\R)$.  For $m\in\N$ we have the limit of complementary series representation given as an irreducible unitary subrepresentation of $I_m=\Ind_P^G(\bar{\pi}_m\ox\rez{2}\rho)$, which we will denote by $\pi_m^c$.  The representations $\pi_m$ are the Langlands quotients of the $I_m$.

Let $\La_{\pi_m}$ be the infinitesimal character of $\pi_m$ and $\La_{\pi_m^c}$ the infinitesimal character of $\pi_m^c$.  Clearly $\La_{\pi_m}=\La_{\pi_m^c}=\La_{\bar{\pi}_m}+\rez{2}\rho_P$, so we need only consider $\La_{\pi_m}$.  The value of
$$
\La_{\pi_m}\matrixfour{a}{b}{c}{-a-b-c}
$$
is equal to
$$
\La_{\bar{\pi}_m}\matrixfour{\frac{a-b}{2}}{\frac{b-a}{2}}{\frac{a+b+2c}{2}}{\frac{a+b-2c}{2}} + \rez{2}\rho_P \matrixfour {\frac{a+b}{2}}{\frac{a+b}{2}}{-\frac{a+b}{2}}{-\frac{a+b}{2}}
$$
\begin{eqnarray*}
  & = & (m-1)\frac{a-b}{2} + (m-1)(c+\frac{a+b}{2}) + (a+b) \\
  & = & ma+b+(m-1)c.
\end{eqnarray*}

Restricting to $\a$ we get
$$
\La_{\pi_m}\matrixfour{a}{a}{-a}{-a} = 2a = \frac{1}{2}\rho_P.
$$
The Weyl group $W$ acts on $\La_{\pi_m}$.  Now for all elements $w\in W$ except one we get $w\La_{\pi_m}|_{\a}\geq-\frac{1}{2}\rho_P$ for all $m\geq 1$.  The exception is $w_1$ which swaps the first and fourth diagonal entries.  For this we have
$$
w_1\La_{\pi_m}\matrixfour{a}{a}{-a}{-a} = -2(m-1)a.
$$
Thus $w_1\La_{\pi_m}|_{\a}\geq -\frac{1}{2}\rho_P$ if $m=1,2$ and if $m\geq 3$ then $w_1\La_{\pi_m}|_{\a}\leq-\rho_P$.

This completes the proof of the proposition.
\qed

\section{Analysis of the Ruelle Zeta Function}
    \label{ch:Ruelle}

Let $G=\SL_4(\R)$ \index{$G$} and $\Ga\subset G$ \index{$\Ga$} be discrete and cocompact.  Let $\g_{\R}=\sl_4(\R)$ and $\g=\sl_4(\C)$ \index{$\g$} be respectively the Lie algebra and complexified Lie algebra of $G$.  As in the previous section, all Haar measures will be normalised as in \cite{HarishChandra75}, \S 7.  Recall that this normalisation depends on the choice of an invariant bilinear form $b$ on $\g$.  Let $b$ be the following multiple of the trace form on $\g$:
\begin{equation}
\label{eqn:norm}
\index{$b(X,Y)$}
b(X,Y)=16\tr(XY).
\end{equation}
We choose this normalisation in order to get the first zero of the Ruelle zeta function at the point $s=1$ in Theorem \ref{thm:RuelleMain} below.  Let $K \subset G$ \index{$K$} be the maximal compact subgroup $\SO(4)$.  Let $\k_{\R} \subset \g_{\R}$ be its Lie algebra and let $\p_{\R} \subset \g_{\R}$ be the orthogonal space of $\k_{\R}$ with respect to the form $b$.  Then $b$ is positive definite and Ad($K$)-invariant on $\p_{\R}$ and thus defines a $G$-invariant metric on $X=G/K$, the symmetric space attached to $G$.

Let
$$
A=\left\{ \matrixfour{a}{a}{a^{-1}}{a^{-1}}:a>0\right\},\index{$A$}
$$
$$
B=\matrixtwo {\SO(2)}{\SO(2)}.\index{$B$}
$$
$B$ is a compact Cartan subgroup of
\begin{eqnarray*}
M & = & \SS\matrixtwo{\SL_2^{\pm}(\R)}{\SL_2^{\pm}(\R)} \\
  & \cong & \left\{(x,y)\in \Mat_2(\R)\times \Mat_2(\R)|\begin{array}{c}\det(x),\det(y)=\pm 1 \\ \det(x)\det(y)=1 \end{array}\right\}.
\index{$M$}
\end{eqnarray*}
We also define
$$
N=\matrix{I_2}{\Mat_2(\R)}{0}{I_2}\textrm{ and }\ \bar{N}=\matrix{I_2}{0}{\Mat_2(\R)}{I_2},\index{$N$}\index{$\bar{N}$}
$$
and set $K_M=K\cap M$.\index{$K_M$}

Let $\m$ denote the complexified Lie algebra of $M$ and let $\m=\k_M\oplus\p_M$ \index{$\p_M$} be its polar decomposition, where $\k_M$ is the complexified Lie algebra of $K_M=K\cap M$.  Let $\h$ \index{$\h$} be the Cartan subalgebra of $\g$ consisting of all diagonal matrices, and let $\a$\index{$\a$}, $\n$ \index{$\n$} and $\bar{\n}$ \index{$\bar{\n}$} be the  complexified Lie algebras of $A$, $N$ and $\bar{N}$ respectively.  Let $P$ denote the parabolic with Langlands decomposition $P=MAN$ \index{$P$} and $\bar{P}$ the opposite parabolic with Langlands decomposition $\bar{P}=MA\bar{N}$.  Let $\rho_P$ be the half-sum of the positive roots of the system $(\g,\a)$, so that
$$
\rho_P \matrixfour{a}{a}{-a}{-a} = 4a.\index{$\rho_P$}
$$

Let
$$
H_1=\frac{1}{8}\matrixfour{-1}{-1}{1}{1} \in \a_{\R}.\index{$H_1$}
$$
Then it follows that $b(H_1)=b(H_1,H_1)=1$ and $\rho_P(H_1)=-\rez{2}$.  Let $A^- = \{\exp(tH_1):t>0\}$ \index{$A^-$} be the negative Weyl chamber in $A$.  Let $\CE_P(\Ga)$ \index{$\CE_P(\Ga)$} be the set of all conjugacy classes $[\ga]$ in $\Ga$ such that $\ga$ is conjugate in G to an element $a_{\ga}b_{\ga}$ of $A^-B$.

An element $\ga\in\Ga$ is called \emph{primitive} if for $\d\in\Ga$ and $n\in\N$ the equation $\d^n =\ga$ implies that $n=1$.  Clearly primitivity is invariant under conjugacy, so we may view it as a property of conjugacy classes.  Let $\CE_P^p(\Ga)\subset\CE_P(\Ga)$ \index{$\CE_P^p(\Ga)$} be the subset of primitive classes.

For $[\ga]\in\CE_P(\Ga)$ we define the length of $\ga$ to be $l_{\ga}=b(\log a_{\ga},\log a_{\ga})^{1/2}$\index{$l_{\ga}$}.  Let $G_{\ga}$ \index{$G_{\ga}$} be the centraliser of $\ga$ in $G$, let $\Ga_{\ga}=\Ga\cap G_{\ga}$ \index{$\Ga_{\ga}$} be the centraliser of $\ga$ in $\Ga$ and let $\chi_1(\Ga_{\ga})$ be the first higher Euler characteristic as in the previous section.

\subsection{The Selberg trace formula}
\label{sec:STF}

Let $\hat{G}$ be the unitary dual of $G$.  The group $G$ acts on the Hilbert space $L^2 (\Ga\bs G)$ by translations from the right.  Since $\Ga\bs G$ is compact this representation decomposes discretely:
$$
L^2 (\Ga\bs G)=\bigoplus_{\pi\in\hat{G}} N_{\Ga}(\pi)\pi
$$
with finite multiplicities $N_{\Ga}(\pi)$\index{$N_{\Ga}(\pi)$}, (see \cite{Gelfand69}).

Recall that for $\ga\in\Ga$ we denote by $G_{\ga}$ and $\Ga_{\ga}$ the centraliser of $\ga$ in $G$ and $\Ga$ respectively.  For a function $f$ on $G$, denote by $\CO_{\ga}(f)$\index{$\CO_{\ga}(f)$} the orbital integral\index{orbital integral}
$$
\CO_{\ga}(f)=\int_{G/G_{\ga}}f(x\ga x^{-1})\,dx.
$$

The \emph{Selberg trace formula} \index{Selberg trace formula}(see \cite{Wallach76}, Theorem 2.1) says that for suitable functions $f$ on $G$ the following identity holds:
\begin{equation}
\label{eqn:STF}
\sum_{\pi\in\hat{G}}N_{\Ga}(\pi)\tr\pi(f)=\sum_{[\ga]}\vol(\Ga_{\ga}\bs G_{\ga})\CO_{\ga}(f),
\end{equation}
where the sum on the right is over all conjugacy classes in $\Ga$.  The set of suitable functions includes, but is not limited to, all $\dim G+1$ times continuously differentiable functions of compact support on $G$.  In fact, we shall need to extend the set of test functions for which the trace formula is valid.

\begin{lemma}
Let $d\in\N$, $d\geq 16$, that is $d=2d'$ for some $d'>\dim G/2$.  Let $f$ be a $d$-times differentiable function on $G$ such that $Df\in L^1(G)$ for all left invariant differentiable operators $D$ on $G$ with complex coefficients and of degree $\leq d$.  Then the trace formula is valid for $f$.
\end{lemma}
\prf
This is Lemma 1.3 of \cite{Deitmar04} in the case $G=\SL_4(\R)$.
\qed

For $g\in G$ and $V$ any complex vector space on which $g$ acts linearly let $E(g|V)\subset\R_+$ be defined by
$$
E(g|V)=\{|\mu|:\mu\textrm{ is an eigenvalue of }g\textrm{ on }V\}.
$$
Let $\la_{\min}(g|V)=\min(E(g|V))$ and $\la_{\max}(g|V)=\max(E(g|V))$.

For $am\in AM$ define
$$
\la(am)=\frac{\la_{\min}(a|\n)}{\la_{\max}(m|\g)^2},
$$
where we are considering the adjoint action of $G$ on $\n$ and $\g$ resepectively.  Define the set
$$
(AM)^{\sim}=\{am\in AM:\la(am)<1\}.
$$
An element of $G$ is said to be \emph{elliptic} \index{elliptic element} if it is conjugate to an element of a compact torus.  Let $M_{\ell}$ denote the set of elliptic elements in $M$.

\begin{lemma}
The set $(AM)^{\sim}$ has the following properties:

(1) $A^-M_{\ell}\subset(AM)^{\sim}$;

(2) $am\in(AM)^{\sim}\Ra a\in A^-$;

(3) $am,a'm'\in(AM)^{\sim},g\in G$ with $a'm'=gamg^{-1}\Ra a=a',g\in AM$.
\end{lemma}
\prf
See \cite{geom}, Lemma 2.4.
\qed

For the construction of our test function we shall need, for given $s\in\C$ and $j\in\N$, a smooth, conjugation invariant, $j$-times continuously differentiable function on $AM$, with support in $(AM)^{\sim}$.  Further, we require that at each point $ab\in A^-B$ the function takes the value $l_a^{j+1}e^{-sl_a}$.  In \cite{geom} a function $g_s^j$ \index{$g_s^j$} is constructed with these properties, with the one difference that there the positive Weyl chamber $A^+$ is used instead of the negative Weyl chamber $A^-$.  With only very minor modification the construction in \cite{geom} will yield a function with our required properties, which we shall also call $g_s^j$. 

Let $\eta:N\ra\R$ be a smooth, non-negative function of compact support, which is invariant under conjugation by elements of $K_M$ and satisfies
$$
\int_N\eta(n)\ dn=1.
$$
Let $f:M\ra\C$ be a smooth, compactly supported function, invariant under conjugation by $K_M$.  Suppose further that the orbital integrals of $f$ on $M$ satisfy
$$
\CO_m^M(f)=\int_{M/M_m} f(xmx^{-1})\ dx=0
$$
whenever $m$ is not conjugate to an element of $B$, where $M_m$ denotes the centraliser of $m$ in $M$.  The group $AM$ acts on $\n$ according to the adjoint representation.

Given these data we define $\Phi=\Phi_{f,j,s}:G\ra\C$ by
\begin{equation}
\label{eqn:TestFunction}
\index{$\Phi$}
\Phi(knam(kn)^{-1})=\eta(n)f(m)\frac{g_s^j(am)}{\det(1-(am)^{-1}|\bar{\n})},
\end{equation}
for $k\in K,n\in N,am\in AM$.

To see that $\Phi$ is well defined we recall that, by the decomposition $G=KNAM$, every $g\in G$ which is conjugate to an element of $(AM)^{\sim}$ can be written in the form $knam(kn)^{-1}$.  By the properties of $(AM)^{\sim}$ we see that two of these representations can only differ by an element of $K_M$.  The components of the function $\Phi$ are all invariant under conjugation by $K_M$, and we note that $\det(1-(am)^{-1}|\n)\neq 0$ for all $am\in(AM)^{\sim}$, so we can see that $\Phi$ is well-defined.

\begin{proposition}
\label{pro:STF2}
The function $\Phi$ is $(j-14)$-times continuously differentiable.  For $j$ and $\Re(s)$ large enough it goes into the trace formula and we have:
\begin{equation}
\label{eqn:STF2}
\sum_{\pi\in\hat{G}} N_{\Ga}(\pi)\tr\pi(\Phi)=\sum_{[\ga]\in\CE_P(\Ga)}\vol(\Ga_{\ga}\bs G_{\ga})\CO^M_{b_{\ga}}(f)\frac{l_{a_{\ga}}^{j+1}e^{-sl_{a_{\ga}}}}{\det(1-(a_{\ga}b_{\ga})^{-1}|\bar{\n})}.
\end{equation}
\end{proposition}
\prf
Noting that $2\dim\n+\dim\k=14$ we see that this follows from Proposition 2.5 of \cite{geom}.  This proposition was proved in the case that the function $f$ is an Euler-Poincar$\eac$ function for some finite dimensional representation of $M$.  However the only properties of Euler-Poincar$\eac$ functions used in the proof are those given above for $f$, namely that it is smooth, of compact support, invariant under conjugation by $K$, and the orbital integrals satisfy the given condition.
\qed

\subsection{The Selberg zeta function}
\label{sec:SelbergZetaFn}

An element $g\in G$ is said to be \emph{weakly neat} \index{weakly neat} if the adjoint $\Ad(g)$ has no non-trivial roots of unity as eigenvalues.  A subgroup $H$ of $G$ is said to be weakly neat if every element of $H$ is weakly neat.  Let $[\ga]\in\CE^p_P(\Ga)$ so that $\ga$ is conjugate in $G$ to $a_{\ga}b_{\ga}\in A^-B$.  We want to know which roots of unity can occur as eigenvalues of $\Ad(\ga)$.  These are the same as the roots of unity which occur as eigenvalues of $\Ad(a_{\ga}b_{\ga})$.  If
$$
b_{\ga}=\matrixtwo{R(\th)}{R(\phi)},
$$
where
$$
R(\th)=\matrix{\cos\th}{-\sin\th}{\sin\th}{\cos\th},\index{$R(\th)$}
$$
then the eigenvalues of $\Ad(a_{\ga}b_{\ga})$ are $e^{\pm2i\th}$ and $e^{\pm2i\phi}$.  Define the sets
$$
R_{\th}=\{n\in\N:n\th\in\pi\Z\}\ \textrm{ and }\ R_{\phi}=\{n\in\N:n\phi\in\pi\Z\}
$$
and
$$
R_{\ga}=\{\min R_{\th},\min R_{\phi}\}.\index{$R_{\ga}$}
$$
Then $R_{\ga}$ contains either $0,1$ or $2$ elements.  We can see that $\ga$ weakly neat is equivalent to $R_{\ga}=\varnothing$.  For $\ga\in\CE_P(\Ga)$, where $\ga=\ga_0^n$ for $\ga_0$ primitive, the value of $\chi_1(\Ga_{\ga})$ depends on whether or not $n\in R_{\ga_0}$.

For $I\subset R_{\ga}$ with $I\neq\varnothing$ define $n_I$ \index{$n_I$} to be the least common multiple of the elements of $I$ and set $n_{\varnothing}=1$.  Further, define
$$
\chi_I(\ga)=\frac{(-1)^{|I|}}{n_I}\sum_{J\subset I}(-1)^{|J|}\chi_1(\Ga_{\ga^{n_J}}).\index{$\chi_I(\ga)$}
$$

Let $z\in\C\smallsetminus\{0\}$ and $q\in\Q$.  We define $z^q$ to be equal to $e^{q\log z}$, where we take the branch of the logarithm with imaginary part in the interval $(-\pi,\pi]$.

For any finite-dimensional virtual representation $\si$ of $M$ we define, for $\Re(s)$ large, the \emph{generalised Selberg zeta function}\index{generalised Selberg zeta function}\index{Selberg zeta function, generalised}
\begin{equation}
\label{eqn:GenSelb}
\index{$Z_{P,\si}(s)$}
Z_{P,\si}(s)= \prod_{[\ga]\in\CE_P^p(\Ga)}\prod_{n\geq 0}\prod_{I\subset R_{\ga}}\det\left(1-e^{-sn_Il_{\ga}}\ga^{n_I}\,|V_{\si}\ox S^n(\n)\right)^{\chi_I(\ga)},
\end{equation}
where $S^n (\n)$ denotes the $n^{{\rm th}}$ symmetric power of the space $\n$ and $\ga$ acts on $V_{\si}\ox S^n(\n)$ via $\si(b_{\ga})\ox \Ad^n(a_{\ga}b_{\ga})$.  In the case that $\Ga$ is weakly neat this simplifies to
$$
Z_{P,\si}(s)= \prod_{[\ga]\in\CE_P^p(\Ga)}\prod_{n\geq 0}\det\left(1-e^{-sl_{\ga}}\ga\,|V_{\si}\ox S^n(\n)\right)^{\chi_1(\Ga_{\ga})},
$$
(see \cite{geom}).

Let $\pi\in\hat G$.
The Lie algebra $\n$ acts on $\pi_K$ and we denote by $H^{\bullet}(\n,\pi_K)$ \index{$H^{\bullet}(\n,\pi_K)$} (resp. $H_{\bullet}(\n,\pi_K)$) the corresponding Lie algebra cohomology (resp. homology) (see \cite{BorelWallach80},\cite{CartanEilenberg56}).  We have the following isomorphism of $AM$-modules (see \cite{HechtSchmid83}, p57):
\begin{equation}
\label{eqn:AM-iso}
H_p(\n,\pi_K)\cong H^{4-p}(\n,\pi_K)\otimes \bwedge^4 \n.
\end{equation}

For $\la\in\a^*$ let
\begin{equation}
\label{eqn:MLamdaPi}
\index{$m_{\la}(\pi)$}
m_{\la}(\pi)=\sum_{p,q\geq 0} (-1)^{p+q} \dim\left( H^q(\n, \pi_K)^{\la}\otimes\bwedge^p \p_M \otimes V_{\si}\right)^{K_M},
\end{equation}
where for an $\a$-module V and $\la\in\a^*$, $V^{\la}$ \index{$V^{\la}$} denotes the generalised $\la$-eigenspace 
$$
\{v\in V|\ \exists n\in\N\ {\rm such\ that}\ (a-\la(a))^n v = 0\ \forall a\in\a\}
$$ 
and the superscript $K_M$ denotes the subspace of $K_M$ invariants.

We say that an admissible representation $\pi$ of a linear connected reductive group $G'$ has a \emph{global character} \index{global character} $\Th=\Th^{G'}_{\pi}$ \index{$\Th^G_{\pi}$} if for all smooth, compactly supported functions $f$ on $G'$ the operator $\pi(f)$ is of trace class and $\Th^{G'}_{\pi}$ is a locally $L^2$ function on $G'$ satisfying
$$
\tr\pi(f)=\int_{G'}\Th^{G'}_{\pi}(g)f(g)\,dg
$$
for each such $f$.  By \cite{Knapp}, Theorem 10.2, every irreducible admissible, and in particular every irreducible unitary, representation of $G'$ has a global character.

\begin{theorem}
\label{thm:SelbergZeta}
Let $\si$ be a finite dimensional virtual representation of $M$.  The function $Z_{P,\si}(s)$ extends to a meromorphic function on the whole of $\C$.  For $\la\in\a^*$, the vanishing order of $Z_{P,\si}(s)$ at the point $s=\la(H_1)$ is
\begin{equation}
\label{eqn:v-order}
\sum_{\pi\in\hat{G}}N_{\Ga}(\pi)m_{\la}(\pi).
\end{equation}
Further, all the poles and zeros of $Z_{P,\si}(s)$ lie in $\R\cup(\frac{1}{2}+\textit{i}\R)$.
\end{theorem}
\prf
The analogue of this theorem in the case that $G$ has trivial centre and $\Ga$ is weakly neat is proved in \cite{geom}, Theorem 2.1.  If $G$ has trivial centre then $\Ga$ weakly neat implies $\Ga$ torsion free, since non weakly neat torsion elements must be central.  With a few modfications the proof of \cite{geom}, Theorem 2.1 becomes valid in our case also.  In fact the assumption that $\Ga$ is weakly neat is used in two places.  We sketch the proof here, pointing out the necessary modifications for it to be valid in our case.

The theorem is proved by setting $f=f_{\si}$ in the test function $\Phi=\Phi_{j,s}=\Phi_{\si,j,s}$ defined in (\ref{eqn:TestFunction}), where $f_{\si}$ is an Euler-Poincar$\eac$ function for $\si$ on $M$.  We then get that for $j\in\N$ and $\Re(s)$ sufficiently large, the right hand side of (\ref{eqn:STF2}) is equal to
\begin{equation}
\label{eqn:STFRHS}
\sum_{[\ga]\in\CE_P(\Ga)}\vol(\Ga_{\ga}\bs G_{\ga})\CO^M_{b_{\ga}}(f_{\si})\frac{l_{\ga}^{j+1}e^{-sl_{\ga}}}{\det(1-(a_{\ga}b_{\ga})^{-1}|\bar{\n})}.
\end{equation}

In \cite{geom} the following equation from \cite{Deitmar95}, proved under the assumption that $\Ga$ is torsion free, is used:
\begin{equation}
\label{eqn:ECharVol}
\chi_1(\Ga_{\ga})=\frac{C_{\ga}\vol(\Ga_{\ga}\bs G_{\ga})}{l_{\ga_0}},
\end{equation}
where $C_{\ga}$ is an explicit constant depending only on $\ga$, and $\ga_0$ denotes the primitive geodesic underlying $\ga$.  We have shown in Lemma \ref{lem:ECharEqn} that (\ref{eqn:ECharVol}) holds for all $\ga\in\CE_P(\Ga)$ in our case also.  From \cite{geom}, Proposition 1.4 we take the equation
$$
\CO^M_{b_{\ga}}(f_{\si})=C_{\ga}\tr\si(b_{\ga}),
$$
which, together with (\ref{eqn:ECharVol}) gives us
$$
\CO^M_{b_{\ga}}(f_{\si})\vol(\Ga_{\ga}\bs G_{\ga})=l_{\ga_0}\tr\si(b_{\ga})\chi_1(\Ga_{\ga}).
$$
Substituting this into (\ref{eqn:STFRHS}) we get
$$
\sum_{[\ga]\in\CE_P(\Ga)}l_{\ga_0}\tr\si(b_{\ga})\chi_1(\Ga_{\ga})\frac{l_{\ga}^{j+1}e^{-sl_{\ga}}}{\det(1-(a_{\ga}b_{\ga})^{-1}|\bar{\n})}.
$$
We claim that this is equal to
$$
(-1)^{j+1}\left(\frac{\partial\ }{\partial s}\right)^{j+2}\log Z_{P,\si}(s).
$$
Indeed
$$
\log Z_{P,\si}(s)=-\sum_{[\ga]\in\CE_P^p(\Ga)}\sum_{m=1}^\infty\sum_{I\subset R_{\ga}}\chi_I(\ga)\frac{e^{-sml_{\ga}}}{m}\tr\si(b_{\ga}^m)\sum_{n\geq 0}\tr\left((a_{\ga}b_{\ga})^m|S^n(\n)\right).
$$
We note that
\begin{eqnarray*}
\sum_{n\geq 0}\tr\left((a_{\ga}b_{\ga})^m|S^n(\n)\right) & = & \det\left(1-(a_{\ga}b_{\ga})^m|\n\right)^{-1} \\
  & = & \det\left(1-(a_{\ga}b_{\ga})^{-m}|\bar{\n}\right)^{-1}
\end{eqnarray*}
and the claimed equality follows.

In \cite{geom}, since $\Ga$ is weakly neat it follows that $G_{\ga^n}=G_{\ga}$ for all $\ga\in\Ga$ and $n\in\N$.  In \cite{Duistermaat79} it is shown that $X_{\ga}\cong\Ga_{\ga}\bs G_{\ga}/K_{\ga}$ and so it follows that for all $n\in\N$ we have $X_{\ga^n}\cong X_{\ga}$ and hence $\chi_1(\Ga_{\ga^n})=\chi_1(\Ga_{\ga})$.  Thus the above equality involving the logarithmic derivative of $Z_{P,\si}(s)$ is derived with the simpler Euler product expansion for $Z_{P,\si}(s)$ given above.

In our case we may have an element $[\ga]\in\CE_P^p(\Ga)$ with a non-trivial root of unity as an eigenvalue.  For such a $\ga$ we have $\chi_1(\Ga_{\ga^n})\neq\chi_1(\Ga_{\ga})$ for $n\in R_{\ga}$.  For this reason we have had to introduce the more complicated Euler product expansion for $Z_{P,\si}(s)$ so that the above equality involving the logarithmic derivative of $Z_{P,\si}(s)$ still holds.

On page 909 of \cite{geom} it is shown that
\begin{equation}
\label{eqn:TracePiPhi1}
\tr\pi(\Phi_s) = \int_{MA^-}f_{\si}(m)\Th^{MA}_{H^{\bullet}(\n,\pi_K)}(ma)dm\,g_s^j(a)da.
\end{equation}
Using the property (\ref{eqn:EPFnId}) of $f_{\si}$ we get
\begin{eqnarray}
\label{eqn:TracePiPhi}
\tr\pi(\Phi_s) & = & \int_{A^-}\tr\left(a\left|\left(H^{\bullet}(\n,\pi_K)\ox\bwedge^*\p_M\ox V_{\si}\right)^{K_M}\right.\right)l_a^{j+1}e^{-sl_a}\,da \nonumber \\
  & = & \int_0^{\infty}\sum_{\la\in\a^*}m_{\la}(\pi)\,e^{(\la(H_1)-s)t}\,t^{j+1}\,dt \nonumber \\
  & = & (-1)^{j+1}\left(\frac{\partial\ }{\partial s}\right)^{j+1}\sum_{\la\in\a^*}m_{\la}(\pi)\rez{s-\la(H_1)}.
\end{eqnarray}
Proposition \ref{pro:STF2} then gives us that
\begin{equation}
\label{eqn:HighLogDer}
\left(\frac{\partial\ }{\partial s}\right)^{j+2}\log Z_{P,\si}(s)=\sum_{\pi\in\hat{G}}N_{\Ga}(\pi)\left(\frac{\partial\ }{\partial s}\right)^{j+1}\sum_{\la\in\a^*}m_{\la}(\pi)\rez{s-\la(H_1)},
\end{equation}
from which it follows that the vanishing-order of $Z_{P,\si}(s)$ at $s=\la(H_1)$ is
$$
\sum_{\pi\in\hat{G}}N_\Ga(\pi) m_{\la}(\pi).
$$

Two further comments are in order.  First, in \cite{geom} the positive Weyl chamber $A^+=\{\exp(tH_1):t<0\}$ in $A$ is considered, where we have instead considered the negative chamber $A^-$.  For this reason we have interchanged the positions of the Lie algebras $\n$ and $\bar{\n}$ from the way they are used in \cite{geom}.  This is easily seen to give a precisely equivalent result.

Secondly, the results of \cite{geom} are stated for a finite dimensional representation $\si$ of $M$.  By linearity the results extend in a straightforward way to the case where $\si$ is a virtual representation, which we use here.
\qed

\subsection{A functional equation for $Z_{P,\si}(s)$}

\begin{proposition}
\label{pro:Bounds}
For $\la\in\a^*$ let $\|\la\|$ be the norm given by the form $b$ in \textnormal{(\ref{eqn:norm})}.  There are $m_1\in\N$ and $C>0$ such that for every $\pi\in\hat{G}$ and every $\la\in\a^*$ we have
$$
\left|\sum_{q=0}^4(-1)^q\dim(H^q(\n,\pi_K)^{\la})\right|\leq C(1+\|\la\|)^{m_1}.
$$
Further, let $S$ denote the setof all pairs $(\pi,\la)\in\hat{G}\x\a^*$ such that
$$
\sum_{q=0}^4(-1)^q\dim\left(H^q(\n,\pi_K)^{\la}\right)\neq 0.
$$
Then there is $m_2\in\N$ such that
$$
\sum_{(\pi,\la)\in S}\frac{N_{\Ga}(\pi)}{(1+\|\la\|)^{m_2}}<\infty.
$$
\end{proposition}
\prf
These results follow from \cite{primgeo}, Proposition 2.4 and Lemma 2.7.
\qed

An entire function $f$ is said to be of \emph{finite order} \index{finite order, function of} if there is a positive constant $a$ and an $r>0$ such that $|f(z)|<\exp\left(|z|^a\right)$ for $|z|>r$.  If $f$ is of finite order then the \emph{order of $f$} \index{order!of a function} is defined to be the infimum of such $a$'s.

It is well known that the classical Selberg zeta function is an entire function of order two (see \cite{Selberg56}).  Our next lemma gives an analogous result for the generalised Selberg zeta function we are considering here.

Let $\si$ be a finite dimensional virtual representation of $M$.  Theorem \ref{thm:SelbergZeta} tells us that $Z_{P,\si}(s)$ is meromorphic and hence it may be written as the quotient of two entire functions:
$$
Z_{P,\si}(s)=\frac{Z_1(s)}{Z_2(s)},
$$
where the zeros of $Z_1(s)$ correspond to the zeros of $Z_{P,\si}(S)$ and the zeros of $Z_2(s)$ correspond to the poles of $Z_{P,\si}(s)$.  The orders of the zeros of $Z_1(s)$ (resp. $Z_2(s)$) equal the orders of the corresponding zeros (resp. poles) of $Z_{P,\si}(s)$.  The functions $Z_1(s)$ and $Z_2(s)$ are not uniquely determined, but clearly their respective sets of zeros, together with the orders of the zeros, are.  For $i=1,2$ let $R_i$ denote the set of zeros of $Z_i(s)$ counted with multiplicity.

\begin{lemma}
\label{lem:FiniteOrder}
There exist $Z_i(s)$, for $i=1,2$, with the above properties, which are, in addition, of finite order.
\end{lemma}

\prf
Let $Z_1(s)$, $Z_2(s)$ be as above.  We shall show that we may take them to be of finite order.

For $\pi\in\hat{G}$ let $\La(\pi)$ be the set of all $\la\in\a^*$ with $\la\neq 0$ and $m_{\la}(\pi)\neq 0$.  Let $\hat{G}(\Ga)$ be the set of $\pi\in\hat{G}$ such that $N_{\Ga}(\pi)\neq 0$ and let $S$ denote the set of all pairs $(\pi,\la)$ such that $\pi\in\hat{G}(\Ga)$ and $\la\in\La(\pi)$.  For $\la\in\a^*$ let $\|\la\|$ be the norm given by the form $b$ in (\ref{eqn:norm}).

The expression (\ref{eqn:v-order}) tells us that $s\neq 0$ is a zero or pole of $Z_P(s)$ if and only if $s=\la(H_1)$ for some $\la\in\a^*$, for which there exists $\pi\in\hat{G}(\Ga)$ such that $(\pi,\la)\in S$.

Since the function $Z_{P,\si}(s)$ is meromorphic and non-zero, it follows that there exists $c>0$ such that
\begin{equation}
\label{eqn:Bound1}
|\la(H_1)|\geq c\left(1+\|\la\|\right)
\end{equation}
for all $\la$ such that $(\pi,\la)\in S$ for some $\pi$.

By the definition of $m_{\la}(\pi)$ we deduce immediately from Proposition \ref{pro:Bounds} that there exist $m_1\in\N$ and $C>0$ such that for every $\pi\in\hat{G}$ and every $\la\in\a^*$ we have
\begin{equation}
\label{eqn:Bound2}
|m_{\la}(\pi)|\leq C\left(1+\|\la\|\right)^{m_1},
\end{equation}
and that there exists $m_2\in\N$ such that
\begin{equation}
\label{eqn:Bound3}
\sum_{(\pi,\la)\in S} \frac{N_{\Ga}(\pi)}{\left(1+\|\la\|\right)^{m_2}}<\infty.
\end{equation}
Let $k = m_1 + m_2$.  Then, for $(\pi,\la)\in S$, by (\ref{eqn:Bound1}) and (\ref{eqn:Bound2}) we have,
\begin{eqnarray*}
\left|\frac{m_{\la}(\pi)}{\la(H_1)^k}\right| & \leq & \rez{c^k}\cdot\frac{|m_{\la}(\pi)|}{\left(1+\|\la\|\right)^k} \\
  & \leq & \frac{C}{c^k}\cdot\rez{\left(1+\|\la\|\right)^{m_2}}.
\end{eqnarray*}
It then follows from (\ref{eqn:v-order}) and (\ref{eqn:Bound3}) that, for $i=1,2$,
\begin{equation}
\label{eqn:FiniteRank}
\sum_{\rho\in R_i\smallsetminus\{0\}}\rez{|\rho|^k}\leq\sum_{(\pi,\la)\in S}\frac{N_{\Ga}(\pi)|m_{\la}(\pi)|}{|\la(H_1)|^k}<\infty.
\end{equation}
We say that $Z_1(s)$ and $Z_2(s)$ are of \emph{finite rank}\index{finite rank, function of}.

By the Weierstrass Factorisation Theorem \index{Weierstrass Factorisation Theorem} (\cite{Conway78}, VII.5.14) and (\ref{eqn:FiniteRank}), there exist entire functions $g_1(s)$, $g_2(s)$ such that
\begin{equation}
\label{eqn:Weierstrass}
Z_i(s)=s^{n_i}e^{g_i(s)}\prod_{\rho\in R_i\smallsetminus\{0\}}\left(1-\frac{s}{\rho}\right)\exp\left(\frac{s}{\rho}+\frac{s^2}{2\rho^2}\cdots+\frac{s^k}{k\rho^k}\right),
\end{equation}
where $n_i$ is the order of the zero of $Z_i(s)$ at $s=0$.

From (\ref{eqn:HighLogDer}) we recall the equation
\begin{equation}
\label{eqn:LogDerEqn}
\left(\frac{d\ }{ds}\right)^j\log Z_{P,\si}(s)=\sum_{\pi\in\hat{G}(\Ga)}N_{\Ga}(\pi)\left(\frac{d\ }{ds}\right)^{j-1}\sum_{\la\in\La(\pi)\cup\{0\}}\frac{m_{\la}(\pi)}{s-\la(H_1)},
\end{equation}
where $j\in\N$ is sufficiently large.  Let $J=\max(j,k)$.  It follows from (\ref{eqn:Weierstrass}) and (\ref{eqn:LogDerEqn}) that
$$
\left(\frac{d\ }{ds}\right)^{J-1}\left(g'_1(s)-g'_2(s)\right) + \sum_{\rho\in R_1}\frac{(-1)^{J-1}(J-1)!}{(s-\rho)^J} - \sum_{\rho\in R_2}\frac{(-1)^{J-1}(J-1)!}{(s-\rho)^J}
$$
is equal to
$$
\sum_{\pi\in\hat{G}(\Ga)}\sum_{\la\in\La(\pi)\cup\{0\}}N_{\Ga}(\pi)m_{\la}(\pi)\frac{(-1)^{J-1}(J-1)!}{(s-\la(H_1))^J}.
$$

Remembering that the zeros of $Z_i(s)$ are included with multiplicity in $R_i(s)$ and bearing in mind the expression (\ref{eqn:v-order}) for the vanishing order of $Z_{P,\si}(s)$, we see that this implies that
$$
\left(\frac{d\ }{ds}\right)^{J-1}\left(g'_1(s)-g'_2(s)\right)=0.
$$
It follows that $g_1(s)-g_2(s)$ is a polynomial of degree at most $J$.  Without loss of generality we may take $g_2(s)$ to be zero, so that $g_1(s)$ is itself a polynomial of degree at most $J$.  Finally, Theorem XI.2.6 of \cite{Conway78} tells us that since $Z_1(s)$ and $Z_2(s)$ are of finite rank and $g_1(s)$ and $g_2(s)$ are both polynomials of degree at most $J$, it follows that $Z_1(s)$ and $Z_2(s)$ are both of order at most $J$.
\qed

Before we give the next lemma we make a couple of definitions.  Let $G'$ be a linear connected reductive group with maximal compact subgroup $K'$ and let $\xi$ be a representation of $G'$.  Let $\g'$ be the complexified Lie algebra of $G'$, let $\h'\subset\g'$ be the diagonal subalgebra and let $\rho'\in(\h')^*$ be the half sum of the positive roots of the system $(\g',\h')$.  We say that $\xi$ is \emph{tempered} \index{representation!tempered} if for all $K'$-finite vectors $u,v\in \xi_{K'}$ there exist constants $c_{u,v}$ such that for all $g\in G'$
$$
|\langle\xi(g)u,v\rangle|\leq c_{u,v}\int_{K'}e^{-\rho'(H(g^{-1}k))}dk,
$$
where $H(g^{-1}k)$ denotes the logarithm of the split part of $g^{-1}k$ under the Iwasawa decomposition and $\langle\cdot,\cdot\rangle$ is the inner product on $V_{\xi}$.  An admissible representation of $G'$ is called \emph{standard} \index{representation!standard} if it is induced from an irreducible tempered representation of $M'\subset G'$, where $P'=M'A'N'$ is a parabolic subgroup of $G'$.

The Weyl group $W(G,A)$ has two elements, let $w$ be the nontrivial element therein.  Then $w$ acts on $M$ by conjugation and for a representation $\si$ of $M$ we can define the representation $^w\si$ by $^w\si(m)=\si(wmw^{-1})$.

\begin{lemma}
\label{lem:Zeros}
Suppose that $\si\cong\,^w\!\si$ as $K_M$-modules.  Then the functions $Z_{P,\si}(s)$ and $Z_{P,\si}(1-s)$ have the same poles and zeros with multiplicity.
\end{lemma}
\prf
By Theorem \ref{thm:SelbergZeta} $\la\in\a^*$, the vanishing order of $Z_{P,\si}(s)$ at the point $s=\la(H_1)$ is equal to $\sum_{\pi\in\hat{G}}N_{\Ga}(\pi)m_{\la}(\pi)$.  We shall show that
\begin{equation}
\label{eqn:MLaPiEq}
m_{\la}(\pi)=m_{-\la-2\rho_P}(\pi)
\end{equation}
for all $\pi\in\hat{G}$.  Since $-2\rho_P(H_1)=1$ the lemma will follow.

Let $\pi\in\hat{G}$.  Recall from (\ref{eqn:MLamdaPi}) that
$$
m_{\la}(\pi)=\sum_{p,q\geq 0} (-1)^{p+q} \dim\left( H^q(\n, \pi_K)^{\la}\otimes\bwedge^p \p_M \otimes V_{\si}\right)^{K_M}.
$$
The global character $\Th^G_{\pi}$ of $\pi$ on $G$ can be written as a linear combination with integer coefficients of characters of standard representations (\cite{Knapp}, Chapter X, \S10.2).  We are interested in the values taken by $\Th^G_{\pi}$ on $MA^-$.  By \cite{Knapp}, Proposition 10.19, the only characters which are non-zero on $MA$ are those characters of representations induced from $P=MAN$.  So there exist $n\in\N$ and for all $1\leq i\leq n$ integers $c_i$, tempered representations $\xi_i$ of $M$ and $\nu_i\in\a^*$ such that
\begin{equation}
\label{eqn:TempRepSum}
\Th^G_{\pi}=\sum_{i=1}^n c_i\Th^G_{\pi^i},
\end{equation}
where $\pi^i=\Ind_P^G(\xi_i\ox\nu_i)$.  From \cite{HechtSchmid83}, Theorem 3.6 and (\ref{eqn:AM-iso}) it follows that for all regular $ma\in MA^-$ we have
$$
\Th^G_{\pi}(ma)=\Th^{MA}_{H^{\bullet}(\n,\pi_K)}(ma)\frac{\det\left(a|\bwedge^4\n\right)}{\det(1-ma|\n)},
$$
and the same holds for all $\pi^i$.  Together with (\ref{eqn:TempRepSum}) this implies that
$$
\Th^{MA}_{H^{\bullet}(\n,\pi_K)}(ma)=\sum_{i=1}^n\Th^{MA}_{H^{\bullet}(\n,\pi^i_K)}(ma)
$$
for all regular $ma\in MA^-$.  Substituting this into (\ref{eqn:TracePiPhi1}) and proceeding as in (\ref{eqn:TracePiPhi}) we see that it suffices to show (\ref{eqn:MLaPiEq}) for the representations $\pi^i$.

The Weyl group element $w$ acts on the group $MA$ by conjugation, which has the effect of swapping the two components.  For $\nu\in\a^*$ we therefore have $w\nu=-\nu$.  Recall that for a representation $\xi$ of $M$ we let $^w\xi$ denote the representation $^w\xi(m)=\xi(wmw^{-1})$.  By \cite{Knapp}, Theorem 8.38, the representations $\pi^i$ and $^w\pi^i=\Ind_P^G(^w\xi_i\ox(-\nu_i))$ are equivalent.  From the definition of the induced group action we can see that
$$
H^q(\n,\pi^i_K)^{\la}=H^q(\n,^w\!\pi^i_K)^{-\la-2\rho_P}.
$$
Thus we see that
$$
\dim\left( H^q(\n, \pi^i_K)^{\la}\otimes\bwedge^p \p_M \otimes V_{\si}\right)^{K_M}
$$
is equal to
$$
\dim\left( H^q(\n, ^w\!\pi^i_K)^{-\la-2\rho_P}\otimes\bwedge^p \p_M \otimes V_{\si}\right)^{K_M}.
$$

Using the notation of Proposition \ref{pro:KMDual} we note that the $K_M$-types satisfy $^wtriv_{K_M}=triv_{K_M}$, $^wdet_{K_M}=det_{K_M}$ and $^w\d_{l,m}=\d_{m,l}\cong\d_{-m,-l}$.  Lemma \ref{lem:pMKMTypes1} and the isomorphism of $K_M$-modules
$$
\n\cong\d_{2,2}\oplus\d_{2,-2}
$$
tell us that the $K_M$-modules $\bwedge^p\p_M$ and $\n$ are invariant under the action of $w$ and we have assumed that $^w\si\cong\si$ as $K_M$-modules.  Since $w^2$ is the identity element of $W(G,A)$ we can conclude that
$$
\dim\left( H^q(\n, \pi^i_K)^{\la}\otimes\bwedge^p \p_M \otimes V_{\si}\right)^{K_M}
$$
is equal to
$$
\dim\left( H^q(\n,\pi^i_K)^{-\la-2\rho_P}\otimes\bwedge^p \p_M \otimes V_{\si}\right)^{K_M}.
$$
Hence $m_{\la}(\pi^i)=m_{-\la-2\rho_P}(\pi^i)$ and the lemma follows.
\qed

\begin{theorem}
\label{thm:FunctEqn}
Suppose that $\si\cong\,^w\!\si$ as $K_M$-modules.  Then there exists a polynomial $G(s)$ such that $Z_{P,\si}(s)$ satisfies the functional equation
$$
Z_{P,\si}(s)=e^{G(s)}Z_{P,\si}(1-s).
$$
\end{theorem}
\prf
Let $R_i$, $n_i$ be as above for $i=1,2$.  We saw in the proof of Lemma \ref{lem:FiniteOrder} that there exist entire functions $Z_1(s)$ and $Z_2(s)$ of finite order such that
$$
Z_{P,\si}(s)=\frac{Z_1(s)}{Z_2(s)}
$$
and polynomials $g_1(s)$ and $g_2(s)$ such that
$$
Z_i(s)=s^{n_i}e^{g_i(s)}\prod_{\rho\in R_i\smallsetminus\{0\}}\left(1-\frac{s}{\rho}\right)\exp\left(\frac{s}{\rho}+\frac{s^2}{2\rho^2}\cdots+\frac{s^k}{k\rho^k}\right).
$$
Lemma \ref{lem:Zeros} tells us that $Z_{P,\si}(s)$ and $Z_{P,\si}(1-s)$ have the same poles and zeros.  Hence we can in the same way conclude that there exist entire functions $W_1(s)$ and $W_2(s)$ of finite order such that
$$
Z_{P,\si}(1-s)=\frac{W_1(s)}{W_2(s)}
$$
and polynomials $h_1(s)$ and $h_2(s)$ such that
$$
W_i(s)=s^{n_i}e^{h_i(s)}\prod_{\rho\in R_i\smallsetminus\{0\}}\left(1-\frac{s}{\rho}\right)\exp\left(\frac{s}{\rho}+\frac{s^2}{2\rho^2}\cdots+\frac{s^k}{k\rho^k}\right).
$$
Setting
$$
G(s)=g_1(s)+h_2(s)-g_2(s)-h_1(s)
$$
we can see that the claimed functional equation does indeed hold.
\qed

\subsection{The Ruelle zeta function}
\label{sec:RuelleZetaFn}
For any finite-dimensional virtual representation $\si$ of $M$ we define, for $\Re(s)$ large, the \emph{generalised Ruelle zeta function}\index{generalised Ruelle zeta function}\index{Ruelle zeta function, generalised}
$$
R_{\Ga,\si}(s) = \prod_{[\ga]\in\CE_P^p(\Ga)}\prod_{I\subset R_{\ga}}\det\left(1-e^{-sn_Il_{\ga}}\ga^{n_I}\,|V_{\si}\right)^{\chi_I(\ga)}.\index{$R_{\Ga,\si}(s)$}
$$
We have the following theorem giving a relationship between the generalised Selberg zeta function and the generalised Ruelle zeta function.

\begin{theorem}
\label{thm:RuelleZeta}
Let $\si$ be a finite dimensional virtual representation of $M$.  The function $R_{\Ga,\si}(s)$ extends to a meromorphic function on $\C$.  More precisely
$$
R_{\Ga,\si}(s) = \prod_{q=0}^4 Z_{P,(\bwedge^q\bar{\n}\ox\si)}\left( s+\frac{q}{4}\right)^{(-1)^q}.
$$
\end{theorem}
\prf
Let $\ga\in\CE_P(\Ga)$ let $\mu(\ga)\in\N$ be the least such that $\ga=\ga_0^{\mu(\ga)}$ for some $\ga_0\in\CE_P^p(\Ga)$.  We compute at first
\begin{eqnarray*}
\log R_{\Ga,\si}(s) & = & \sum_{[\ga]\in\CE_P^p(\Ga)}\sum_{I\subset R_{\ga}}\chi_I(\ga)\tr\log\left(1-e^{-sn_Il_{\ga}}\ga^{n_I}\,|V_{\si}\right) \\
  & = & \sum_{[\ga]\in\CE_P^p(\Ga)}\sum_{I\subset R_{\ga}}\chi_I(\ga)\sum_{m=1}^{\infty}\rez{m}e^{-smn_Il_{\ga}}\tr\si(b_{\ga}^{n_I}) \\
  & = & \sum_{[\ga]\in\CE_P(\Ga)}\sum_{I\subset R_{\ga}}\chi_I(\ga)\frac{e^{-sn_Il_{\ga}}}{\mu(\ga)}\tr\si(b_{\ga}^{n_I}).
\end{eqnarray*}
Similarly, for $\tau_q=\bwedge^q\bar{\n}\ox \si$ we have
that $\log Z_{P,\tau_q}(s)$ equals
\begin{eqnarray*}
 &  & \sum_{[\ga]\in\CE^p_P(\Ga)}\sum_{I\subset R_{\ga}}\chi_I(\ga)\sum_{n=0}^{\infty}\tr\log\left((1-e^{-sn_Il_{\ga}}\ga^{n_I}|V_{\tau_q}\ox S^n(\n)\right) \\
  & = & \sum_{[\ga]\in\CE^p_P(\Ga)}\sum_{I\subset R_{\ga}}\chi_I(\ga)\sum_{n=0}^{\infty}\sum_{m=1}^{\infty}\frac{-1}{m}e^{-smn_Il_{\ga}}\tr\left(\ga^{mn_I}|V_{\tau_q}\ox S^n(\n)\right) \\
  & = & -\sum_{[\ga]\in\CE_P(\Ga)}\sum_{I\subset R_{\ga}}\chi_I(\ga)\frac{e^{-sn_Il_{\ga}}}{\mu(\ga)}\tr\si(b_{\ga}^{n_I})\frac{\tr(b_{\ga}^{n_I}|\bwedge^q\bar{\n})}{\det\left(1-(a_{\ga}b_{\ga})^{-n_I}|\bar{\n}\right)} \\\
  & = & -\sum_{[\ga]\in\CE_P(\Ga)}\sum_{I\subset R_{\ga}}\chi_I(\ga)\frac{e^{-sn_Il_{\ga}}}{\mu(\ga)}\tr\si(b_{\ga}^{n_I})\frac{\tr(b_{\ga}^{n_I}|\bwedge^q\bar{\n})}{\tr\left((a_{\ga}b_{\ga})^{-n_I}|\bwedge^*\bar{\n}\right)}.
\end{eqnarray*}
Since $\n$ is an $M$-module defined over the reals we conclude that the trace $\tr\left((a_{\ga}b_{\ga})^{-1}|\bwedge^*\bar{\n}\right)$ is a real number.  Therefore it equals its complex conjugate, which is $\tr\left(a_{\ga}^{-1}b_{\ga}|\bwedge^*\bar{\n}\right)$.  Summing over $q$ we get
$$
\log R_{\Ga,\si}(s)=\sum_{q=0}^4(-1)^q\log Z_{P,\tau_q}\left(s+\frac{q}{4}\right).
$$
The theorem follows.
\qed

We shall be interested in the zeta function $R_{P,\si}(s)$ in the case where $\si=triv=1$ is the trivial representation of $M$ and in the case where $\si$ is the following virtual representation of $M$.  Define
$$
\tilde{\si}^+=15\bwedge^0\m+6\bwedge^2\m+\bwedge^4\m
$$
and
$$
\tilde{\si}^-=10\bwedge^1\m+3\bwedge^3\m
$$
and let $\tilde{\si}=\tilde{\si}^+ -\tilde{\si}^-$\index{$\tilde{\si}$}, where $M$ acts on $\bwedge^n\m$ according to $\Ad^n$.  The reason for this choice of $\tilde{\si}$ is that it allows us to separate the contribution of the non-regular elements of $\CE_P(\Ga)$ from the regular ones, as made clear in the following lemma.

\begin{lemma}
\label{lem:SigmaTilde}
Let $\ga\in\CE_P(\Ga)$ with
$$
b_{\ga}=\matrixtwo{R(\th)}{R(\phi)},
$$
where
$$
R(\th)=\matrix{\cos\th}{-\sin\th}{\sin\th}{\cos\th}.
$$
Then
$$
\tr\tilde{\si}(b_{\ga})=\det(1-b_{\ga}:\m/\b)=4(1-\cos 2\th)(1-\cos 2\phi),
$$
where $\b$ is the complexified Lie algebra of the group $B$.  In particular we have $\tr\tilde{\si}(b_{\ga})\geq 0$ for all $\ga\in\CE_P(\Ga)$, and $\tr\tilde{\si}(b_{\ga})=0$ if and only if $\ga$ is non-regular.
\end{lemma}
\prf
As $B$-modules, we have the isomorphism $\m\cong\b\oplus(\m/\b)$.  The group $B$ is abelian and so acts trivially on the 2-dimensional Lie algebra $\b$.  Hence, for $0\leq n\leq 6$, we have the following isomorphism of $B$-modules
\begin{equation}
\label{eqn:BModIso}
\bwedge^n\m\cong\bigoplus_{p+q=n}\left({2}\atop{p}\right)\bwedge^q\m/\b.
\end{equation}
We set $a_0=1$, $a_1=-3$, $a_2=6$, $a_3=-10$, $a_4=15$ and note that, for $k=0,...,4$, these satisfy
\begin{equation}
\label{eqn:ASum}
\sum_{m=0}^{k}a_{k-m}\left({2}\atop{m}\right)=(-1)^k.
\end{equation}
The $B$-module isomorphism
$$
\bwedge^*(\m/\b)\cong\sum_{n=0}^4 a_{4-n}\bwedge^n\m=\tilde{\si}
$$
then follows from (\ref{eqn:BModIso}) and (\ref{eqn:ASum}).  We can then see that
$$
\tr\tilde{\si}(b_{\ga})=\tr\left(b_{\ga}:\bwedge^*\m/\b\right)=\det(1-b_{\ga}:\m/\b).
$$
The adjoint action of $B$ on $\m/\b$ can easily be computed to give the claimed value.  Finally we note that $\tr\tilde{\si}(b_{\ga})=0$ if and only if $\th,\phi\in\Z\pi$ which is equivalent to $\ga$ being non-regular.
\qed

The main result of this section is the following theorem.

\begin{theorem}
\label{thm:RuelleMain}
The function $R_{\Ga,1}(s)$ has a double zero at $s=1$.  The function $R_{\Ga,\tilde{\si}}(s)$ has a zero of order eight at $s=1$.  Apart from that, for $\si\in\{1,\tilde{\si}\}$, all poles and zeros of $R_{\Ga,\si}(s)$ are contained in the union of the interval $\left[-1,\frac{3}{4}\right]$ with the five vertical lines given by $\frac{k}{4} + i\R$ for $k=-2,-1,0,1,2$.
\end{theorem}

The theorem will follow from the following proposition, which we prove in the remainder of the section.

\begin{proposition}
\label{pro:SelbergPolesZeros}
The function $Z_{P,1}(s)$ has a double zero at $s=1$.  The function $Z_{P,\tilde{\si}}(s)$ has a zero of order eight at $s=1$.  Apart from that, for $\si\in\{1,\tilde{\si}\}$, all poles and zeros of $Z_{P,\si}(s)$ lie in the half-plane $\{\Re(s)\leq \frac{3}{4}\}$.  Further, for $\tau_q$ being the representation of $M$ on $\bwedge^q\bar{\n}\ (q\in\{1,\ldots,4\})$ obtained from the adjoint representation and $\si\in\{\tau_q,\tau_q\ox\tilde{\si}\}$, all poles and zeros of $Z_{P,\si}(s)$ lie in the half-plane $\{\Re(s)\leq 1\}$.
\end{proposition}

Let us assume for a moment that the proposition has been proved.  The representations $\si$ of $M$ considered in the proposition all satisfy the isomorphism of $K_M$-modules $\si\cong\,^w\!\si$, where $w$ is the non-trivial element of the Weyl group $W(G,A)$.  Hence we can apply the functional equation given in Theorem \ref{thm:FunctEqn} to see that the poles and zeros of the functions $Z_{P,\si}(s)$ all lie in the region $0\leq\Re(s)\leq 1$.  We can then apply Theorem \ref{thm:SelbergZeta} to see that the poles and zeros in fact lie in $[0,1]\cup(\rez{2}+i\R)$.  Finally, an application of Theorem \ref{thm:RuelleZeta} completes the proof of Theorem \ref{thm:RuelleMain}.

The proof of the proposition will take up the rest of the section.  We see from (\ref{eqn:v-order}) that a representation $\pi\in\hat{G}$ makes a contribution to the vanishing order of $Z_{P,\si}(s)$ only if $m_{\la}(\pi)\neq 0$ for some $\la\in\a^*$.  From (\ref{eqn:TracePiPhi}) it follows that, if $\tr\pi(\Phi_{s,\si})=0$ for some $\pi\in\hat{G}$ and finite dimensional representation $\si$ of $M$, then $m_{\la}(\pi)=0$ for all $\la\in\a^*$.  Thus, if we can show that $\tr\pi(\Phi_{s,\si})=0$ for some $\pi$ and $\si$, then we will know that $\pi$ makes no contribution to the vanishing order of $Z_{P,\si}$.

Since $\rho_P(H_1)=-\frac{1}{2}$, we shall be able to prove Proposition \ref{pro:SelbergPolesZeros} with the following steps.  Firstly, irrespective of the choice of $\si$, all eigenvalues $\la$ of $\a$ on $H^{\bullet}(\n, \pi_K)$ satisfy
\begin{equation}
\label{eqn:CharIneq2}
\Re(\la)\geq -2\rho_P,
\end{equation}
with equality only if $\pi=triv$.  Secondly, in the case $\si\in\{1,\tilde{\si}\}$, if $\pi$ is the trivial representation on $G$, then $H^4(\n, \pi_K)^{-2\rho_P}\neq 0$.  If $\si=1$ then this gives a double zero at the point $s=1$ and if $\si=\tilde{\si}$ this gives a zero of order eight at $s=1$.  Apart from this, for all $\pi\in\hat{G}$ and all $\la\in\a^*$, we have either $\tr\pi(\Phi)=0$, or that $H^{\bullet}(\n, \pi_K)^{\la}\neq 0$ implies \begin{equation}
\label{eqn:CharIneq1}
\Re(\la)\geq -\frac{3}{2}\rho_P.
\end{equation}

The first step will be proven using the Hochschild-Serre spectral sequence in the following two sections.  Then the second step will be proven using a result of Hecht and Schmid which relates the action of $\a$ on $H^{\bullet}(\n, \pi_K)$ to the infinitesimal character of $\pi$.

\subsection{The Hochschild-Serre spectral sequence}

We shall need the following proposition.

\begin{proposition}
\label{pro:Vanishing}
Let $P'=M'A'N'$ be a parabolic subgroup of $G$ and let $\a'$,$\n'$ be the complexified Lie algebras of $A'$,$N'$ respectively.  Define a partial order $>_{\a'}$ on the dual space $\a'^*$ of $\a'$ by $\mu>\nu$ if and only if $\mu-\nu$ is a non-zero linear combination with positive integral coefficients of positive roots of $(\g,\a)$.

Let $\la\in\a'^*$ and $0\leq p<\dim\n$ be such that $H^p(\n,\pi_K)^{\la}\neq 0$.  Then there exists $\mu\in\a'^*$ with $\la >_{\a'}\mu$ such that $H^{\dim\n}(\n,\pi_K)^{\mu}\neq 0$.
\end{proposition}
\prf
This follows from \cite{HechtSchmid83}, Proposition 2.32 and the isomorphism (\ref{eqn:AM-iso}).
\qed

Let $M_0\subset G$ be the subgroup of diagonal matrices with each diagonal entry equal to $\pm 1$, let $A_0\subset G$ be the subgroup of diagonal matrices with positive entries and let $N_0\subset G$ be the subgroup of upper triangular matrices with ones on the diagonal.  Then $P_0=M_0A_0N_0$ is the minimal parabolic subgroup of $G$.  Let $\a_0$ and $\n_0$ be the complexified Lie algebras of $A_0$ and $N_0$ respectively, and let $\a_0^*$ be the dual of $\a_0$.  Let $\rho_0\in\a_0^*$ be defined as follows:
$$
\rho_0\matrixfour{a}{b}{c}{d}=3a+2b+c.
$$
Then $\rho_0$ is the half-sum of the positive roots of the system $(\g,\a_0)$.  Let $\a_{0,\R}$ be the real Lie algebra of $A_0$, which we may consider as a subalgebra of $\a_0$.  Let $\a_{0,\R}^-$ be the set of all $X\in\a_{0,\R}$ such that $\al(X)<0$ for all positive roots $\al$ of the system $(\g,\a_0)$.  Let $\a_{0,\R}^{*,+}$ be the set of all $\nu\in\a_{0,\R}^*$ such that $\nu(X)<0$ for all $X\in\a_{0,\R}^-$.  Then for all $\nu\in\a_{0,\R}^{*,+}$ we have $\nu=\sum_{\al}\la_{\al}\al$, where the sum is over the positive roots of $(\g,\a_0)$ and every $\la_{\al}>0$.

For $\pi\in\hat{G}$ a \emph{matrix-coefficient} \index{matrix-coefficient} of $\pi$ is any function $G\ra\C$ of the form
$$
f_{u,v}(g)=\langle\pi(g)u,v\rangle
$$
for some $u,v\in\pi$.

\begin{lemma}
\label{lem:MinParAAction}
Let $\pi\in\hat{G}$ and let $\la\in\a_0^*$ be such that $H^{\bullet}(\n_0,\pi_K)^{\la}\neq 0$.  Then $\la\in-2\rho_0+\a_{0,\R}^{*,+}$, except in the case $\pi=triv$, when $H^{6}(\n_0,\pi_K)^{-2\rho_0}\neq 0$ and other than this $H^{\bullet}(\n_0,triv)^{\la}\neq 0$ implies $\la\in-2\rho_0+\a_{0,\R}^{*,+}$.
\end{lemma}
\prf
According to \cite{HechtSchmid83}, Theorem~4.16, there exists a countable set $\CEE(\pi)\subset\a_1^*$ and a collection of polynomial functions $p^{\nu}_{u,v}$ indexed by $u,v\in\pi_K$ and $\nu\in\CEE(\pi)$ such that if $f_{u,v}$, for $u,v\in\pi_K$, is a matrix coefficient of $\pi$ then:
\begin{equation}
\label{eqn:mtxcoeff}
f_{u,v}(\exp\ X)=\sum_{\nu\in\CEE(\pi)} p^{\nu}_{u,v}(X)e^{(\nu +\rho_0)(X)},
\end{equation}
for all $X\in\a_{0,\R}^-$.  We assume that $\CEE(\pi)$ is minimal and hence that each polynomial $p^{\nu}_{u,v}$ is non-zero.

If $\pi=triv$ then every matrix coefficient of $\pi$ is a constant function so in (\ref{eqn:mtxcoeff}) we have $\CEE(\pi)=\{-\rho_0\}$ and $p^{-\rho_0}_{u,v}=C_{u,v}$, where $C_{u,v}$ is a constant for all $u,v\in\pi$.

If $\pi\in\hat{G}\smallsetminus\{triv\}$ then $\pi$ is infinite dimensional so by \cite{HoweMoore79}, Theorem 5.1 all the matrix coefficients of $\pi$ vanish at infinity.  Hence, for all $\nu\in\CEE(\pi)$ and for all $X\in\a_{0,\R}^-$, we have $(\nu+\rho_0)(X)<0$.  In other words, for all $\nu\in\CEE(\pi)$ we have $\nu\in-\rho_0+\a_{0,\R}^{*,+}$.

We have assumed that $H^{\bullet}(\n_0,\pi_K)^{\la}\neq 0$.  Let $>_{\a_0}$ be the partial order defined on $\a_0^*$ by $\mu>_{\a_0}\nu$ if and only if $\mu-\nu$ is a linear combination with positive integral coefficients of positive roots of $(\g,\a_0)$.  By Proposition \ref{pro:Vanishing} there exists $\mu\in\a_0^*$ such that $H^6(\n_0,\pi_K)^{\mu}\neq 0$ and $\la\geq_{\a_0}\mu$.

Let $\La=\{\nu+\rho_0\in\a_0^*:H^6(\n_0,\pi_K)^{\nu}\neq 0\}$ and let $\La^{\min}$ be the set of elements of $\La$ which are minimal with respect to $>_{\a_0}$.  Let $\CEE^{\min}(\pi)$ be the set of $\nu\in\CEE(\pi)$ which are minimal with respect to $>_{\a_0}$.  Theorem~4.25 of \cite{HechtSchmid83} says that $\La^{\min}=\CEE(\pi)^{\min}$.

If $\pi=triv$ then $\La^{\min}=\{-\rho_0\}$ and the claims of the proposition follow from the definition of $\La$ and \cite{HechtSchmid83}, Proposition~2.32 as quoted above.

If $\pi\in\hat{G}\smallsetminus\{triv\}$ then $\La^{\min}=\CEE(\pi)^{\min}$ implies that there exists $\nu\in\CEE(\pi)$ such that $\la\geq_{\a_0}\nu-\rho_0$.  We saw above that $\nu\in-\rho_0+\a_{0,\R}^{*,+}$ so the claim follows.
\qed

\begin{proposition}
\label{pro:AAction}
Let $\pi\in\hat{G}$ and let $\la\in\a^*$ be such that $H^{\bullet}(\n,\pi_K)^{\la}\neq 0$.  Then $\Re(\la)\geq-2\rho_P$, with equality only if $\pi=triv$.
\end{proposition}
\prf
Let $\n_M=\n_0\cap\m$ and $\a_M=\a_0\cap\m$.  Then $\n_0=\n_M\oplus\n$ and $\a_0=\a\oplus\a_M$.  In \cite{HochschildSerre53b, HochschildSerre53a} a filtered complex is constructed so that the spectral sequence derived from it has
$$
E_2^{p,q} \cong H^p(\n_M,H^q(\n,\pi_K))
$$
and
$$
E_{\infty}^{p,q} \cong \Gr^p H^{p+q}(\n_0,\pi_K),
$$
where $H^{p+q}(\n_0,\pi_K)$ is appropriately filtered.  This is the so-called \emph{Hochschild-Serre spectral sequence}\index{spectral sequence!Hochschild-Serre}.

By Proposition \ref{pro:Vanishing} it suffices to prove the proposition in the case
$$
H^4(\n,\pi_K)^{\la}\neq 0.
$$
In this case it then follows that
$$
H^2(\n_M,H^4(\n,\pi_K)^{\la})\neq 0
$$
and so there exists $\la_M\in\a_M^*$ such that
$$
H^2(\n_M,H^4(\n,\pi_K)^{\la})^{\la_M}\neq 0.
$$
Since $A$ acts trivially on $\n_M$, this equals
$$
H^2(\n_M,H^4(\n,\pi_K))^{\la+\la_M}=(E_2^{2,4})^{\la+\la_M},
$$
where we consider $\la+\la_M$ as an element of $\a_0^*=\a^*\oplus\a_M^*$.  For all $r$ we have that $E_r^{p,q}\neq 0$ only if $0\leq p\leq 2$ and $0\leq q\leq 4$, so it follows that $E_2^{2,4}=E_{\infty}^{2,4}$.  Since the action of $A_0$ commutes with the differentials of the spectral sequence it follows that
$$
(E_{\infty}^{2,4})^{\la+\la_M}\neq 0
$$
and hence that
$$
H^6(\n_0,\pi_K)^{\la+\la_M}\neq 0.
$$
The proposition then follows from Lemma \ref{lem:MinParAAction} by projection of $\la+\la_M$ onto the $\a^*$ component.
\qed

\subsection{Contribution of the trivial representation}

For a $K_M$-module $\eta$ let $2\eta$ denote the module $\eta\oplus\eta$.  We shall need the following:
\begin{lemma}
\label{lem:pMKMTypes}
We have the following isomorphisms of $K_M$-modules:
\begin{eqnarray*}
  \bwedge^0 \p_M & \cong & triv \\
  \bwedge^1 \p_M & \cong & \d_{2,0}\oplus\d_{0,2} \\
  \bwedge^2 \p_M & \cong & 2\d\oplus\d_{2,2}\oplus\d_{2,-2} \\
  \bwedge^3 \p_M & \cong & \d_{2,0}\oplus\d_{0,2} \\
  \bwedge^4 \p_M & \cong & triv.
\end{eqnarray*}
and
\begin{eqnarray*}
\bwedge^0 \m & \cong & triv \\
\bwedge^1 \m & \cong & 2\d\oplus\d_{2,0}\oplus\d_{0,2} \\
\bwedge^2 \m & \cong & triv\oplus2\d\oplus2(\d_{2,0}\oplus\d_{0,2})\oplus\d_{2,2}\oplus\d_{2,-2} \\
\bwedge^3 \m & \cong & 2triv\oplus2triv\oplus2(\d_{2,0}\oplus\d_{0,2}\oplus\d_{2,2}\oplus\d_{2,-2}) \\
\bwedge^4 \m & \cong & triv\oplus2\d\oplus2(\d_{2,0}\oplus\d_{0,2})\oplus\d_{2,2}\oplus\d_{2,-2}
\end{eqnarray*}
\end{lemma}

\prf
The isomorphisms for $p_M$ were given in Lemma \ref{lem:pMKMTypes1}.  For $\m$ note that $K_M$ acts on $\m$ by the adjoint representation and we can compute
$$
\m\cong 2\d\oplus\d_{2,0}\oplus\d_{0,2}.
$$
The other isomorphisms follow straightforwardly from this.
\qed

\begin{proposition}
\label{pro:TrivAAction}
Let $\pi$ be the trivial representation on $G$.

For $\si$ the trivial representation on $M$, the representation $\pi$ contributes a double zero of $Z_{P,\si}(s)$ at the point $s=1$.  For $\si=\tilde{\si}$ the representation $\pi$ contributes a zero of $Z_{P,\si}(s)$ of order eight at the point $s=1$.  In both cases, all other poles and zeros contributed by $\pi$ are in $\left\{\Re(s)\leq\frac{3}{4}\right\}$.
\end{proposition}
\prf
The space $H_0(\n,\pi_K)=\pi_K /\n\pi_K$ is one dimensional with trivial $\a$-action.  The action of $\a$ on $\n$ is given by $\rez{2}\rho_P$, so the isomorphism (\ref{eqn:AM-iso}) tells us that $\a$ acts on the one dimensional space $H^4(\n,\pi_K)$ according to $-2\rho_P$.  Proposition \ref{pro:Vanishing} tells us that for $q=0,1,2,3$ and $\la\in\a^*$, $H^q(\n,\pi_K)^{\la}\neq 0$ implies $\la\geq -\frac{3}{2}\rho_P$.  Since $-2\rho_P(H_1)=1$ this gives a pole or zero at $s=1$ and evaluation of the relevant characters at $H_1$ also gives the other poles and zeros contributed by $\pi$ in $\{\Re(s)\leq\frac{3}{4}\}$.

It remains to compute the vanishing order of $Z_{P,\si}(s)$ for $\si\in\{1,\tilde{\si}\}$ at the point $s=1$.  Since $\dim H^4(\n,\pi_K)^{-2\rho_P}=1$, we get from (\ref{eqn:v-order}) the following expression for the vanishing order:
$$
N_{\Ga}(\pi) \sum_{p\geq 0} (-1)^{p} \dim\left(\bwedge^p \p_M\ox V_{\si}\right)^{K_M}.
$$
For $\si=1$ this is equal to
$$
N_{\Ga}(\pi) \sum_{p\geq 0} (-1)^{p} \dim\left(\bwedge^p \p_M\right)^{K_M}
$$
and for $\si=\tilde{\si}$ to
$$
N_{\Ga}(\pi) \sum_{p,q\geq 0} (-1)^{p+q} \dim\left(\bwedge^p \p_M\ox a_q\bwedge^q\m\right)^{K_M},
$$
where $a_0=15$, $a_1=10$, $a_2=6$, $a_3=3$, $a_4=1$ and $a_q=0$ for $q\geq 5$.  The only functions on $L^2(\Ga\bs G)$ invariant under the action of $G$ are the constant functions, hence $N_{\Ga}(\pi)=1$.  We can then use Lemma \ref{lem:pMKMTypes} to see that the above expressions take the claimed values.
\qed

\subsection{Contribution of the other representations}

The following proposition gives a relationship between the action of $\a$ on $H^{\bullet}(\n,\pi_K)$ and the infinitesimal character of $\pi$.

\begin{proposition}
\label{pro:aCohomologyAction}
Let $\pi\in\hat{G}$.

Suppose $H^{\bullet}(\n, \pi_K)^{\la}\neq 0$ for some $\la\in\a^*$.  Then $\la=w\La_{\pi}|_{\a}-\rho_P$, where $w\in W(\g,\h)$ and $\La_{\pi}$ is a representative of the infinitesimal character of $\pi$.  Moreover, for $p\in\Z$ we have
$$
H^p(\n,\pi_K)=\bigoplus_{w\in W(\g,\h)} H^p(\n,\pi_K)^{w\La_{\pi}-\rho_P}.
$$
\end{proposition}
\prf
This follows from \cite{HechtSchmid83}, Corollary~3.32 and the isomorphism (\ref{eqn:AM-iso}).
\qed

In light of this proposition, in order to complete the proof of Proposition \ref{pro:SelbergPolesZeros}, and hence of Theorem \ref{thm:RuelleMain}, it will be sufficient to show that for all $\pi\in\hat{G}\smallsetminus\{triv\}$ and for $\si\in\{1,\tilde{\si}\}$ either $\tr\pi(\Phi_{\si})=0$ or the infinitesimal character $\La_{\pi}$ of $\pi$ satisfies
$$
\Re(w\La_{\pi})|_{\a}\geq-\rez{2}\rho_P\ \textrm{ or }\ -\rho_P\geq\Re(w\La_{\pi})|_{\a}
$$
for all $w\in W(\g,\h)$.

We first consider the elements of $\hat{G}$ induced from parabolic subgroups other that $P=MAN$.

We say that two parabolics $P'=M'A'N'$ and $P''=M''A''N''$, where 
$$
A', A''\subset A_0 =\left\{\diag(a,b,c,(abc)^{-1})|a,b,c>0\right\}
$$ 
are \emph{associate} if there exists a member $w$ of the Weyl group $W(\g,\h)$ such that $w^{-1}MAw=M'A'$.  Representations of $G$ induced from associate parabolics are equivalent.

Up to association G has four parabolic subgroups: $P$ defined above; the minimal parabolic $P_0=M_0 A_0 N_0$ where $M_0 \cong \{\pm 1\}$ and $N_0$ is the group of real, upper triangular matrices with ones on the diagonal; the parabolic $P'$ with Langlands decomposition $P'=M'A'N'$ where $M'\cong \SL_2^{\pm}(\R)\times\{\pm 1\}$, $A'=\left\{\diag(a,a,b,a^{-2}b^{-1})|a,b>0\right\}$ and $N'$ is the group of real, upper triangular matrices with ones on the diagonal and otherwise with zeros in the second column; and the parabolic $P''$ with Langlands decomposition $P''=M''A''N''$ where $M''\cong \SL_3^{\pm}(\R)$, $A''=\left\{\diag(a,a,a,a^{-3})|a>0\right\}$ and $N''$ is the group of real, upper triangular matrices with ones on the diagonal whose only non-zero entries above the diagonal are in the rightmost column.  

\begin{proposition}
\label{pro:parabolics}
Let $\pi$ be a principal series or complementary series representation induced from the parabolic $\bar{P}=\bar{M}\bar{A}\bar{N}$, where $\bar{P} = P_0$, $P'$ or $P''$ and let $\si$ be a finite dimensional representation on $M$.  Then $\tr\pi(\Phi_{\si})=0$ so $\pi$ contributes no zeros or poles to $Z_{P,\si}(s)$.
\end{proposition}
\prf
Let $\Th_{\pi}$ be the global character of the irreducible unitary representation $\pi$, so that
$$
\tr\pi(\Phi_{\si}) = \int_G \Phi_{\si}(x)\Th_{\pi}(x)\ dx.
$$
The Weyl integration formula can be applied (see \cite{geom}, p908) to give us
$$
\tr\pi(\Phi_{\si}) = \sum_L \rez{|W(L)|} \int_{A^+ L} \int_{M/L} f_{\si}\left( mlm^{-1}\right)dm\ \Th_{\pi}(al)d(al)\ dl,
$$
where the sum is over conjugacy classes of Cartan subgroups $L$ of $M$, and we denote by $W(L)=W(L,M)$ the Weyl group of $L$ in $M$, by $f_{\si}$ an Euler-Poincar$\eac$ function for $\si$ and $d(al)$ is an explicitly given function on $A^+ L$.  Proposition~1.4 of \cite{geom} tells us that for $L\neq B$
$$
\int_{M/L} f_{\si}\left( mlm^{-1}\right)dm=0,
$$
hence
$$
\tr\pi(\Phi_{\si}) = \rez{|W(B)|} \int_{A^+ B} \int_{M/B} f_{\si}\left( mbm^{-1}\right)dm\ \Th_{\pi}(ab)d(ab)\ db.
$$

The character $\Th_{\pi}$ of $\pi$ is non-zero only on Cartan subgroups of $G$ that are $G$-conjugate to Cartan subgroups of $\bar{M}\bar{A}$ (see \cite{Knapp}, Proposition~10.19).  The subgroup $BA$ is not G-conjugate to any Cartan subgroup of $\bar{M}\bar{A}$, so it follows that $\tr\pi(\Phi_{\si})=0$.
\qed

Now let $\pi=\Ind_P^G(\xi\ox\nu)$ for some $\xi\in\hat{M}$ and $\nu\in\a^*$.  For $\tau \in \widehat{K_M}$ let $P_{\tau}:V_{\xi}\rightarrow V_{\xi}(\tau)$ be the projection onto the $\tau$-isotype.  For any function $f$ on $G$ which is sufficiently smooth and of sufficient decay the operator $\pi(f)$ is of trace class.  Its trace is
$$
\sum_{\tau \in \widehat{K_M}}\int_K \int_{MAN} a^{\nu + \rho_P} f(k^{-1}mank)\ \tr P_{\tau}\xi(m)P_{\tau}\ dman\ dk.
$$
Plugging in the test function $f=\Phi_{\si}$, where $\si\in\{triv,\tilde{\si}\}$, this gives us, as in \cite{class}:
$$
\tr\pi(\Phi_{\si})=\int_{A^+} C(a)\ \tr\xi(f_{\si})\ da,
$$
where $C(a)$ depends only on $a$ and $f_{\si}$ is an Euler-Poincar\'{e} function on M attached to the representation $\si$.  We can see that $\tr\pi(\Phi_{\si})$ is non-zero only if $\tr\xi(f_{\si})$ is also.

Theorem \ref{thm:RuelleMain} now follows from Proposition \ref{pro:InfChars}.

\section{A Prime Geodesic Theorem for $\SL_4(\R)$}
    \label{ch:PGT}

We continue using the notation defined in the previous sections, in particular we take $G=\SL_4(\R)$ and take $\Ga\subset G$ to be a discrete, cocompact subgroup.  For $\ga\in\Ga$ let $N(\ga)=e^{l_{\ga}}$ \index{$N(\ga)$} and define for $x>0$
$$
\pi(x)=\sum_{{[\ga]\in\CE^{p}_P(\Ga)}\atop{N(\ga)\leq x}}\chi_1(\Ga_{\ga}),\index{$\pi(x)$}\ \textrm{ and }\ \ 
\tilde{\pi}(x)=\sum_{{[\ga]\in\CE^{p,\reg}_P(\Ga)}\atop{N(\ga)\leq x}}\chi_1(\Ga_{\ga})\tr\tilde{\si}(b_{\ga}),\index{$\tilde{\pi}(x)$}
$$
where $\Ga_{\ga}$ is the centraliser of $\ga$ in $\Ga$, we denote by $\chi_1(\Ga_{\ga})$ the first higher Euler characteristic and $\tilde{\si}$ is the virtual representation defined in Section \ref{sec:RuelleZetaFn}.  We recall from Proposition \ref{thm:ECharPos} that the first higher Euler characteristics are all positive and from Lemma \ref{lem:SigmaTilde} that the traces $\tr\tilde{\si}(b_{\ga})$ are also all positive, so both $\pi(x)$ and $\tilde{\pi}(x)$ are monotonically increasing functions.

\begin{theorem}\textnormal{(Prime Geodesic Theorem)}
\label{thm:PGT}
\index{Prime Geodesic Theorem}
For $x\rightarrow\infty$ we have
$$
\pi(x)\sim\frac{2x}{\log x}\ \textrm{ and }\ \tilde{\pi}(x)\sim\frac{8x}{\log x}.
$$
More sharply,
$$
\pi(x)=2\,\li(x)+O\left(\frac{x^{3/4}}{\log x}\right)\ \textrm{ and }\ \tilde{\pi}(x)=8\,\li(x)+O\left(\frac{x^{3/4}}{\log x}\right)
$$
as $x\rightarrow\infty$, where $\li(x)=\int_2^x \rez{\log t} dt$ \index{$\li(x)$} is the integral logarithm.
\end{theorem}

We also prove the following Prime Geodesic Theorem, which will be needed for our application to class numbers.

Let $B^0$ \index{$B^0$} be a closed subset of the compact group $B$ with the following properties: it is of measure zero; it is invariant under the map $b\mapsto b^{-1}$ and contains all fixed points of this map; and its complement $B^1=B\smallsetminus B^0$ \index{$B^1$} in $B$ is homeomorphic to an open subset of Euclidean space each of whose connected components is contractible.  The assumption that $B^0$ contains all fixed points of the map $b\mapsto b^{-1}$ is equivalent to the assumption that $B^0$ contains all non-regular elements of $B$.  Let $\CE_P^{p,1}(\Ga)$ be the subset of all $[\ga]\in\CE_P^p(\Ga)$ such that $\ga$ is conjugate in $G$ to an element of $A^-B^1$.  We define for $x>0$
$$
\pi^1(x)=\sum_{{[\ga]\in\CE^{p,1}_P(\Ga)}\atop{N(\ga)\leq x}}\chi_1(\Ga_{\ga}).\index{$\pi^1(x)$}
$$
Then as in the case of $\pi(x)$ above, $\pi^1(x)$ is a monotonically increasing function.

\begin{theorem}
\label{thm:PGT2}
For $x\rightarrow\infty$ we have
$$
\pi^1(x)\sim\frac{2x}{\log x}.
$$
\end{theorem}

The proof of these two theorems will occupy the rest of the section.

\subsection{Analytic properties of $R_{\Ga,\si}(s)$}

This subsection and the following proceed according to the methods of \cite{Randol77} and \cite{Randol78}, in which the Selberg zeta function for quotients of the hyperbolic plane are considered.  In this section the analogs of a series of lemmas are proved in our context.  The next section translates the main theorem of \cite{Randol77} into our context.

Recall from Theorem \ref{thm:FunctEqn} that for a finite dimensional virtual representation $\si$ of $M$ we have the functional equation
\begin{equation}
\label{eqn:FunctionalEqn}
Z_{P,\si}(1-s)=e^{-G(s)}Z_{P,\si}(s),
\end{equation}
where $G(s)$ is a polynomial.  Let $D$ be the degree of the polynomial $G(x)$.

\begin{lemma}
\label{lem:RuelleLogDerEst1}
Let $H$ be a half-plane of the form $\left\{\Re(s)<-(1+\ep)\right\}$ for some $\ep>0$ and let $\si$ be a finite dimensional virtual representation of $M$.  Then there exists a constant $C>0$ such that for $s\in H$ we have
$$
|R_{\Ga,\si}'(s)/R_{\Ga,\si}(s)|\leq C|s|^{D-1}.
$$
\end{lemma}

\prf
From Theorem \ref{thm:RuelleZeta} we get the identity
$$
R_{\Ga,\si}(s) = \prod_{q=0}^4 Z_{P,(\bwedge^q\bar{\n}\ox V_{\si})}\left( s+\frac{q}{4}\right)^{(-1)^q},
$$
which implies
\begin{equation}
\label{eqn:RuelleSelbergLogDer}
\frac{R'_{\Ga,\si}(s)}{R_{\Ga,\si}(s)} = \sum_{q=0}^4 (-1)^q \frac{Z'_{P,(\bwedge^q\bar{\n}\ox V_{\si})}\left( s+\frac{q}{4}\right)}{Z_{P,(\bwedge^q\bar{\n}\ox V_{\si})}\left( s+\frac{q}{4}\right)}.
\end{equation}
Considering this identity, it will suffice to prove that when $K$ is a half-plane of the form $\left\{\Re(s)<-\ep\right\}$ for some $\ep>0$, there exists a constant $C>0$ such that for $s\in K$ we have
$$
|Z'_{P,(\bwedge^q\n\ox V_{\si})}(s)/Z_{P,(\bwedge^q\bar{\n}\ox V_{\si})}(s)|\leq C|s|^{D-1}
$$
for all $q=0,\ldots,4$.  Since the proof does not depend on the value of $q$ or on $\si$ we shall abbreviate our notation for the zeta function to $Z_P(s)$, which notation we shall use for the rest of the section.

It follows that
$$
-\frac{Z'_P(1-s)}{Z_P(1-s)}=\frac{Z'_P(s)}{Z_P(s)} - G'(s).
$$
From the definition (\ref{eqn:GenSelb}) of $Z_P(s)$ and Proposition \ref{pro:SelbergPolesZeros}, we can see that $Z_P(s)$ is both bounded above and bounded away from zero on the half plane $K'=\left\{\Re(s)>1+\ep\right\}$.  It follows that $Z'_P(s)/Z_P(s)$ is bounded on $K'$.  This implies the lemma.
\qed

For $t>0$, let $N(t)$\index{$N(t)$}, denote the number of poles and zeros of $Z_P(s)$ at points $s=\rez{2}+x$, where $0<x<t$.

\begin{lemma}
\label{lem:NEstimate}
$N(t)=O(t^D)$
\end{lemma}

\prf
Define $\xi(s)=\left(Z_P(s)\right)^2 e^{-G(s)}$, where $G(s)$ is the polynomial in the functional equation (\ref{eqn:FunctionalEqn}) for $Z_P(s)$.  We note that, in light of its role in the functional equation, the polynomial $G(s)$ must satisfy $G(s)=-G(1-s)$.  It then follows that
$$
\xi(1-s)=\xi(s).
$$
Note that $\xi(s)$ is real on the real axis and so $\xi(\bar{s})=\overline{\xi(s)}$.

Fix a real number $1<a<5/4$ and let $t>0$.  Let $R$ be the rectangle defined by the inequalities $1-a\leq\Re s\leq a$ and $-t\leq\Im a\leq t$.  We assume that $t$ has been chosen so that no zero or pole occurs on the boundary of $R$.  Then
$$
N(t)=\rez{4}\cdot\rez{2\pi i}\int_{\partial R}\frac{\xi'(s)}{\xi(s)}ds - N_0=\rez{4}\cdot\rez{2\pi}\Im\left(\int_{\partial R}\frac{\xi'(s)}{\xi(s)}ds\right) - N_0,
$$
where $N_0$ is the number of poles and zeros of $Z_P(s)$ on the real line.  It then follows from the functional equation for $\xi(s)$ and from the fact that $\xi(\bar{s})=\overline{\xi(s)}$, that we have
$$
N(t)=\rez{2\pi}\Im\left(\int_{C}\frac{\xi'(s)}{\xi(s)}ds\right) - N_0,
$$
where $C$ is the portion of $\partial R$ consisting of the vertical segment from $a$ to $a+it$ plus the horizontal segment from $a+it$ to $\rez{2}+it$.

Now the definition of $\xi(s)$ gives us
$$
\frac{\xi'(s)}{\xi(s)}=2\frac{Z'_P(s)}{Z_P(s)} - G'(s),
$$
so that
\begin{eqnarray*}
\Im\left(\int_{C}\frac{\xi'(s)}{\xi(s)}ds\right) & = & 2.\Im\left(\int_{C}\frac{Z'_P(s)}{Z_P(s)}ds\right) - \Im\left(\int_{C}G'(s)ds\right) \\
  \\
  & = & 2.\Im\left(\int_{C}\frac{Z'_P(s)}{Z_P(s)}ds\right) - \Im\left(G\left(\rez{2}+it\right)+G(a)\right) \\
  \\
  & = & 2.\Im\left(\int_{C}\frac{Z'_P(s)}{Z_P(s)}ds\right) + O(t^D).
\end{eqnarray*}

It thus remains to show that
$$
S(t)=\Im\left(\int_{C}\frac{Z'_P(s)}{Z_P(s)}ds\right)=O(t^D).
$$
Note that $S(t)$ is the variation of the argument of $Z_P(s)$ along $C$.  We may extend the definition of $S(t)$ to those values of $t$ for which $\rez{2}+it$ is a pole or zero of $Z_P(s)$ by defining it to be $\lim_{\ep\ra 0}\rez{2}(S(t+\ep)+S(t-\ep))$.

From the definition of $Z_P(s)$ we can see that $S(t)=h(t)+O(1)$, where $h(t)$ is the variation of the argument of $Z_P(s)$ along the segment from $a+it$ to $\rez{2}+it$.  The value of $h(t)$ is bounded by a multiple of the number of zeros of $\Re (Z_P(s))$ on this segment, since the point $Z_P(s)$ cannot move between the right and left half-planes without crossing the imaginary axis.  On the segment, the real part of $Z_P(s)$ coincides with
$$
f_P(w)=\rez{2}(Z_P(w+it)+Z_P(w-it)),
$$
where $w$ runs from $\rez{2}$ to $a$ on the real axis.  Since we have assumed that $\rez{2}+it$ is not a zero or pole of $Z_P(s)$, the function $f_P(w)$ is holomorphic in a neighbourhood of the closed disc $S$, centred at $a$, of radius $a-\rez{2}$.  As this disc contains the interval from $\rez{2}$ to $a$, we may apply Jensen's Formula (\cite{Conway78}, XI.1.2) to conclude that
\begin{eqnarray*}
h(t) & = & O\left(\int_{\partial S} \log|f_P(w)|dw\right) \\
  & = & O\left(\int_{\partial S}\log |Z_P(w+it)+Z_P(w-it)|dw\right) \\
  & = & O\left(\int_{\partial S}\log |Z_P(w+it)|+\log|Z_P(w-it)|dw\right) \\
  & = & O\left(\sum_{i=1,2}\int_{\partial S}\log |Z_i(w+it)|dw+\int_{\partial S}\log |Z_i(w-it)|dw\right).
\end{eqnarray*}

There remains one final step in the proof of the lemma:
\begin{equation}
\label{eqn:Estimate}
|Z_i(x+it)|=e^{O(|t|^D)}
\end{equation}
uniformly in the strip $1-a\leq x\leq a$ for $i=1,2$.  We know that the function $Z_P(s)$ is bounded and bounded away from zero on the line $\Re s=a$.  Hence, by Lemma \ref{lem:FiniteOrder} and the functional equation (\ref{eqn:FunctionalEqn}), we can apply the Phragm\'en Lindel\"of Theorem \index{Phragm\'en Lindel\"of Theorem} (\cite{Conway78}, Theorem VI.4.1) to give (\ref{eqn:Estimate}).
\qed

\begin{lemma}
\label{lem:RuelleLogDerEst2}
Given $a<b\in\R$, there exists a sequence $(y_n)$ tending to infinity such that
$$
\left|\frac{R_{\Ga}'(x+iy_n)}{R_{\Ga}(x+iy_n)}\right|=O(y_n^{2D})
$$
for $a<x<b$.
\end{lemma}

\prf
As in Lemma \ref{lem:RuelleLogDerEst1} it will suffice to prove the result for $Z_P(s)$.

Using the notation of Lemma \ref{lem:FiniteOrder} we have that
$$
\frac{Z'_P(s)}{Z_P(s)}=\rez{s}(n_1-n_2)+g'_1(s)-g'_2(s)+\sum_{i=1,2}(-1)^{i-1}\sum_{\rho\in R_i\smallsetminus\{0\}}s^k\rho^{-k}(s-\rho)^{-1}.
$$
Let $t_0>2$ be fixed and consider the segment of the line $\Re s=\rez{2}$ given by $\rez{2}+it$ for $t_0-1<t\leq t_0+1$.  Let $N(t)$ be as above, then by Lemma \ref{lem:NEstimate} we know that $N(t)=O(t^D)$.  It follows immediately that the number of roots on the segment is $O(t_0^D)$.

By the Dirichlet principle, there exists a $\rez{2}+iy$ in the segment whose distance from any pole or zero is greater that $C/T^D$, for some fixed $C>0$.  We conclude that the portion of the sum 
$$
\sum_{i=1,2}(-1)^{i-1}\sum_{\rho}s^k\rho^{-k}(s-\rho)^{-1}
$$ 
corresponding to poles and zeros in the segment for $s_x=x+iy$ is $O(y^{2D})$, since $|s_x^k\rho^{-k}|=O(1)$ for these $\rho$, when $a<x<b$.

To deal with the segments $\rez{2}+it$ for $0<t\leq t_0-1$ and $t_0+1<t<\infty$, we proceed as follows.  The portions of the sum $\sum_{i=1,2}(-1)^{i-1}\sum_{\rho}s^k\rho^{-k}(s-\rho)^{-1}$ corresponding to the first and second segments respectively, can be written
$$
s_x^k\int_0^{t_p-1}\left(\rez{2}+it\right)^{-k}\left(s_x-\rez{2}-it\right)^{-1} dN(t)
$$
and
$$
s_x^k\int_{t_p+1}^{\infty}\left(\rez{2}+it\right)^{-k}\left(s_x-\rez{2}-it\right)^{-1} dN(t).
$$
Recalling that $N(t)=O(t^D)$, we easily conclude that both of these expressions are $O(T^{2D})$.
\qed

\subsection{Estimating $\psi(x)$ and $\tilde{\psi}(x)$}

To simplify notation, in what follows we write $\ga$ for an element of $\CE_P(\Ga)$ and $\ga_0$ for a primitive element and recall that $\ga\in\CE_P(\Ga)$ implies that $\ga$ is conjugate in $\Ga$ to an element $a_{\ga}b_{\ga}\in A^-B$.  Unless otherwise specified, sums involving $\ga$ or $\ga_0$ will be taken over conjugacy classes in $\CEE_P(\Ga)$ and $\CEE_P^p(\Ga)$ respectively.  If $\ga$ and $\ga_0$ occur in the same formula it is understood that $\ga_0$ wil be the primitive element underlying $\ga$.  We denote by $\CE_P^{\reg}(\Ga)$ \index{$\CE_P^{\reg}(\Ga)$} the subset of regular elements in $\CE_P(\Ga)$.  For $x>0$ let
$$
\psi(x)=\sum_{{[\ga]\in\CE_P(\Ga)}\atop{N(\ga)\leq x}} \chi_1(\Ga_{\ga})l_{\ga_0}
$$\index{$\psi(x)$}
and
$$
\tilde{\psi}(x)=\sum_{{[\ga]\in\CE_P^{\reg}(\Ga)}\atop{N(\ga)\leq x}}\chi_1(\Ga_{\ga_0})\tr\tilde{\si}(b_{\ga})l_{\ga_0}.\index{$\tilde{\psi}(x)$}
$$

Let $\si$ be a finite dimensional virtual representation of $M$.  We have for $\Re(s)>1$:
\begin{equation}
\label{eqn:Ruellelogder}
\frac{R_{\Ga,\si}'(s)}{R_{\Ga,\si}(s)}=\sum_{\ga}\chi_1(\Ga_{\ga})\tr\si(b_{\ga})l_{\ga_0}e^{-sl_{\ga}}.
\end{equation}

The following propositions are analogs of Theorem 2 of \cite{Randol77}, from which, in the next subsection, we prove the Prime Geodesic Theorem using standard techniques of analytic number theory.

\begin{proposition}
\label{pro:PsiHatEst}
$\psi(x)=2x+O\left(x^{3/4}\right)$
\end{proposition}
\prf
Let $D$ be the degree of the polynomial $G(s)$, as in the previous section, and suppose $k\geq 2D$ is an integer and $x,c>1$.  Then, by (\ref{eqn:Ruellelogder}) and \cite{HardyRiesz64}, Theorem~40,
\begin{eqnarray}
\label{eqn:Ruellelogderint}
\lefteqn{\rez{2\pi i}\int_{c-i\infty}^{c+i\infty} \frac{R_{\Ga,1}'(s)}{R_{\Ga,1}(s)} s^{-1}(s+1)^{-1}\cdots(s+k)^{-1}x^s ds} \nonumber \\
  \nonumber \\
  & = & \rez{2\pi i}\int_{c-i\infty}^{c+i\infty} \left(\sum_{\ga} \chi_1(\Ga_{\ga})l_{\ga_0}e^{-sl_{\ga}}\right) s^{-1}(s+1)^{-1}\cdots(s+k)^{-1}x^s ds \nonumber \\
  \nonumber \\
  & = & \rez{k!} \sum_{N(\ga)\leq x} \chi_1(\Ga_{\ga})l_{\ga_0}\left(1-\frac{N(\ga)}{x}\right)^k.
\end{eqnarray}

Theorem~\ref{thm:RuelleMain} tells us that all poles of $R_{\Ga}'(s)/R_{\Ga}(s)$ lie in the strip $-1\leq\Re(s)\leq 1$.  By virtue of Lemmas~\ref{lem:RuelleLogDerEst1} and \ref{lem:RuelleLogDerEst2} it is permissible to shift the line of integration in (\ref{eqn:Ruellelogderint}) into the half plane $\Re(s)<-1$, taking into account the residues of the poles of $R_{\Ga}'(s)/R_{\Ga}(s)$.  Hence, for $c'<-1$
\begin{eqnarray}
\label{eqn:RuellePolesSum}
\lefteqn{\rez{k!} \sum_{N(\ga)\leq x} \chi_1(\Ga_{\ga})l_{\ga_0}\left(1-\frac{N(\ga)}{x}\right)^k} \\
 & = & \sum_{{\al\in S_k}\atop{\Re(\al)>c'}} c_k(\al)x^{\al} + \rez{2\pi i}\int_{c'-i\infty}^{c'+i\infty} \frac{R_{\Ga,1}'(s)}{R_{\Ga,1}(s)} s^{-1}(s+1)^{-1}\cdots(s+k)^{-1}x^s ds, \nonumber
\end{eqnarray}
where $S_k$ denotes the set of poles of $(R_{\Ga,1}'(s)/R_{\Ga,1}(s))s^{-1}(s+1)^{-1}\cdots(s+k)^{-1}$ and $c_k(\al)$ denotes the residue at $\al$.

For $x>1$, Lemma~\ref{lem:RuelleLogDerEst1} implies the integral in (\ref{eqn:RuellePolesSum}) tends to zero as $c'\rightarrow -\infty$.  If we define
$$
\psi_0(x) =\psi(x),\ \ \psi_j(x)=\int_0^x \psi_{j-1}(t)dt,\ j\in\N,
$$
it is well known that
$$
\psi_j(x)=\rez{j!}\sum_{N(\ga)\leq x} \chi_1(\Ga_{\ga})l_{\ga_0} (x-N(\ga))^j
$$
and we deduce from (\ref{eqn:RuellePolesSum}) that
\begin{equation}
\label{eqn:Psik}
\psi_k(x)=\sum_{\al\in S_k} c_k(\al) x^{k+\al}.
\end{equation}

Let $d>0$.  For a function $f:\R\rightarrow\R$ define the operator $\D$ by setting
$$
\D f(x)=\sum_{i=0}^{2D}(-1)^i\left({2D}\atop{i}\right)f(x+(2D-i)d).
$$
It follows from (\ref{eqn:Psik}) and Theorem~\ref{thm:RuelleMain}, setting $k=2D$, that
\begin{equation}
 \psi_{2D}(x)= \frac{2}{(2D+1)!}x^{2D+1} + \sum_{\al\in S_{2D}^{\R}} c_{2D}(\al)x^{2D+\al} + \sum_{p=-2}^2 \sum_{\al\in S_{2D}^{p/4}} c_{2D}(\al)x^{2D+\al},
\end{equation}
where $S_{2D}^{\R}=S_{2D}\cap (\R\smallsetminus\{1\})$, the real elements of $S_{2D}$ not including $\al=1$, and $S_{2D}^q=S_{2D}\cap (q+i(\R\smallsetminus \{0\}))$, the non-real elements of $S_{2D}$ on the line $\Re(s)=q$.  The coefficient $((2D+1)!)^{-1}$ on the leading term comes from the fact that $(R_{\Ga,1}'(s)/R_{\Ga,1}(s))$ has a double pole at $s=1$ and from the factors $s^{-1}...(s+2D)^{-1}$.

In general if $f$ is at least $2D$ times differentiable,
\begin{equation}
 \label{eqn:DeltaId}
 \D f(x)=\int_x^{x+d} \int_{t_{2D}}^{t_{2D}+d}\cdots\int_{t_2}^{t_2+d} f^{(2D)}(t_1)\ dt_1...dt_{2D}.
\end{equation}
By applying the Mean Value Theorem we get
\begin{equation}
 \label{eqn:DeltaMVT}
 \D x^r = d^{2D} r(r-1)...(r-(2D-1))\tilde{x}^{r-2D},
\end{equation}
where $\tilde{x}\in [x,x+2Dd]$.  In particular we notice that $\D(x^{2D+1})=O(x)$, hence
$$
\D\left(\frac{2}{(2D+1)!}x^{2D+1}\right)=ax+b
$$
for some $a,b\in\R$.  Computations show that
\begin{equation}
\label{eqn:Deltax3}
a=2\sum_{j=0}^{2D}(-1)^j\rez{j!(2D-j)!}((2D-j)d)^{2D}=2d^{2D}.
\end{equation}

By definition $\psi_0(x) =\psi_{2D}^{(2D)}(s)$ so from (\ref{eqn:DeltaId}) we have
\begin{equation}
 \label{eqn:DeltaIdPsi2}
 \D\psi_{2D}(x) = \int_x^{x+d} \int_{t_{2D}}^{t_{2D}+d}\cdots\int_{t_2}^{t_2+d} \psi_0(t_1)\ dt_1...dt_{2D}.
\end{equation}
By Proposition \ref{thm:ECharPos}, the Euler characteristics $\chi_1(\Ga_{\ga})$ are all positive.  Hence $\psi_0(x)$ is non-decreasing and it follows from (\ref{eqn:DeltaIdPsi2}) that
\begin{equation}
\label{eqn:PsiIneq}
\psi_0(x)\leq d^{-2D}\D\psi_{2D}(x)\leq\psi_0(x+2Dd).
\end{equation}
It also follows from (\ref{eqn:DeltaMVT}) and (\ref{eqn:Deltax3}) that
$$
d^{-2D}\D\left(\frac{2}{(2D+1)!}x^{2D+1} + \sum_{\al\in S_{2D}^{\R}} c_{2D}(\al)x^{2D+\al}\right) = 2x + O(x^{3/4}).
$$
Thus it remains to show that for $p\in\{ -2,-1,0,1,2\}$
$$
d^{-2D}\D\left(\sum_{\al\in S_{2D}^{p/4}} c_{2D}(\al)x^{2D+\al}\right) = O(x^{3/4})
$$
and the proposition follows by (\ref{eqn:PsiIneq}).

Let $p\in\{ -2,-1,0,1,2\}$.  In order to estimate
$$
d^{-2D}\D\left(\sum_{\al\in S_{2D}^{p/4}} c_{2D}(\al)x^{2D+\al}\right)
$$
we need two estimates for $\D\left(c_{2D}(\al)x^{2D+\al}\right)$, where $\al\in S_{2D}^{p/4}$.  The residues at the poles of $R_{\Ga,1}'(s)/R_{\Ga,1}(s)$ are $O(1)$, so for $\al\in S_{2D}^{p/4}$
$$
d^{-2D}\D\left(c_{2D}(\al)x^{2D+\al}\right) = O\left(d^{-2D}|\al|^{-(2D+1)}x^{2D+p/4}\right).
$$
On the other hand, it follows from (\ref{eqn:DeltaMVT}) that
$$
d^{-2D}\D\left(c_{2D}(\al)x^{2D+\al}\right) = O\left(|\al|^{-1}x^{p/4}\right).
$$
Define $N_p(t)$ to be the number of poles of $R_{\Ga,1}'(s)/R_{\Ga,1}(s)$ on the interval $\frac{p}{4} + i(0,t]$.  From Lemma \ref{lem:NEstimate} we have that $N_p(t)=O\left( t^D\right)$.  Thus
\begin{eqnarray}
\label{eqn:NonRealPoles}
\lefteqn{d^{-2D}\D\left(\sum_{\al\in S_{2D}^{p/4}} c_{2D}(\al)x^{2D+\al}\right)} \nonumber \\
  & = & O\left(x^{p/4}\int_1^{K^D} t^{-1}dN(t) + d^{-2D}x^{2D+p/4}\int_{K^D}^{\infty}t^{-(2D+1)}dN(t)\right) \nonumber \\
  & = & O\left(K^{D-1}x^{p/4} + K^{-(D+2)}d^{-2D}x^{2+p/4}\right)
\end{eqnarray}
If we choose $K$ and $d$ appropriately, then we can conclude from (\ref{eqn:NonRealPoles}) that $d^{-2}\D\left(\sum_{\al\in S_2^{p/4}} c_2(\al)x^{2+\al}\right) = O(x^{3/4})$, as required.
\qed

\begin{proposition}
\label{pro:PsiTildeEst}
$\tilde{\psi}(x)=8x+O\left(x^{3/4}\right)$
\end{proposition}
\prf
The proof follows exactly as for the previous proposition, replacing $R_{\Ga,1}(s)$ with $R_{\Ga,\tilde{\si}}(s)$ throughout.  It follows in particular from the fact (Theorem \ref{thm:RuelleMain}) that $R_{\Ga,\tilde{\si}}(s)$ has a zero of order eight at the point $s=1$ and from the fact (Lemma \ref{lem:SigmaTilde}) that $\tr\tilde{\si}(b_{\ga})\geq 0$ for all $\ga\in\CE_P^p(\Ga)$ and $\tr\tilde{\si}(b_{\ga})=0$ if and only if $\ga$ is non-regular.
\qed

\subsection{The Wiener-Ikehara Theorem}

We shall use the following version of the Wiener-Ikehara theorem (see also \cite{Chandrasekharan68}, Chapter XI, Theorem 2, and \cite{Deitmar04}, Theorem 3.2).

Let $R_k(s)$, $k\in\N$ be a sequence of rational functions on $\C$.  for an open set $U\subset\C$ let $\N(U)$ be the set of natural numbers $k$ such that the pole divisor of $R_k$ does not intersect $U$.  We say that the series
$$
\sum_{k\in\N}R_k(s)
$$
\emph{converges weakly locally uniformly on $\C$} if for every open $U\subset\C$ the series
$$
\sum_{k\in\N(U)}R_k(s)
$$
converges locally uniformly on $U$.

\begin{theorem}\textnormal{(Wiener-Ikehara)}
\label{thm:WienerIkehara}
\index{Wiener-Ikehara Theorem}
Let $A(x)\geq 0$ be a monotonic measurable function on $\R_+$.  Suppose that the integral
$$
L(s)=\int_0^{\infty}A(x)e^{-sx}\ dx
$$
converges for $s\in\C$ with $\Re(s)>1$.  Suppose further that there are $j\in\N$, $r\in\R$ and a countable set $I$, and for each $i\in I$ there is $\th_i\in\C$ with $\Re(\th_i)<1$ and $c_i\in\Z$, such that the function
$$
L(s)-\left(\frac{\partial\ }{\partial s}\right)^{j+1}\frac{r}{s-1}-\sum_{i\in I}c_i\left(\frac{\partial\ }{\partial s}\right)^{j+1}\rez{s-\th_i}
$$
extends to a holomorphic function on the half-plane $\Re(s)\geq 1$.  Here we assume the sum converges weakly locally uniformly absolutely on $\C$.  Then
$$
\lim_{x\ra\infty}A(x)x^{-(j+1)}e^{-x}=r.
$$
\end{theorem}
\prf
This is a straightforward variation of the proof of the Wiener-Ikehara-Theorem in \cite{Chandrasekharan68}.
\qed

\subsection{The Dirichlet series}
\label{sec:DirSeries}
Let $B^0$ \index{$B^0$} be a closed subset of the compact group $B$ with the following properties: it is of measure zero; it is invariant under the map $b\mapsto b^{-1}$ and contains all fixed points of this map; and its complement $B^1=B\smallsetminus B^0$ \index{$B^1$} in $B$ is homeomorphic to an open subset of Euclidean space each of whose connected components is contractible.

The Weyl group $W=W(M,B)$ contains two elements and the non-trivial element acts on $B$ by $b\mapsto b^{-1}$, so the invariance condition above says that both $B^0$ and $B^1$ are invariant under the action of $W$.  The fixed points under the action of the Weyl group are precisely the non-regular elements of $B$, so the assumption that $B^0$ contains all these fixed points is equivalent to $B^{\nreg}\subset B^0$, where $B^{\nreg}$ \index{$B^{\nreg}$} denotes the subset of nonregular elements of $B$.  The action of $W$ on $B^1$ permutes the connected components and the assumption $B^{\nreg}\subset B^0$ implies that the quotient space $B^1/W$ is also homeomorphic to an open subset of Euclidean space each of whose connected components is simply connected.

For subsets $S$ and $T$ of $G$ we denote by $S.T$ the subset
$$
S.T=\{sts^{-1}:s\in S,t\in T\}
$$
of $G$.  Let $M_{\ell}^1=M.B^1\subset M$.  Let $\CE_P^0(\Ga)$ \index{$\CE_P^0(\Ga)$} and $\CE_P^1(\Ga)$ \index{$\CE_P^1(\Ga)$} be the subsets of $\CE_P(\Ga)$ consisting of all conjugacy classes $[\ga]$ such that $b_{\ga}\in B^0$ or $b_{\ga}\in B^1$ respectively.  The assumption that $B^0$ is of measure zero is not required for the following lemma, but will be needed later on.

\begin{lemma}
There exist a set $M_c\subset M_{\ell}^1$ with compact closure in $M$ and a monotonically increasing sequence $(g_n)$ \index{$g_n$} of smooth functions on $M$, supported on $M_c$, which are invariant under conjugation by elements of $K_M$, and whose orbital integrals satisfy
$$
\CO_{b_{\ga}}^M(g_n)=\int_{M/B}g_n(xb_{\ga}x^{-1})\ dx\ra 1\ \ \ \textrm{ as }n\ra\infty
$$
for all $\ga\in\CE_P^1(\Ga)$.
\end{lemma}
\prf
We view $M_{\ell}^1$ as a fibre bundle with base space $B^1/W$ and fibres homeomorphic to $M/B$.

Let $d(\cdot,\cdot)$ denote the metric on $B$ given by the form $b$ in (\ref{eqn:norm}).  For $n\in\N$ let $\tilde{B}_n\subset B^1$ be the set $\tilde{B}_n=\{b\in B:d(b,B^0)\geq 1/n\}$ and let $B_n=\tilde{B}_n/W$.  Then, by \cite{Warner83}, Corollary 1.11, for each $n\in\N$ there exists a function $h_n:B/W\ra\R$ such that $0\leq h_n(b)\leq 1$ for all $b\in B/W$, for all $b\in B_n$ we have $h_n(b)=1$ and $h_n$ is supported on $B^1/W$.  We may assume that the $h_n$'s form an increasing series.

Let $U$ be a compact neighbourhood of a point in $M/B$ homeomorphic to a subset of Euclidean space.  Let $k:M/B\ra\R$ be a smooth positive function supported on $U$ and satisfying $\int_{M/B}k(m)dm=1$.

We now define, for $n\in\N$, functions $\tilde{g}_n:M\ra\R$ as follows.  On each connected component $V$ of $B^1/W$ we fix a trivialisation of the restriction to $V$ of the bundle $M/B\ra M_{\ell}^1\ra B^1/W$.  Then for $v\in V$ and $m\in M/B$ define $\tilde{g}_n(mvm^{-1})=h_n(v)k(m)$.  Since $V$ is simply connected there are no problems with global agreement of this definition.  For $m\in M\smallsetminus M_{\ell}^1$ define $\tilde{g}_n(m)=0$.  Finally we define the functions $g_n:M\ra\R$ by
$$
g_n(m)=\rez{2}\int_{K_M}\tilde{g}_n(kmk^{-1})\,dk
$$
for all $m\in M$.  Then the functions $g_n$ form an increasing sequence of smooth $K_M$-invariant functions supported on the compact set $M_c$, which we define to be the closure in $M$ of $K_M.(U.B^1)$.  We will show that their orbital integrals have the required properties.  Let $b\in B^1$.
\begin{eqnarray*}
\CO_{b}^M(g_n) & = & \int_{M/B}g_n(mbm^{-1})\ dm \\
  & = & \rez{2}\int_{M/B}\int_{K_M}\tilde{g}_n(kmbm^{-1}k^{-1})\ dk\,dm \\
  & = & \rez{2}\int_{M/B}\int_{K_M}\tilde{g}_n(mkbk^{-1}m^{-1})\ dk\,dm \\
  & = & \rez{2}\int_{M/B}\int_{B}\tilde{g}_n(mbm^{-1}) + \tilde{g}_n(mb^{-1}m^{-1})\ dk\,dm \\
  & = & \rez{2}\left(h_n(b)+h_n(b^{-1})\right)\int_{M/B}k(m)\ dm \\
  & = & h_n(b).
\end{eqnarray*}
The last equality holds since $bW=b^{-1}W$ in $B^1/W$.  There exists $N\in\N$ such that $b\in B_n$ for all $n\geq N$, hence $\CO_{b}^M(g_n)\ra 1$ as $n\ra\infty$, by the definition of the functions $h_n$.
\qed

For $n\in\N$ let
$$
L^j_n(s)=\sum_{[\ga]\in\CE_P^1(\Ga)}\chi_1(\Ga_{\ga_0})\CO_{b_{\ga}}(g_n)l_{\ga_0}\frac{l_{\ga}^{j+1}e^{-sl_{\ga}}}{\det(1-(a_{\ga}b_{\ga})^{-1}|\n)},\index{$L^j_n(s)$}
$$
let
$$
r_n=\int_M g_n(x)dx,\index{$r_n$}
$$
and let
$$
m_{n,\la}=\sum_{q=0}^4 (-1)^q r_n\dim H^q(\n,\pi_K)^{\la}.\index{$m_{n,\la}$}
$$

\begin{proposition}
\label{pro:DirSeries}
For all $n\in\N$ and for $j\in\N$ large enough the series $L^j_n(s)$ converges locally uniformly in the set $\{s\in\C:\Re(s)>1\}$.

The function $L^j_n(s)$ can be written as a Mittag-Leffler series
$$
L_n^j(s) = \left(\frac{\partial\ }{\partial s}\right)^{j+1}\frac{r_n}{s-1} + \sum_{\pi\in\hat{G}}N_{\Ga}(\pi)\sum_{\la\in\a^*}m_{n,\la}\left(\frac{\partial\ }{\partial s}\right)^{j+1}\rez{s-\la(H_1)},
$$
where the summand of the double series corresponding to $\pi=triv$, $\la=-2\rho_P$ is excluded.  The double series converges weakly locally uniformly on $\C$.  For $\pi\in\hat{G}$, $\la\in\a^*$ such that $(\pi,\la)\neq (triv,-2\rho_P)$ we have $m_{n,\la}\neq 0$ only if $\Re(\la)>-2\rho_P$.  Thus, in particular, the double series converges locally uniformly on $\{s\in\C:\Re(s)>1\}$.
\end{proposition}
\prf
We can put $f=g_n$ in (\ref{eqn:TestFunction}) to define the function $\Phi_{g_n,j,s}$, which by Proposition \ref{pro:STF2}, for $j$ and $\Re(s)$ large enough, goes into the Selberg trace formula to give the equation
$$
\sum_{\pi\in\hat{G}} N_{\Ga}(\pi)\tr\pi(\Phi_{g_n})=\sum_{[\ga]\in\CE_P(\Ga)}\vol(\Ga_{\ga}\bs G_{\ga})\CO_{b_{\ga}}(g_n)\frac{l_{\ga}^{j+1}e^{-sl_{\ga}}}{\det(1-(a_{\ga}b_{\ga})^{-1}|\n)}.
$$
By the definition of the functions $g_n$, the orbital integrals $\CO_{b_{\ga}}(g_n)$ are equal to zero when $\ga\in\CE_P^0(\Ga)$.  Also, for $\ga\in\CE_P^1(\Ga)$, since we have assumed that $B^{\nreg}\subset B^0$ we have that $\ga$ is regular and we get from Lemma \ref{lem:ECharEqn} that
$$
\vol(\Ga_{\ga}\bs G_{\ga})=\chi_1(\Ga_{\ga})l_{\ga_0}.
$$
Thus we can see that the geometric side of the trace formula is equal to $L_n^j(s)$.  It follows that for $j$ and $\Re(s)$ large enough the Dirichlet series $L_n^j(s)$ converges absolutely.

\begin{lemma}
$$
\tr\pi(\Phi_{g_n})=(-1)^{j+1}\left(\frac{\partial\ }{\partial s}\right)^{j+1}r_n\sum_{\la\in\a^*}\dim H^{\bullet}(\n,\pi_K)^{\la}\rez{s-\la(H_1)}.
$$
\end{lemma}
\prf
Replacing $f_{\si}$ with $g_n$ in (\ref{eqn:TracePiPhi1}) we get
\begin{eqnarray*}
\tr\pi(\Phi_{g_n}) & = & \int_{MA^-}g_n(m)\Th^{MA}_{H^{\bullet}(\n,\pi_K)}(ma)dm\,g_s^j(a)da \\
  & = & \int_{A^-}\int_M g_n(m)dm\ \tr\!\left(a|H^{\bullet}(\n,\pi_K)\right)l_a^{j+1}e^{-sl_a}da. \\
  & = & r_n\int_0^{\infty}\sum_{\la\in\a^*}\dim H^{\bullet}(\n,\pi_K)^{\la}e^{(\la(H_1)-s)t}\,t^{j+1}\,dt \\
  & = & (-1)^{j+1}\left(\frac{\partial\ }{\partial s}\right)^{j+1}r_n\sum_{\la\in\a^*}\dim H^{\bullet}(\n,\pi_K)^{\la}\rez{s-\la(H_1)}.
\end{eqnarray*}
\qed

By the above lemma and Proposition \ref{pro:AAction}, and using similar arguments to those used in the proof of Proposition \ref{pro:TrivAAction}, we can see that the spectral side of the trace formula is equal to the Mittag-Leffler series given in the propostion.

Since $L_n^j(s)$ is a Dirichlet series with positive coefficients it will converge locally uniformly for $s$ in some open set.  By proving the convergence of the Mittag-Leffler series we shall show that the Dirichlet series extends to a holomorphic function on $\{\Re(s)>1\}$ and hence, since it has positive coefficients, converges locally uniformly there.

We recall that for $\la\in\a^*$ we write $\|\la\|$ for the norm given by the form $b$ in (\ref{eqn:norm}).  It then follows from Proposition \ref{pro:Bounds} that there exist $m_1\in\N$ and $C>0$ such that for every $\pi\in\hat{G}$ and every $\la\in\a^*$ we have
$$
|m_{n,\la}|\leq C(1+\|\la\|)^{m_1}.
$$
If $S$ denotes the set of all pairs $(\pi,\la)\in\hat{G}\x\a^*$ such that $m_{n,\la}\neq 0$, then there exists $m_2\in\N$ such that
$$
\sum_{(\pi,\la)\in S}\frac{N_{\Ga}(\pi)}{(1+\|\la\|)^{m_2}}<\infty.
$$

Now, let $U\subset\C$ be open and let $S(U)$ be the set of all pairs $(\pi,\la)\in\hat{G}\x\a^*$ such that $m_{n,\la}\neq 0$ and $\la(H_1)\notin U$. Let $V$ be a compact subset of $U$.  We have to show that for some $j\in\N$, which does not depend on $U$ or $V$,
$$
\sup_{s\in V}\sum_{(\pi,\la)\in S(U)}\left|\frac{N_{\Ga}(\pi)m_{n,\la}}{(s-\la(H_1))^{j+2}}\right|<\infty.
$$
Let $m_1$, $m_2$ be as above and let $j\geq m_1+m_2-2$.  Since $V\subset U$ and $V$ is compact there is $\ep>0$ such that $s\in V$ and $(\pi,\la)\in S(U)$ implies $|s-\la(H_1)|\geq\ep$.  Hence there is $c>0$ such that for every $s\in V$ and every $(\pi,\la)\in S(U)$,
$$
|(s-\la(H_1)|\geq c(1+\|\la\|).
$$
Putting this all together we get
\begin{eqnarray*}
\sup_{s\in V}\sum_{(\pi,\la)\in S(U)}\left|\frac{N_{\Ga}(\pi)m_{n,\la}}{(s-\la(H_1))^{j+2}}\right| & \leq & \sup_{s\in V}\sum_{(\pi,\la)\in S(U)} \rez{c^{j+2}}\frac{N_{\Ga}(\pi)|m_{n,\la}|}{(1+\|\la\|)^{j+2}} \\
  & \leq & \sup_{s\in V}\sum_{(\pi,\la)\in S(U)}\frac{C}{c^{j+2}}\frac{N_{\Ga}(\pi)}{(1+\|\la\|)^{m_2}} \\
  & < & \infty,
\end{eqnarray*}
which proves the proposition.
\qed

\subsection{Estimating $\psi_n(x)$}

Let
$$
\phi^j_n(x)=\sum_{{[\ga]\in\CE_P^1(\Ga)}\atop{N(\ga)\leq x}} \chi_1(\Ga_{\ga_0})\CO_{b_{\ga}}(g_n)l_{\ga_0}\frac{l_{\ga}^{j+1}}{\det(1-(a_{\ga}b_{\ga})^{-1}|\n)}.\index{$\phi^j_n(x)$}
$$

\begin{lemma}
\label{lem:Abel}
$$
\int_0^{\infty}\phi^j_n(e^x)e^{-sx}\ dx=L^j_n(s)
$$
\end{lemma}
\prf
This follows from Abel's summation formula (\cite{Chandrasekharan68}, Chapter VII, Theorem 6).
\qed

Let
$$
\phi_n(x)=\sum_{{[\ga]\in\CE_P^1(\Ga)}\atop{N(\ga)\leq x}} \chi_1(\Ga_{\ga_0})\CO_{b_{\ga}}(g_n)l_{\ga_0}\rez{\det(1-(a_{\ga}b_{\ga})^{-1}|\n)}.\index{$\phi_n(x)$}
$$

\begin{lemma}
\label{lem:Phin}
For each $n\in\N$ we have
$$
\phi_n(x)\sim r_n x.
$$
\end{lemma}
\prf
By Propostion \ref{pro:DirSeries} and Lemma \ref{lem:Abel} we can apply Theorem \ref{thm:WienerIkehara} to the series $L_n^j(s)$ and the function $\phi_n^j(x)$ to deduce that
\begin{equation}
\label{eqn:PhijCgence}
\lim_{x\ra\infty}\frac{\phi^j_n(x)}{x(\log x)^{j+1}}=r_n.
\end{equation}
Also it is clear that
$$
\phi_n(x)\geq\frac{\phi^j_n(x)}{(\log x)^{j+1}},
$$
so it follows that
$$
\liminf_{x\ra\infty}\frac{\phi_n(x)}{x}\geq r_n.
$$
Let $0<\mu<1$.  Then
\begin{eqnarray*}
\phi^j_n(x) & \geq & \sum_{{[\ga]\in\CE_P^1(\Ga)}\atop{x^{\mu}<N(\ga)\leq x}} \chi_1(\Ga_{\ga_0})\CO_{b_{\ga}}(g_n)l_{\ga_0}\frac{l_{\ga}^{j+1}}{\det(1-(a_{\ga}b_{\ga})^{-1}|\n)} \\
  & \geq & \mu^{j+1}(\log x)^{j+1}\sum_{{[\ga]\in\CE_P^1(\Ga)}\atop{x^{\mu}<N(\ga)\leq x}} \chi_1(\Ga_{\ga_0})\CO_{b_{\ga}}(g_n)l_{\ga_0}\rez{\det(1-(a_{\ga}b_{\ga})^{-1}|\n)} \\
  & = & \mu^{j+1}(\log x)^{j+1}(\phi_n(x)-\phi_n(x^{\mu})).
\end{eqnarray*}
From this we get
\begin{equation}
\label{eqn:PhinBound}
\frac{\phi_n(x)}{x}\leq\mu^{-(j+1)}\frac{\phi^j_n(x)}{x(\log x)^{j+1}}+\frac{\phi_n(x^{\mu})}{x^{\mu}x^{1-\mu}}.
\end{equation}

Assume first that $\phi_n(x)/x$ tends to infinity as $x\ra\infty$.  Then there is a sequence $x_m$ of positive real numbers, tending to infinity such that $\phi(x_m)/x_m$ tends to infinity as $m\ra\infty$ and
$$
\frac{\phi_n(x_m)}{x_m}\geq\frac{\phi_n(x_m^{\mu})}{x_m^{\mu}}
$$
for all $m$.  Then
$$
\frac{\phi_n(x_m)}{x_m}\leq\mu^{-(j+1)}\frac{\phi_n^j(x_m)}{x_m(\log x_m)^{j+1}}+\frac{\phi_n(x_m)}{x_m}.\rez{x_m^{1-\mu}},
$$
so that
$$
\frac{\phi_n(x_m)}{x_m}\leq\frac{\mu^{-(j+1)}\frac{\phi_n^j(x_m)}{x_m(\log x_m)^{j+1}}}{1-\rez{(x_m)^{1-\mu}}}.
$$
By (\ref{eqn:PhijCgence}) the right hand side converges, so we have a contradiction.  This implies that
$
\lim_{x\ra\infty}\frac{\phi_n(x)}{x}
$
exists and is finite and we set
$$
L=\limsup_{x\ra\infty}\frac{\phi_n(x)}{x}=\limsup_{x\ra\infty}\frac{\phi_n(x^{\mu})}{x^{\mu}}.
$$
From (\ref{eqn:PhijCgence}) and (\ref{eqn:PhinBound}) we get
$$
L \ \leq \ r_n\mu^{-(j+1)}+L\limsup_{x\ra\infty}\rez{x^{1-\mu}}
  \ = \ r_n\mu^{-(j+1)}.
$$
Since this holds for any value of $\mu$ in the interval $0<\mu<1$ we get that $L\leq r$ and the lemma follows.
\qed

For $n\in\N$ and $x>0$ let
$$
\psi_n(x)=\sum_{{[\ga]\in\CE_P^1(\Ga)}\atop{N(\ga)\leq x}} \chi_1(\Ga_{\ga_0})\CO_{b_{\ga}}(g_n)l_{\ga_0}.\index{$\psi_n(x)$}
$$

\begin{proposition}
\label{pro:PsiNEst}
For each $n\in\N$ we have
$
\psi_n(x)\sim r_n x.
$
\end{proposition}
\prf
We note that $0<\det(1-a_{\ga}b_{\ga}|\n)<1$ for all $\ga\in\CE_P^{\reg}(\Ga)$ and that the value of the determinant tends to $1$ as $l_{\ga}$ tends to infinity, hence for $0<\ep<1$ there are only finitely many $\ga\in\CE_P^{\reg}(\Ga)$ such that $\det(1-a_{\ga}b_{\ga}|\n)\leq 1-\ep$.  Since $\CE_P^1(\Ga)\subset\CE_P^{\reg}(\Ga)$ the same holds for $\ga\in\CE_P^1(\Ga)$.  We fix $0<\ep<1$ and define for each $n\in\N$ the functions $\phi_{n,\ep}$ and $\psi_{n,\ep}$ to be the same sums as for $\phi_n$ and $\psi_n$ respectively but restricted to those classes $[\ga]\in\CE_P^1(\Ga)$ such that $1-\ep<\det(1-a_{\ga}b_{\ga}|\n)<1$.  It then follows that both
\begin{equation}
\label{eqn:PhiEp}
\frac{\phi_n(x)-\phi_{n,\ep}(x)}{x}\ra 0
\end{equation}
and
\begin{equation}
\label{eqn:PsiEp}
\frac{\psi_n(x)-\psi_{n,\ep}(x)}{x}\ra 0
\end{equation}
as $x\ra\infty$.  Now
$$
1-\ep<\det(1-a_{\ga}b_{\ga}|\n)<1
$$
immediately implies
$$
\frac{1-\ep}{\det(1-a_{\ga}b_{\ga}|\n)}<1<\rez{\det(1-a_{\ga}b_{\ga}|\n)},
$$
so summing up we get
$$
\frac{\phi_{n,\ep}(x)}{x}(1-\ep)<\frac{\psi_{n,\ep}(x)}{x}<\frac{\phi_{n,\ep}(x)}{x}.
$$
By Lemma \ref{lem:Phin} and (\ref{eqn:PhiEp}), it then follows that
$$
r_n(1-\ep)\leq\liminf_{x\ra\infty}\frac{\psi_{n,\ep}(x)}{x}\leq\limsup_{x\ra\infty}\frac{\psi_{n,\ep}(x)}{x}\leq r_n.
$$
It then follows from (\ref{eqn:PsiEp}) that
$$
r_n(1-\ep)\leq\liminf_{x\ra\infty}\frac{\psi_n(x)}{x}\leq\limsup_{x\ra\infty}\frac{\psi_n(x)}{x}\leq r_n.
$$
Since this holds for any value of $\ep$ in the interval $0<\ep<1$ the proposition is proven.
\qed

\subsection{Estimating $\psi^1(x)$}

\begin{lemma}
\label{lem:RnLimit}
The sequence $(r_n)$ is monotonically increasing and $r_n\ra 2$ as $n\ra\infty$.
\end{lemma}
\prf
That the sequence $(r_n)$ is monotonically increasing is clear since the sequence of functions $(g_n)$ is monotonically increasing.

By Weyl's integral formula (\cite{Knapp}, Proposition 5.27), and since the functions $g_n$ are only non-zero on elliptic elements of $M$, we have
\begin{eqnarray*}
r_n & = & \int_M g_n(x)\ dx \\
  & = & \rez{|W(B:G)|}\int_B\int_{M/B}g_n(xbx^{-1})|D(b)|^2\ dx\,db \\
  & = & \rez{|W(B:G)|}\int_B\CO^M_b(g_n)|D(b)|^2\ db,
\end{eqnarray*}
where $W(B:G)$ is the Weyl group and $D(b)$ is the Weyl denominator.  Since the sequence of functions $g_n$ is monotonically increasing and the functions are all supported within a given compact subset of $M$ we can interchange integral and limit to get
\begin{eqnarray*}
\lim_{n\ra\infty}r_n & = & \lim_{n\ra\infty}\rez{|W(B:G)|}\int_B\CO^M_b(g_n)|D(b)|^2\ db \\
  & = & \rez{|W(B:G)|}\int_B\lim_{n\ra\infty}\CO^M_b(g_n)|D(b)|^2\ db.
\end{eqnarray*}
Furthermore, the orbital integrals $\CO^M_b(g_n)$ tend to one as $n\ra\infty$, except for $b$ in a set of measure zero, so it follows that
$$
\lim_{n\ra\infty}r_n=\rez{|W(B:G)|}\int_B|D(b)|^2\ db.
$$

The Weyl group $W(B:G)$ consists of the identity element and the element $\diag(1,-1,1,-1)$.
We can also compute the Weyl denominator $D(b)$ for $b\in B$.  Let
$$
R(\th)=\matrix{\cos\th}{-\sin\th}{\sin\th}{\cos\th}\in\SO(2)
$$
and let
$$
R(\th,\phi)=\matrixtwo{R(\th)}{R(\phi)}\in B.
$$
Then
\begin{eqnarray*}
D(R(\th,\phi)) & = & e^{i(\th+\phi)}(1-e^{-2i\th})(1-e^{-2i\phi}) \\
  & = & (e^{i(\th+\phi)}+e^{-i(\th+\phi)})-(e^{i(\th-\phi)}+e^{-i(\th-\phi)}) \\
  & = & 2\cos(\th+\phi)-2\cos(\th-\phi) \\
  & = & 2(\cos\th\cos\phi-\sin\th\sin\phi-\cos\th\cos\phi-\sin\th\sin\phi) \\
  & = & -4\,\sin\th\,\sin\phi
\end{eqnarray*}
We have normalised the Haar measure on $B$ so that $\int_B db=1$, hence
\begin{eqnarray*}
\rez{|W(B:G)|}\int_B|D(b)|^2\ db & = & \rez{2}\int_0^{2\pi}\int_0^{2\pi}16\,\sin\!^2\th\,\sin\!^2\phi\ \frac{d\th\,d\phi}{4\pi^2} \\
  & = & \frac{2}{\pi^2}\int_0^{2\pi}\int_0^{2\pi}\sin\!^2\th\,\sin\!^2\phi\ d\th\,d\phi \\
  & = & \frac{2}{\pi^2}.\pi^2 \= 2.
\end{eqnarray*}
This proves the lemma.
\qed

Let
$$
\psi^1(x)=\sum_{{[\ga]\in\CE_P^1(\Ga)}\atop{N(\ga)\leq x}}\chi_1(\Ga_{\ga_0})l_{\ga_0}.\index{$\psi(x)$}
$$
\begin{proposition}
\label{pro:PsiEst}
$\psi^1(x)\sim\frac{2x}{\log x}$.
\end{proposition}
\prf
For all $n\in\N$ and all $x>0$ we have
$
\psi_n(x)\leq\psi^1(x)\leq\psi(x).
$
It then follows, using Proposition \ref{pro:PsiHatEst} and Proposition \ref{pro:PsiNEst}, that for all $n\in\N$
$$
r_n=\liminf_{x\ra\infty}\frac{\psi_n(x)}{x}\leq\liminf_{x\ra\infty}\frac{\psi^1(x)}{x}\leq\limsup_{x\ra\infty}\frac{\psi^1(x)}{x}\leq\limsup_{x\ra\infty}\frac{\psi(x)}{x}=2.
$$
We can then deduce from Lemma \ref{lem:RnLimit} that
$
\lim_{x\ra\infty}\frac{\psi^1(x)}{x}=2.
$
\qed

\subsection{Estimating $\pi(x)$, $\tilde{\pi}(x)$ and $\pi^1(x)$}

To keep the notation less cluttered, in the sums over conjugacy classes in $\Ga$ which appear in this section we shall not specify which set of classes the sum is being taken over, since this will be clear from the context.  We shall always use $\ga_0$ to denote primitive elements and where $\ga$ and $\ga_0$ appear together in the same formula we shall mean that $\ga_0$ is the primitive element underlying $\ga$.

\begin{proposition}
\label{pro:PiEst}
$$
\displaystyle\lim_{x\ra\infty}\frac{\pi(x)}{2x/\log x} = \lim_{x\ra\infty}\frac{\psi(x)}{2x} = 1.
$$
\end{proposition}
\prf
We can write
$$
\psi(x)=\sum_{N(\ga_0)\leq x} n_{\ga_0}\chi_1(\Ga_{\ga_0})l_{\ga_0},
$$
where $n_{\ga_0}\in\N$ is maximal such that $N(\ga_0)^{n_{\ga_0}}\leq x$.  By definition $N(\ga_0)=e^{l_{\ga_0}}$, so $N(\ga_0)^n\leq x$ implies $nl_{\ga_0}\leq \log x$ and we can see that
\begin{equation}
\label{eqn:PsiPiIneq}
\psi(x)\leq\log(x)\pi(x).
\end{equation}

Next we fix a real number $0<a<1$.  By Proposition \ref{thm:ECharPos} the Euler characteristic $\chi_1(\Ga_{\ga_0})>0$ for all $\ga_0\in\CE_P^p(\Ga)$, so for $x>1$,
$$
\psi(x) \geq \sum_{x^a <N(\ga_0)\leq x} \chi_1(\Ga_{\ga_0})l_{\ga_0}.
$$
As above, $N(\ga_0)>x^a$ implies $l_{\ga_0}>\log x^a$, hence
$$
\psi(x)\geq a\log x\sum_{x^a <N(\ga_0)\leq x}\chi_1(\Ga_{\ga_0}) = a\log x(\pi(x)-\pi(x^a)).
$$
Since $\Ga\subset G$ is discrete, $\pi(x^a)<Cx^a$ for some constant $C$ so
$$
\psi(x)>a\pi(x)\log x - aCx^a\log x,
$$
which gives
$$
\frac{\psi(x)}{2x}>a\pi(x)\frac{\log x}{2x} - aC\frac{\log x}{2x^{1-a}}.
$$
Since $0<a<1$, it follows that $(\log x)/x^{1-a}\ra 0$ as $x\ra\infty$ and
$$
\lim_{x\ra\infty}\frac{\psi(x)}{2x} \geq a\lim_{x\ra\infty}\frac{\pi(x)}{2x/\log x}
$$
for all $0<a<1$.  Hence
$$
\lim_{x\ra\infty}\frac{\psi(x)}{2x} \geq \lim_{x\ra\infty}\frac{\pi(x)}{2x/\log x}.
$$
Together with (\ref{eqn:PsiPiIneq}) and Proposition \ref{pro:PsiHatEst} this proves the proposition.
\qed

\begin{proposition}
\label{pro:PiTildeEst}
$$
\displaystyle\lim_{x\ra\infty}\frac{\tilde{\pi}(x)}{8x/\log x} = \lim_{x\ra\infty}\frac{\tilde{\psi}(x)}{8x} = 1.
$$
\end{proposition}
\prf
Exactly as for the previous proposition, making use of Proposition \ref{pro:PsiTildeEst} and Lemma \ref{lem:SigmaTilde}.
\qed

\begin{proposition}
\label{pro:PiEst2}
$$
\pi(x)=2\li(x)+O\left(\frac{x^{3/4}}{\log x}\right).
$$
\end{proposition}
\prf
We consider the function
$$
S(x)  \=  \sum_{N(\ga)\leq x}\chi_1(\Ga_{\ga_0})\frac{l_{\ga_0}}{l_{\ga}}
  \=  \sum_{N(\ga_0)\leq x} \chi_1(\Ga_{\ga_0}) + \sum_{k\geq 2}\sum_{N(\ga_0)\leq x^{1/k}} \chi_1(\Ga_{\ga_0})\rez{k}.
$$
We consider the double sum on the right.  Since $\Ga\subset G$ is discrete there is a geodesic $\ga_{\rm min}$ of minimum length.  For a given $x$ the inner sum contains at least one summand only for $k\leq \log x/l_{\ga_{\rm min}}$.  For each such $k\geq 2$, by Proposition \ref{pro:PiEst}, the inner sum is equal to $O(\sqrt{x}/\log x)$.  Therefore we have
\begin{equation}
\label{eqn:SPiEst}
S(x) = \pi(x)+O\left(\sqrt{x}\right).
\end{equation}
Now
\begin{eqnarray*}
\int_2^x \frac{\psi(t)}{t\log^2\! t}\ dt & = & \int_2^x \sum_{N(\ga)\leq t} \chi_1(\Ga_{\ga_0})\frac{l_{\ga_0}}{t\log^2\! t}\ dt \\
  & = & \sum_{N(\ga)\leq x} \int_{N(\ga)}^x \chi_1(\Ga_{\ga_0})\frac{l_{\ga_0}}{t\log^2 t}\ dt \\
  & = & \sum_{N(\ga)\leq x} \chi_1(\Ga_{\ga_0})l_{\ga_0}\left(\rez{l_{\ga}}-\rez{\log x}\right) \\
  & = & S(x) - \frac{\psi(x)}{\log x}.
\end{eqnarray*}
Hence
$
S(x)=\int_2^x \frac{\psi(t)}{t\log^2\! t}\ dt + \frac{\psi(x)}{\log x}.
$
By Proposition~\ref{pro:PsiHatEst} we have
\begin{eqnarray*}
S(x) & = & \int_2^x \frac{2}{\log^2 t}\ dt + \frac{2x}{\log x} + O\left(\int_2^x \rez{t^{1/4}\log^2 t}\ dt\right) + O\left(\frac{x^{3/4}}{\log x}\right) \\
     & = & \left[-\frac{2t}{\log t}\right]_2^x + \int_2^x\frac{2}{\log t}\ dt + \frac{2x}{\log x} + O\left(\frac{x^{3/4}}{\log x}\right) \\
     & = & \int_2^x\frac{2}{\log t}\ dt + O\left(\frac{x^{3/4}}{\log x}\right).
\end{eqnarray*}
Together with (\ref{eqn:SPiEst}) this proves the proposition.
\qed

\begin{proposition}
\label{pro:PiTildeEst2}
$$
\tilde{\pi}(x)=8\li(x)+O\left(\frac{x^{3/4}}{\log x}\right).
$$
\end{proposition}

\prf
Exactly as for the previous proposition, making use of Proposition \ref{pro:PsiTildeEst} and Proposition \ref{pro:PiTildeEst}.
\qed

The proof of Theorem \ref{thm:PGT2} proceeds exactly as for the proof of Proposition \ref{pro:PiEst}, making use of Proposition \ref{pro:PsiEst}.

\end{document}